\documentclass{siamonline190516}
\usepackage[T1]{fontenc}
\usepackage[utf8]{inputenc}
\usepackage{lmodern}
\usepackage{verbatim}
\usepackage{amsfonts}
\usepackage{graphicx}
\usepackage{hyperref}
\usepackage{amssymb}
\usepackage{amsmath}
\usepackage{bussproofs}
\usepackage{pdfpages}
\usepackage{float}
\usepackage{cite}
\usepackage{tikz}
\usepackage{tikzscale}
\usetikzlibrary{automata,positioning,arrows,calc,decorations.shapes,decorations.markings,quotes,shapes.misc}
\usepackage{pgfplots}
\usepackage{makeidx}
\usepackage{mathrsfs}
\usepackage{dsfont}
\usepackage{svg}
\usepackage{array}
\usepackage{multirow}
\usepackage{bibentry}
\usepackage{comment}
\usepackage{faktor}
\usepackage[mathlines]{lineno}
\usepackage{subcaption}
\usepackage{algorithm}
\usepackage{algpseudocode}
\usepackage{float}
\usepackage[shortlabels]{enumitem}

\usepackage[normalem]{ulem}
\newcommand{\nt}[1]{{\color{black} #1}}

\definecolor{col}{rgb}{0,0,0}

\algdef{SE}[SUBALG]{Indent}{EndIndent}{}{\algorithmicend\ }%
\algtext*{Indent}
\algtext*{EndIndent}

\title{Random Multi-Type Spanning Forests for Synchronization on Sparse Graphs\footnote{This work is an extension of the conference proceedings~\cite{jaquard_smoothing_2023}.}}
\author{Hugo Jaquard, Pierre-Olivier Amblard, Simon Barthelmé, Nicolas Tremblay \\ \small CNRS, Univ. Grenoble Alpes, Grenoble INP, GIPSA-lab, France}

\theoremstyle{plain}
\newtheorem{thm}{Theorem}[section]
\newtheorem{prop}[thm]{Proposition}

\newtheorem{lem}[thm]{Lemma}

\newtheorem{rem}[thm]{Remark}
\newtheorem{cond}[thm]{Condition}

\renewcommand{\qedsymbol}{$\blacksquare$}

\newcommand\R{\mathbf{R}}

\newcommand\N{\mathbf{N}}
\newcommand\C{\mathbf{C}}
\newcommand\prb{\mathbf{P}}
\newcommand\E{\mathbf{E}}

\newcommand\eps{\varepsilon}

\newcommand\id{\mathrm{id}}

\newcommand{\argmin}{\operatornamewithlimits{argmin}}

\usetikzlibrary{arrows.meta}
\usetikzlibrary{backgrounds}
\usepgfplotslibrary{patchplots}
\usepgfplotslibrary{fillbetween}
\pgfplotsset{
     layers/standard/.define layer set={
         background,axis background,axis grid,axis ticks,axis lines,axis tick labels,pre main,main,axis descriptions,axis foreground
     }{
         grid style={/pgfplots/on layer=axis grid},
         tick style={/pgfplots/on layer=axis ticks},
         axis line style={/pgfplots/on layer=axis lines},
         label style={/pgfplots/on layer=axis descriptions},
         legend style={/pgfplots/on layer=axis descriptions},
         title style={/pgfplots/on layer=axis descriptions},
         colorbar style={/pgfplots/on layer=axis descriptions},
         ticklabel style={/pgfplots/on layer=axis tick labels},
         axis background@ style={/pgfplots/on layer=axis background},
         3d box foreground style={/pgfplots/on layer=axis foreground},
     },
 }
 
\usepackage{caption}
\begin{document}

    \maketitle

    \begin{abstract}
        Random diffusions are a popular tool in Monte-Carlo estimations, with well established algorithms such as Walk-on-Spheres (WoS) going back several decades. 
        In this work, we introduce diffusion estimators for the problems of angular synchronization and smoothing on graphs, in the presence of a rotation associated to each edge. Unlike classical WoS algorithms \nt{that are point-wise estimators, our diffusion} estimators allow for \emph{global} estimations by propagating along the branches of \nt{random spanning subgraphs} \nt{called} \emph{multi-type spanning forests}. \nt{Building upon efficient samplers based on variants of Wilson's algorithm}, we show that our estimators outperform standard numerical-linear-algebra solvers in challenging instances, depending on the topology and density of the graph.
    \end{abstract}
    


    
    \section*{Introduction}
    \label{sect:intro}

    Data processing over graphs is of interest in many fields of research, such as discrete geometry, signal processing or graph machine learning. A common setting is that of nodes supporting some kind of data, typically numerical values, whereas the edges describe the topology of the underlying space. In a number of situations, the edges can carry more precise geometric information regarding \emph{how} two data points defined on adjacent nodes should be compared, and thus processed. In geometrical parlance, this information is known as a \emph{connection} (an idea first formalized in~\cite{cartan1926groupes}).
    
    We focus here on a connection describing rotations of angle $\theta_{i,j}$ associated to each edge $(i,j)$, a setting which has found a growing number of applications over the last decade, including \emph{angular synchronization}~\cite{singer2011angular,yu_angular_2012}, signal processing over directed graphs~\cite{furutani_gsp_2019,zhang_magnet_2021}, discrete geometry processing~\cite{sharp2019vector}, or even direction of arrival estimation~\cite{moreira_doa_2019,raia2020gsp}. Solving these problems requires computing the solution of large numerical-linear-algebra (NLA) problems, whose formulation relies on the \emph{connection Laplacian} $\mathsf{L}_\theta$ (the formal definition of $\mathsf{L}_\theta$ is postponed to Section~\ref{sect:graphs:connections}), a linear operator encoding both the topology of the graph \emph{and} the rotations $\theta_{i,j}$.
    
    Computing an exact solution to many of these  problems is prohibitive beyond moderately large graphs ($n \simeq 10^4$)\footnote{For instance, Tikhonov smoothing requires solving a linear system, with generic time complexity in $\mathcal{O}(n^3)$.}, so that approximations via iterative Krylov-subspace-based algorithms \nt{(such as the popular Conjugate Gradient algorithm)} are typically used instead~\cite{saad2003iterative,saad2011numerical}. The main drawback of these methods is their spectrum-dependent convergence rate~\cite{kuijlaars2006convergence}, \nt{that often requires hard-to-find, good-quality preconditioners to ensure fast convergence.}

    Randomized numerical linear algebra (RandNLA) is a successful and modern alternative~\cite{drineas_randnla_2016,martinsson_randomized_2020}:  Monte-Carlo estimators allow flexible schemes of computation (\emph{e.g.}, parallelized or distributed) and can exhibit both advantageous complexities and  state-of-the-art practical performances, \emph{in spite of} their slow convergence rates in $O\left( \frac{\sigma}{\sqrt{m}} \right)$ (with $\sigma^2$ the variance of the estimator and $m$ the number of Monte-Carlo samples). Methods specialized in graphs, based on random walks and decompositions of the graph, have also appeared  (\emph{e.g.}~\cite{spielman2008graph,kyng2018matrix,pilavci_graph_2021}).
    
    In this paper, we leverage novel connection-aware random decompositions of the graph, and propose RandNLA estimators for two connection-Laplacian-based problems:
    \begin{itemize}
        \item A graph Tikhonov smoothing problem in the presence of a connection.
        \item The angular synchronization problem on graphs.
    \end{itemize}

    \begin{paragraph}{Tikhonov smoothing}
        Connection-aware graph Tikhonov smoothing amounts to solving:
        \begin{equation}
            f_*=\argmin_{f \in \C^n} \ q \Vert f - g \Vert_2^2 + \frac{1}{2} \sum_{i,j} \vert f(j) - e^{\iota \theta_{i,j}} f(i) \vert^2,
            \label{eq:tikhonov_intro}
        \end{equation}
        where $g \in \C^n$ \nt{is a complex graph signal one wishes to smooth (\textit{e.g.,} denoise)}, $q \in \R_+^*$ is an \emph{a priori} known regularization parameter, and the sum is over a subset of ordered pairs of indices $(i,j)$ describing the edges of the graph. This problem appears in different contexts, such as vector field extension in discrete geometry processing~\cite{sharp2019vector}, or directed graph signal processing~\cite{furutani_gsp_2019}. Here, the regularization term \nt{(the sum)} penalizes functions that are not \emph{locally coherent} with respect to the connection. The solution \nt{to Eq.~\eqref{eq:tikhonov_intro}} can be expressed using the connection Laplacian $\mathsf{L}_\theta$: $f_*=q(\mathsf{L}_\theta + q\mathsf{I})^{-1} g$. \nt{This inverse operation has a cubic complexity in $n$ and is prohibitive for large graphs, requiring the use of efficient estimators.} 
        \color{col}{Connection-aware Tikhonov smoothing also appears as an intermediate step in angular synchronization (the second problem we consider) as well as in iterative solvers for yet other problems\footnote{This is for instance the case for the so-called ``edge-lasso'' over graphs~\cite{sharpnack2012sparsistency} (which only appears in the literature on graphs \emph{without} a connection, but is relevant on graphs with a connection as well), which corresponds to an $l_1$-regularization problem over graph signals. The solution to this problem can for instance be obtained from iteratively-reweighted least-squares (see, \emph{e.g.},~\cite{pilavci_graph_2021}), or from an ADMM strategy (\emph{e.g.},~\cite{gabay1976dual}), or other proximal-operator-based algorithms~\cite{parikh2014proximal}; all of these strategies may involve solving graph Tikhonov regularization problems.}.} 
    \end{paragraph}
    
    \begin{paragraph}{Angular Synchronization}
        The objective is to recover a set of $n$ unknown angles $\omega = (\omega)_i \subseteq [0,2\pi)^n$ from measured pairwise offset measurements $\{\theta_{i,j}\}_{i,j}$~\cite{singer2011angular}~:
        \begin{equation}
            \theta_{i,j} = \omega_j - \omega_i + \eps_{i,j} \mod 2\pi,
        \end{equation}
        where $\eps_{i,j}$ represents some unknown degradation of the measurement. 
        This task appears in many structured signal processing problems, where it is often a key component in state-of-the-art recovery methods, such as perceived luminance reconstruction~\cite{yu_angular_2009}, ptychography~\cite{filbir_recovery_2021}, ranking~\cite{cucuringu_syncrank_2016}, clock synchronization~\cite{giridhar_distributed_2006} or phase reconstruction~\cite{alexeev_phase_2014}, and also appears in the statistical physics literature (\emph{e.g.}~\cite{chen_planted_2022}).
        \\
        In practice, we may only observe a subset of all such measurements $\theta_{i,j}$, and the problem is naturally formulated on a graph $\mathcal{G}$ with $n$ nodes whose set of edges $\mathcal{E}$ is indexed by the number of measurements, and where edge $(i,j)$ (resp. $(j,i)$) carries the offset $\theta_{i,j}$ (resp. ($\theta_{j,i} = -\theta_{i,j}$)). 
        \\
        If noiseless measurements $\theta_{i,j} = \omega_j - \omega_i$ are available, exact recovery can be trivially performed up to a global phase shift, by \emph{propagating} values according to the offset measurements along a spanning tree of $\mathcal{G}$ (a procedure explained in Figures~\ref{fig:ex_propagation:basics} and~\ref{fig:ex_propagation:tree}).
        \\
        The problem becomes more involved when considering imperfect measurements $\theta_{i,j}$, with non-zero noise $\eps_{i,j}$. In general, \emph{exact} angular synchronization can no longer be performed in this noisy regime, as long as even one \emph{incoherent} cycle is present in the graph (see Figure~\ref{fig:ex_propagation:incoherent}). 
        
        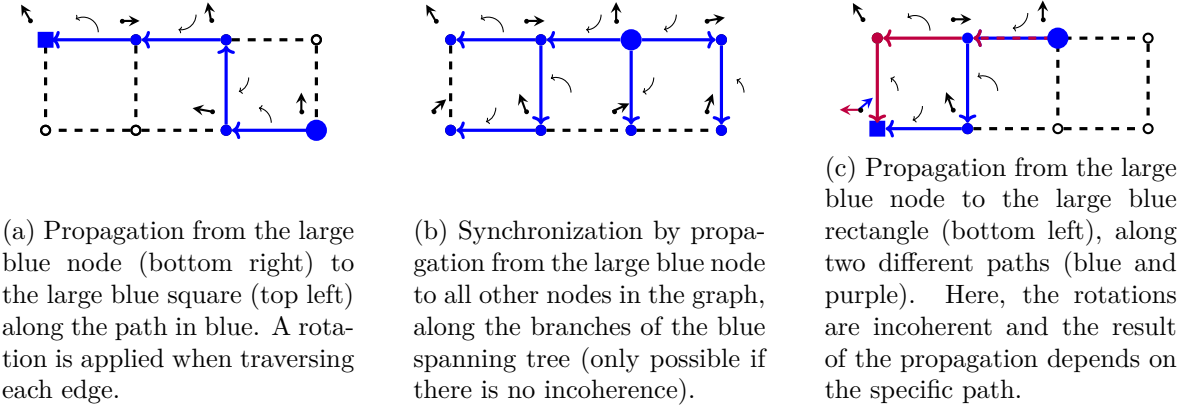
\begin{figure}[ht]
        \centering
        
        \begin{subfigure}{0.3\linewidth}
        \centering
        \scalebox{1}{\begin{tikzpicture}[every node/.style={circle,inner sep=1pt,thick,draw,scale=1.2},
	                        every edge/.style={very thick,draw}]
	                        
	       \def\centerarc[#1](#2)(#3:#4:#5)
    { \draw[#1] ($(#2)+({#5*cos(#3)},{#5*sin(#3)})$) arc (#3:#4:#5); }

            \node (1) at (0,0){};
            \node (2) [right of=1]{};
            \node (3) [right of=2]{};
            \node (4) [right of=3]{};
            \node (5) [below of=1]{};
            \node (6) [right of=5]{};
            \node (7) [right of=6]{};
            \node (8) [right of=7]{};
            
            \path (1) edge[dashed] (5);
            \path (6) edge[dashed] (7);
            \path (8) edge[-{>[scale=0.75]},color=blue] (7);

            \node (3) [right of=2,color=blue,fill]{};
            \path (3) edge[-{>[scale=0.75]},color=blue] (2); \node (2) [right of=1,color=blue,fill]{};
            \path (2) edge[dashed,color=black] (6); \node (6) [right of=5,color=black]{};
            \path (6) edge[dashed,color=black] (5); \node (5) [below of=1,color=black]{};
            
            \path (2) edge[-{>[scale=0.75]},color=blue] (1); \node (1) [left of=2, color=blue,fill,shape=rectangle,very thick,inner sep=2pt,fill]{};
            
            \path (7) edge[-{>[scale=0.75]},color=blue] (3); \node (7) [below of=3,color=blue,fill]{};

            \path (3) edge[dashed,color=black] (4); \node (4) [right of=3,color=black]{};
            \path (4) edge[dashed,color=black] (8); \node (8) [below of=4,color=blue,very thick,inner sep=2pt,fill]{};

            \centerarc[-{>[scale=0.5]}](1/2,1/10)(0:120:0.2);
            
            \centerarc[-{>[scale=0.5]}](1.8,1/3)(0:(-100):0.2);
            
            \centerarc[-{>[scale=0.5]}](2.5,(-1/2)(0:(-70):0.2);
            
            \centerarc[-{>[scale=0.5]}](2.8,-1.1)(0:(80):0.2);

            \draw[thick,-{stealth[scale=0.3]}] (1)++(-1/5,1/4) -- +(120:0.25);
            \node[above left =0.125 and 0.075 of 1,scale=0.3,fill,thick] (c1){};

            \draw[thick,-{stealth[scale=0.3]}] (2)++(-1/5,1/4) -- +(0:0.25);
            \node[above left=0.175 and 0.11 of 2,scale=0.3,fill,thick] (c2){};

            \draw[thick,-{stealth[scale=0.3]}] (3)++(-1/5,1/4) -- +(100:0.25);
            \node[above left=0.185 and 0.12 of 3,scale=0.3,fill,thick] (c3){};
            
            \draw[thick,-{stealth[scale=0.3]}] (7)++(-1/5,1/4) -- +(170:0.25);
            \node[above left=0.175 and 0.11 of 7,scale=0.3,fill,thick] (c7){};
            
            \draw[thick,-{stealth[scale=0.3]}] (8)++(-1/5,1/4) -- +(90:0.25);
            \node[above left=0.125 and 0.075 of 8,scale=0.3,fill,thick] (c8){};
            
            \node[color=white,scale=0.001] (20) at (0,-(2.125){};
            
        \end{tikzpicture}}
        \caption{Propagation from the large blue node (bottom right) to the large blue square (top left) along the path in blue. A rotation is applied when traversing each edge.}
        \label{fig:ex_propagation:basics}
        \end{subfigure}
        \hfill
        \begin{subfigure}{0.3\linewidth}
        \centering
        \scalebox{1}{\begin{tikzpicture}[every node/.style={circle,inner sep=1pt,thick,draw,scale=1.2},
	                        every edge/.style={very thick,draw}]
	                        
	       \def\centerarc[#1](#2)(#3:#4:#5)
    { \draw[#1] ($(#2)+({#5*cos(#3)},{#5*sin(#3)})$) arc (#3:#4:#5); }

            \node (1) at (0,0){};
            \node (2) [right of=1]{};
            \node (3) [right of=2]{};
            \node (4) [right of=3]{};
            \node (5) [below of=1]{};
            \node (6) [right of=5]{};
            \node (7) [right of=6]{};
            \node (8) [right of=7]{};
            
            \path (1) edge[dashed] (5);
            \path (6) edge[dashed] (7);
            \path (7) edge[dashed] (8);

            \node (3) [right of=2,color=blue,very thick, inner sep=2pt,fill]{};
            \path (3) edge[-{>[scale=0.75]},color=blue] (2); \node (2) [right of=1,color=blue,fill]{};
            \path (2) edge[-{>[scale=0.75]},color=blue] (6); \node (6) [right of=5,color=blue,fill]{};
            \path (6) edge[-{>[scale=0.75]},color=blue] (5); \node (5) [below of=1,color=blue,fill]{};
            
            \path (2) edge[-{>[scale=0.75]},color=blue] (1); \node (1) [left of=2, color=blue,fill]{};
            
            \path (3) edge[-{>[scale=0.75]},color=blue] (7); \node (7) [below of=3,color=blue,fill]{};

            \path (3) edge[-{>[scale=0.75]},color=blue] (4); \node (4) [right of=3,color=blue,fill]{};
            \path (4) edge[-{>[scale=0.75]},color=blue] (8); \node (8) [below of=4,color=blue,fill]{};

            \centerarc[-{>[scale=0.5]}](1/2,1/10)(0:120:0.2);
            
            \centerarc[-{>[scale=0.5]}](1.8,1/3)(0:(-90):0.2);
            
            \centerarc[-{>[scale=0.5]}](2.8,1/3)(0:(-80):0.2);
            
            \centerarc[-{>[scale=0.5]}](3.7,-(1/5 + 1/2)(0:(60):0.2);
            
            \centerarc[-{>[scale=0.5]}](2.5,(-1/2)(0:(-75):0.2);
            
            \centerarc[-{>[scale=0.5]}](1.4,-(1/5 + 1/2)(0:(110):0.2);
            
            \centerarc[-{>[scale=0.5]}](1/2-1/10,-9/10)(0:(-70):0.2);

            \draw[thick,-{stealth[scale=0.3]}] (1)++(-1/5,1/4) -- +(120:0.25);
            \node[above left =0.175 and 0.12 of 1,scale=0.3,fill,thick] (c1){};

            \draw[thick,-{stealth[scale=0.3]}] (2)++(-1/5,1/4) -- +(0:0.25);
            \node[above left=0.175 and 0.12 of 2,scale=0.3,fill,thick] (c2){};
            
            \draw[thick,-{stealth[scale=0.3]}] (3)++(-1/5,1/4) -- +(90:0.25);
            \node[above left=0.125 and 0.075 of 3,scale=0.3,fill,thick] (c3){};
            
            \draw[thick,-{stealth[scale=0.3]}] (4)++(-1/5,1/4) -- +(10:0.25);
            \node[above left=0.175 and 0.12 of 4,scale=0.3,fill,thick] (c4){};
            
            \draw[thick,-{stealth[scale=0.3]}] (8)++(-1/5,1/4) -- +(70:0.25);
            \node[above left=0.175 and 0.12 of 8,scale=0.3,fill,thick] (c8){};
            
            \draw[thick,-{stealth[scale=0.3]}] (7)++(-1/4,1/4) -- +(25:0.25);
            \node[above left=0.175 and 0.15 of 7,scale=0.3,fill,thick] (c7){};
            
            \draw[thick,-{stealth[scale=0.3]}] (6)++(-1/5,1/4) -- +(110:0.25);
            \node[above left=0.175 and 0.12 of 6,scale=0.3,fill,thick] (c6){};
            
            \draw[thick,-{stealth[scale=0.3]}] (5)++(-1/4,1/4) -- +(40:0.25);
            \node[above left=0.175 and 0.15 of 5,scale=0.3,fill,thick] (c5){};

            \node[color=white,scale=0.001] (20) at (0,-(2.125){};

        \end{tikzpicture}}
        \caption{Synchronization by propagation from the large blue node to all other nodes in the graph, along the branches of the blue spanning tree (only possible if there is no incoherence).}
        \label{fig:ex_propagation:tree}
        \end{subfigure}
        \hfill
        \begin{subfigure}{0.3\linewidth}
        \centering
        \scalebox{1}{\begin{tikzpicture}[every node/.style={circle,inner sep=1pt,thick,draw,scale=1.2},
	                        every edge/.style={very thick,draw}]

            \def\centerarc[#1](#2)(#3:#4:#5)
    { \draw[#1] ($(#2)+({#5*cos(#3)},{#5*sin(#3)})$) arc (#3:#4:#5); }

            \node (1) at (0,0){};
            \node (3) [right of=2]{};
            \node (4) [right of=3]{};
            \node (5) [below of=1]{};
            \node (6) [right of=5]{};
            \node (7) [right of=6]{};
            \node (8) [right of=7]{};
            
            \path (3) edge[dashed] (7);
            \path (3) edge[dashed] (4);
            \path (4) edge[dashed] (8);
            \path (6) edge[dashed] (7);
            \path (7) edge[dashed] (8);

            \node (5) [below of=1,color=blue,very thick, inner sep=2pt,fill,shape=rectangle]{};
            
            \path (3) edge[-{>[scale=0.75]},color=blue] (2); 
            \path (3) edge[dashed,-{>[scale=0.75]},color=purple] (2);
            \node (2) [right of=1,color=blue,fill]{};
            \path (2) edge[-{>[scale=0.75]},color=blue] (6); \node (6) [right of=5,color=blue,fill]{};
            \path (6) edge[-{>[scale=0.75]},color=blue] (5); \node (5) [below of=1,color=blue,fill]{};

            \path (2) edge[-{>[scale=0.75]},color=purple] (1); \node (1) [left of=2,color=purple,fill]{};
            \path (1) edge[-{>[scale=0.75]},color=purple] (5); \node (5) [below of=1,color=blue,fill]{};

            \node (3) [right of=2,color=blue,very thick, inner sep=2pt,fill]{};

            \centerarc[-{>[scale=0.5]}](1/2,1/10)(0:120:0.2);
            
            \centerarc[-{>[scale=0.5]}](1.8,1/3)(0:(-90):0.2);
            
            \centerarc[-{>[scale=0.5]}](1.4,-(1/5 + 1/2)(0:(110):0.2);
            
            \centerarc[-{>[scale=0.5]}](1/2-1/10,-9/10)(0:(-70):0.2);
            
            \centerarc[-{>[scale=0.5]}](0.075,-(1/2 + 1/5)(0:(60):0.2);
            
            \draw[thick,-{stealth[scale=0.3]}] (1)++(-1/5,1/4) -- +(120:0.25);
            \node[above left =0.175 and 0.12 of 1,scale=0.3,fill,thick] (c1){};

            \draw[thick,-{stealth[scale=0.3]}] (2)++(-1/5,1/4) -- +(0:0.25);
            \node[above left=0.175 and 0.12 of 2,scale=0.3,fill,thick] (c2){};
            
            \draw[thick,-{stealth[scale=0.3]}] (3)++(-1/5,1/4) -- +(90:0.25);
            \node[above left=0.125 and 0.075 of 3,scale=0.3,fill,thick] (c3){};
            
            \draw[thick,-{stealth[scale=0.3]}] (6)++(-1/5,1/4) -- +(110:0.25);
            \node[above left=0.175 and 0.12 of 6,scale=0.3,fill,thick] (c6){};
            
            \draw[color=blue,thick,-{stealth[scale=0.3]}] (5)++(-1/4,1/4) -- +(40:0.25);
            \node[above left=0.175 and 0.15 of 5,scale=0.3,fill,thick] (c5){};
            
            \draw[color=purple,thick,-{stealth[scale=0.3]}] (5)++(-1/4,1/4) -- +(180:0.25);
            \node[above left=0.175 and 0.15 of 5,scale=0.3,fill,thick] (c5b){};

        \end{tikzpicture}}
        \caption{Propagation from the large blue node to the large blue rectangle (bottom left), along two different paths (blue and purple). Here, the rotations are incoherent and the result of the propagation depends on the specific path.}
        \label{fig:ex_propagation:incoherent}
        \end{subfigure}
        
        \caption{Propagations on the $4 \times 2$ grid graph, along different paths. The rotation angles $\theta_{i,j}$'s associated to each edge are represented as circular rotating arrows (only drawn along the paths we consider). Propagations always start from the large blue node, and propagated angles are represented as straight arrows (top left of each node).  Left~\ref{fig:ex_propagation:basics}: propagation along an arbitrary path. Center~\ref{fig:ex_propagation:tree}: exact synchronization achieved by propagation along a spanning tree, in the absence of noise. Right~\ref{fig:ex_propagation:incoherent}: synchronization impossible due to noise and incoherence along cycles.}
        \label{fig:ex_propagation}
        \end{figure}

        \noindent A common workaround, first introduced in~\cite{yu_angular_2009}, consists in quantifying the incoherence of an angular assignment $s \in [0,2\pi)^n$ as 
        \begin{equation}
            \mathcal{I}(s) = \sum_{\{i,j\} \in \mathcal{E}} \left( 2 - 2 \cos{\left( (s_j - s_i) - \theta_{i,j} \right)} \right),
            \label{eq:cosine_loss}
        \end{equation}
        before minimizing this incoherence over all possible assignments. While this problem is non-convex and $\mathrm{NP}$-hard in general~\cite{zhang2006complex, boumal_nonconvex_2016}, different techniques allowing to recover a solution have been proposed~\cite{yu_angular_2012,singer2011angular,boumal_nonconvex_2016,perry2018message,he_robust_2023}, be it approximately via relaxations or for specific noise regimes or topologies. Most of these methods rely on the \emph{connection Laplacian} $\mathsf{L}_\theta$.
    \end{paragraph}

    \subsection*{Our contributions}
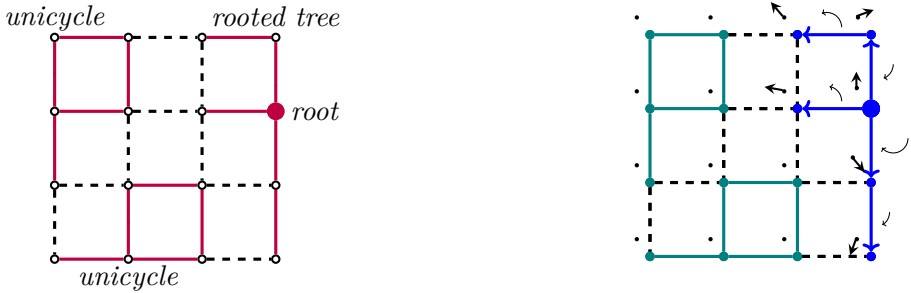
\begin{figure}[t]
    \centering
        
        \begin{subfigure}{0.49\textwidth}
            \centering
            \scalebox{0.98}{\begin{tikzpicture}[every node/.style={circle,inner sep=1pt,thick,draw},
	                        every edge/.style={very thick,draw}]
	        \node (1) at (0,0){};
            \node (2) [right of=1]{};
            \node (3) [right of=2]{};
            \node (4) [right of=3]{};
            \node (5) [below of=1]{};
            \node (6) [right of=5]{};
            \node (7) [right of=6]{};
            \node (8) [right of=7,color=purple,fill,very thick,inner sep=2pt,label={right:\emph{root}}]{};
            \node (9) [below of=5]{};
            \node (10) [right of=9]{};
            \node (11) [right of=10]{};
            \node (12) [right of=11]{};
            \node (13) [below of=9]{};
            \node (14) [right of=13]{};
            \node (15) [right of=14]{};
            \node (16) [right of=15]{};
            \path [-,color=purple] (1) edge (2);
            \path [-,color=black] (2) edge[dashed] (3);
            \path [-,color=purple] (3) edge (4);
            \path [-,color=purple] (5) edge (6);
            \path [-,color=black] (6) edge[dashed] (7);
            \path [-,color=purple] (7) edge (8);
            \path [-,color=black] (9) edge[dashed] (10);
            \path [-,color=purple] (10) edge (11);
            \path [-,color=black] (11) edge[dashed] (12);
            \path [-,color=purple] (13) edge (14);
            \path [-,color=purple] (14) edge (15);
            \path [-,color=black] (15) edge[dashed] (16);
            \path [-,color=purple] (1) edge (5);
            \path [-,color=purple] (2) edge (6);
            \path [-,color=black] (3) edge[dashed] (7);
            \path [-,color=purple] (4) edge (8);
            \path [-,color=purple] (5) edge (9);
            \path [-,color=black] (6) edge[dashed] (10);
            \path [-,color=black] (7) edge[dashed] (11);
            \path [-,color=purple] (8) edge (12);
            \path [-,color=black] (9) edge[dashed] (13);
            \path [-,color=purple] (10) edge (14);
            \path [-,color=purple] (11) edge (15);
            \path [-,color=purple] (12) edge (16);
            \textbf{}
            \node[label={above:\emph{unicycle}},fill=white,scale=0.0001] (a) at (0,-0.5){};
            \node[label={below:\emph{unicycle}},fill=white,scale=0.0001] (b) at (1,-2.5){};
            \node[label={above:\emph{rooted tree}},fill=white,scale=0.0001] (c) at (3,-0.65){};
        \end{tikzpicture}}
            \caption{A MTSF in purple, with two unicycles (top left, bottom) and one tree (right). The large purple node is the root of the tree.}
        \end{subfigure}
        \hfill
        \begin{subfigure}{0.49\textwidth}
            \centering
            \scalebox{0.98}{\begin{tikzpicture}[every node/.style={circle,inner sep=1pt,thick,draw},
	                        every edge/.style={very thick,draw}]
         
	        \node[color=teal,fill,label={above:\small}] (1) at (0,0){};
            \node[color=teal,fill] (2) [right of=1]{};
            \node (3) [right of=2]{};
            \node[label={right:\small}] (4) [right of=3]{};
            \node[color=teal,fill] (5) [below of=1]{};
            \node[color=teal,fill] (6) [right of=5]{};
            \node (7) [right of=6]{};
            \node (8) [right of=7]{};
            \node[color=teal,fill] (9) [below of=5]{};
            \node[color=teal,fill] (10) [right of=9]{};
            \node[color=teal,fill] (11) [right of=10]{};
            \node (12) [right of=11]{};
            \node[color=teal,fill] (13) [below of=9]{};
            \node[color=teal,fill] (14) [right of=13]{};
            \node[color=teal,fill] (15) [right of=14]{};
            \node (16) [right of=15]{};   

            \node (8) [right of=7,color=blue,fill,very thick,inner sep=2pt]{};
            \path [-,color=black] (2) edge[dashed] (3);
            \path [-,color=black] (6) edge[dashed] (7);
            \path [-,color=black] (9) edge[dashed] (10);
            \path [-,color=black] (11) edge[dashed] (12);
            \path [-,color=black] (15) edge[dashed] (16);
            \path [-,color=black] (3) edge[dashed] (7);
            \path [-,color=black] (6) edge[dashed] (10);
            \path [-,color=black] (7) edge[dashed] (11);
            \path [-,color=black] (9) edge[dashed] (13);
            
            \node (8) [right of=7,color=blue,fill,very thick,inner sep=2pt,label={right:\footnotesize}]{};
            \path [-,color=black] (2) edge[dashed] (3);
            \path [-,color=black] (6) edge[dashed] (7);
            \path [-,color=black] (9) edge[dashed] (10);
            \path [-,color=black] (11) edge[dashed] (12);
            \path [-,color=black] (15) edge[dashed] (16);
            \path [-,color=black] (3) edge[dashed] (7);
            \path [-,color=black] (6) edge[dashed] (10);
            \path [-,color=black] (7) edge[dashed] (11);
            \path [-,color=black] (9) edge[dashed] (13);

            \node (8) [right of=7,color=blue,fill,very thick,inner sep=2pt]{};
            \path [{<[scale=0.75]}-,color=blue] (4) edge (8);\node (4) [right of=3,color=blue,fill]{};
            \path [-{>[scale=0.75]},color=blue] (4) edge (3);\node (3) [right of=2,color=blue,fill]{};
            \path [{<[scale=0.75]}-,color=blue] (7) edge (8);\node (7) [right of=6,color=blue,fill]{};
            \path [{<[scale=0.75]}-,color=blue] (12) edge (8);\node (12) [right of=11,color=blue,fill]{};
            \path [-{>[scale=0.75]},color=blue] (12) edge (16);\node (16) [right of=15,color=blue,fill]{};

            \node (1) at (0,0)[color=teal,fill]{};
            \node (2) [right of=1,color=teal,fill]{};
            \node (5) [below of=1,color=teal,fill]{};
            \node (6) [right of=5,color=teal,fill]{};
            \node (9) [below of=5,color=teal,fill]{};
            \node (10) [right of=9,color=teal,fill]{};
            \node (11) [right of=10,color=teal,fill]{};
            \node (13) [below of=9,color=teal,fill]{};
            \node (14) [right of=13,color=teal,fill]{};
            \node (15) [right of=14,color=teal,fill]{};
            
            \path [-,color=teal] (1) edge (2);
            \path [-,color=teal] (5) edge (6);
            \path [-,color=teal] (10) edge (11);
            \path [-,color=teal] (13) edge (14);
            \path [-,color=teal] (14) edge (15);
            \path [-,color=teal] (1) edge (5);
            \path [-,color=teal] (2) edge (6);
            \path [-,color=teal] (5) edge (9);
            \path [-,color=teal] (10) edge (14);
            \path [-,color=teal] (11) edge (15);
            
            \node[above left=0.175 and 0.12 of 1,scale=0.3,fill,thick] (c5){};
            \node[above left=0.175 and 0.12 of 2,scale=0.3,fill,thick] (c5){};
            \node[above left=0.175 and 0.12 of 5,scale=0.3,fill,thick] (c5){};
            \node[above left=0.175 and 0.12 of 6,scale=0.3,fill,thick] (c5){};
            \node[above left=0.175 and 0.12 of 9,scale=0.3,fill,thick] (c5){};
            \node[above left=0.175 and 0.12 of 10,scale=0.3,fill,thick] (c5){};
            \node[above left=0.175 and 0.12 of 11,scale=0.3,fill,thick] (c5){};
            \node[above left=0.175 and 0.12 of 13,scale=0.3,fill,thick] (c5){};
            \node[above left=0.175 and 0.12 of 14,scale=0.3,fill,thick] (c5){};
            \node[above left=0.175 and 0.12 of 15,scale=0.3,fill,thick] (c5){};
            
            \node[above left=0.175 and 0.12 of 3,scale=0.3,fill,thick] (c5){};
            \node[above left=0.175 and 0.12 of 4,scale=0.3,fill,thick] (c5){};
            \node[above left=0.175 and 0.12 of 7,scale=0.3,fill,thick] (c5){};
            \node[above left=0.18 and 0.095 of 8,scale=0.3,fill,thick] (c5){};
            \node[above left=0.28 and 0.19 of 12,scale=0.3,fill,thick] (c5){};
            \node[above left=0.175 and 0.135 of 16,scale=0.3,fill,thick] (c5){};
            
            \draw[thick,-{stealth[scale=0.3]}] (3)++(-1/5,1/4) -- +(130:0.25);
            \draw[thick,-{stealth[scale=0.3]}] (4)++(-1/5,1/4) -- +(20:0.25);
            \draw[thick,-{stealth[scale=0.3]}] (7)++(-1/5,1/4) -- +(170:0.25);
            \draw[thick,-{stealth[scale=0.3]}] (8)++(-1/5,1/4) -- +(90:0.25);
            \draw[thick,-{stealth[scale=0.3]}] (12)++(-1/4,1/3) -- +(-50:0.25);
            \draw[thick,-{stealth[scale=0.3]}] (16)++(-1/5,1/4) -- +(-110:0.25);
        
            \def\centerarc[#1](#2)(#3:#4:#5)
    { \draw[#1] ($(#2)+({#5*cos(#3)},{#5*sin(#3)})$) arc (#3:#4:#5); }
    
            \centerarc[-{>[scale=0.5]}](3.1,-0.4)(0:(-70):0.2);
            \centerarc[-{>[scale=0.5]}](2.4,1/10)(0:(110):0.2);
            \centerarc[-{>[scale=0.5]}](3.3,-1.4)(0:(-140):0.2);
            \centerarc[-{>[scale=0.5]}](3.05,-2.4)(0:(-60):0.2);
            \centerarc[-{>[scale=0.5]}](2.4,-0.9)(0:(80):0.2);

            \node[color=white,scale=0.001] (20) at (0,-3.6){};
        \end{tikzpicture}}
            \caption{\nt{Estimation. Tree: the value of the root (blue circle) is propagated to the rest of the tree via the blue arrows. Unicycles: the estimation is $0$ (arrowless dots at the top left of the nodes).}}
        \end{subfigure}
        
        \caption{Theorem~\ref{th:teaser} illustrated on the $4 \times 4$ grid graph. Left: a MTSF. Right: estimation by propagations along the branches of the trees of the MTSF. The estimation in $\mathbf{C}$ is represented using angled arrows (top left of the nodes).}
        \label{fig:ex_th}
    \end{figure}

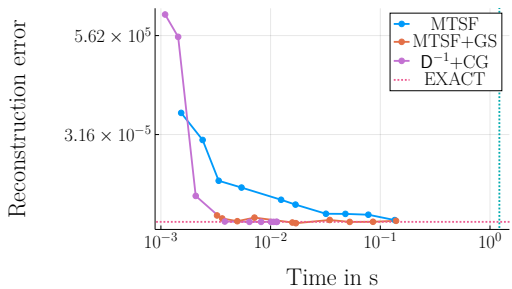
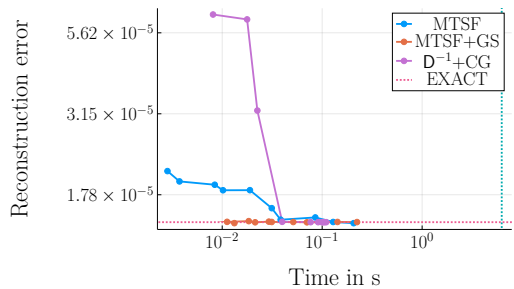
\begin{figure}[ht]
        \centering
        
        \begin{subfigure}{0.49\textwidth}
            \scalebox{0.35}{

\begin{tikzpicture}[/tikz/background rectangle/.style={fill={rgb,1:red,1.0;green,1.0;blue,1.0}, fill opacity={1.0}, draw opacity={1.0}}, show background rectangle]
\begin{axis}[point meta max={nan}, point meta min={nan}, legend cell align={left}, legend columns={1}, title={}, title style={at={{(0.5,1)}}, anchor={south}, font={{\fontsize{14 pt}{18.2 pt}\selectfont}}, color={rgb,1:red,0.0;green,0.0;blue,0.0}, draw opacity={1.0}, rotate={0.0}, align={center}}, legend style={color={rgb,1:red,0.0;green,0.0;blue,0.0}, draw opacity={1.0}, line width={1}, solid, fill={rgb,1:red,1.0;green,1.0;blue,1.0}, fill opacity={1.0}, text opacity={1.0}, font={{\fontsize{20 pt}{10.4 pt}\selectfont}}, text={rgb,1:red,0.0;green,0.0;blue,0.0}, cells={anchor={center}}, at={(0.98, 0.98)}, anchor={north east}}, axis background/.style={fill={rgb,1:red,1.0;green,1.0;blue,1.0}, opacity={1.0}}, anchor={north west}, xshift={1.0mm}, yshift={-1.0mm}, width={150.4mm}, height={99.6mm}, scaled x ticks={false}, xlabel={Time in s}, x tick style={color={rgb,1:red,0.0;green,0.0;blue,0.0}, opacity={1.0}}, x tick label style={color={rgb,1:red,0.0;green,0.0;blue,0.0}, opacity={1.0}, rotate={0}}, xlabel style={at={(ticklabel cs:0.5)}, anchor=near ticklabel, at={{(ticklabel cs:0.5)}}, anchor={near ticklabel}, font={{\fontsize{26 pt}{14.3 pt}\selectfont}}, color={rgb,1:red,0.0;green,0.0;blue,0.0}, draw opacity={1.0}, rotate={0.0}}, xmode={log}, log basis x={10}, xmajorgrids={true}, xmin={0.0008850285277387312}, xmax={1.5031185137995866}, xticklabels={{$10^{-3}$,$10^{-2}$,$10^{-1}$,$10^{0}$}}, xtick={{0.001,0.01,0.1,1.0}}, xtick align={inside}, xticklabel style={font={{\fontsize{20 pt}{10.4 pt}\selectfont}}, color={rgb,1:red,0.0;green,0.0;blue,0.0}, draw opacity={1.0}, rotate={0.0}}, x grid style={color={rgb,1:red,0.0;green,0.0;blue,0.0}, draw opacity={0.1}, line width={0.5}, solid}, axis x line*={left}, x axis line style={color={rgb,1:red,0.0;green,0.0;blue,0.0}, draw opacity={1.0}, line width={1}, solid}, scaled y ticks={false}, ylabel={Reconstruction error}, y tick style={color={rgb,1:red,0.0;green,0.0;blue,0.0}, opacity={1.0}}, y tick label style={color={rgb,1:red,0.0;green,0.0;blue,0.0}, opacity={1.0}, rotate={0}}, ylabel style={at={(ticklabel cs:0.5)}, anchor=near ticklabel, at={{(ticklabel cs:0.5)}}, anchor={near ticklabel}, font={{\fontsize{26 pt}{14.3 pt}\selectfont}}, color={rgb,1:red,0.0;green,0.0;blue,0.0}, draw opacity={1.0}, rotate={0.0}}, ymode={log}, log basis y={10}, ymajorgrids={true}, ymin={1.8208601816326163e-5}, ymax={6.592756393953906e-5}, yticklabels={{\(3.16 \times 10^{-5}\),\(5.62 \times 10^{5}\)}}, ytick={{3.1622776601683795e-5,5.623413251903491e-5}}, ytick align={inside}, yticklabel style={font={{\fontsize{20 pt}{10.4 pt}\selectfont}}, color={rgb,1:red,0.0;green,0.0;blue,0.0}, draw opacity={1.0}, rotate={0.0}}, y grid style={color={rgb,1:red,0.0;green,0.0;blue,0.0}, draw opacity={0.1}, line width={0.5}, solid}, axis y line*={left}, y axis line style={color={rgb,1:red,0.0;green,0.0;blue,0.0}, draw opacity={1.0}, line width={1}, solid}, colorbar={false}]
    \addplot[color={rgb,1:red,0.0;green,0.6056;blue,0.9787}, name path={63337f17-7245-4525-8ff5-69c45ef63830}, draw opacity={1.0}, line width={2}, solid, mark={*}, mark size={3.0 pt}, mark repeat={1}, mark options={color={rgb,1:red,0.0;green,0.6056;blue,0.9787}, draw opacity={1.0}, fill={rgb,1:red,0.0;green,0.6056;blue,0.9787}, fill opacity={1.0}, line width={0.75}, rotate={0}, solid}]
        table[row sep={\\}]
        {
            \\
            0.0015308698232000002  3.58683847044972e-5  \\
            0.0024008618785999994  3.064220908553685e-5  \\
            0.0033631265204  2.416880809224033e-5  \\
            0.005431243193599998  2.3228929799882224e-5  \\
            0.012462797532799998  2.1639479813323863e-5  \\
            0.0168686023766  2.1037033945762114e-5  \\
            0.0319691218302  1.994763611480616e-5  \\
            0.04808730418979999  1.99348635284733e-5  \\
            0.07780924720319998  1.98373931909197e-5  \\
            0.1353063975734  1.9229280453678935e-5  \\
        }
        ;
    \addlegendentry {MTSF}
    \addplot[color={rgb,1:red,0.8889;green,0.4356;blue,0.2781}, name path={15743a2d-486c-4320-8b0f-5ca1425ea110}, draw opacity={1.0}, line width={2}, solid, mark={*}, mark size={3.0 pt}, mark repeat={1}, mark options={color={rgb,1:red,0.8889;green,0.4356;blue,0.2781}, draw opacity={1.0}, fill={rgb,1:red,0.8889;green,0.4356;blue,0.2781}, fill opacity={1.0}, line width={0.75}, rotate={0}, solid}]
        table[row sep={\\}]
        {
            \\
            0.0032501763943999997  1.9754309854612547e-5  \\
            0.003622600276000001  1.9408258272925898e-5  \\
            0.004974376870400002  1.91048482906192e-5  \\
            0.0071287508732000016  1.950773105041196e-5  \\
            0.015786409657799996  1.895868487034726e-5  \\
            0.016996112962399997  1.888388827908881e-5  \\
            0.034595695165400016  1.9238859772626545e-5  \\
            0.052620975981799996  1.902530208278843e-5  \\
            0.085794776625  1.903927603236167e-5  \\
            0.1391596934482  1.9145526042606387e-5  \\
        }
        ;
    \addlegendentry {MTSF+GS}
    \addplot[color={rgb,1:red,0.7644;green,0.4441;blue,0.8243}, name path={d6c1c1c2-6109-422a-ad51-aa30acd30bec}, draw opacity={1.0}, line width={2}, solid, mark={*}, mark size={3.0 pt}, mark repeat={1}, mark options={color={rgb,1:red,0.7644;green,0.4441;blue,0.8243}, draw opacity={1.0}, fill={rgb,1:red,0.7644;green,0.4441;blue,0.8243}, fill opacity={1.0}, line width={0.75}, rotate={0}, solid}]
        table[row sep={\\}]
        {
            \\
            0.0010923790983999998  6.356999907825005e-5  \\
            0.0014361057674000001  5.582493698922853e-5  \\
            0.0020768403252  2.213021217430731e-5  \\
            0.0038239560488000004  1.9067786945840237e-5  \\
            0.0064090582356  1.90311705710683e-5  \\
            0.008190323317200002  1.9031260483443813e-5  \\
            0.010073858436199996  1.9031260423904773e-5  \\
            0.010642968164399998  1.9031260423904773e-5  \\
            0.0112713690074  1.9031260423904773e-5  \\
            0.011425916066200003  1.9031260423904773e-5  \\
        }
        ;
    \addlegendentry {$\mathsf{D}^{-1}$+CG}
    \addplot[color={rgb,1:red,0.0;green,0.6658;blue,0.681}, name path={ad9f83a4-3a28-411b-8b15-faa83625017f}, draw opacity={1.0}, line width={2}, dotted, forget plot]
        table[row sep={\\}]
        {
            \\
            1.2178032033323996  5.02905249782288e-6  \\
            1.2178032033323996  0.00023870276975934037  \\
        }
        ;
    \addplot[color={rgb,1:red,0.9308;green,0.3675;blue,0.5758}, name path={c8641ab2-b31a-4aed-bdec-bc52bcebe883}, draw opacity={1.0}, line width={2}, dotted]
        table[row sep={\\}]
        {
            \\
            5.211002909753405e-7  1.903126042498225e-5  \\
            2552.872812246866  1.903126042498225e-5  \\
        }
        ;
    \addlegendentry {EXACT}
\end{axis}
\end{tikzpicture}}
            \caption{DC-SBM with mean degree $\overline{d} \simeq 30$.}
        \end{subfigure}
        \hfill
        \begin{subfigure}{0.49\textwidth}
            \scalebox{0.35}{

\begin{tikzpicture}[/tikz/background rectangle/.style={fill={rgb,1:red,1.0;green,1.0;blue,1.0}, fill opacity={1.0}, draw opacity={1.0}}, show background rectangle]
\begin{axis}[point meta max={nan}, point meta min={nan}, legend cell align={left}, legend columns={1}, title={}, title style={at={{(0.5,1)}}, anchor={south}, font={{\fontsize{14 pt}{18.2 pt}\selectfont}}, color={rgb,1:red,0.0;green,0.0;blue,0.0}, draw opacity={1.0}, rotate={0.0}, align={center}}, legend style={color={rgb,1:red,0.0;green,0.0;blue,0.0}, draw opacity={1.0}, line width={1}, solid, fill={rgb,1:red,1.0;green,1.0;blue,1.0}, fill opacity={1.0}, text opacity={1.0}, font={{\fontsize{20 pt}{10.4 pt}\selectfont}}, text={rgb,1:red,0.0;green,0.0;blue,0.0}, cells={anchor={center}}, at={(0.98, 0.98)}, anchor={north east}}, axis background/.style={fill={rgb,1:red,1.0;green,1.0;blue,1.0}, opacity={1.0}}, anchor={north west}, xshift={1.0mm}, yshift={-1.0mm}, width={150.4mm}, height={99.6mm}, scaled x ticks={false}, xlabel={Time in s}, x tick style={color={rgb,1:red,0.0;green,0.0;blue,0.0}, opacity={1.0}}, x tick label style={color={rgb,1:red,0.0;green,0.0;blue,0.0}, opacity={1.0}, rotate={0}}, xlabel style={at={(ticklabel cs:0.5)}, anchor=near ticklabel, at={{(ticklabel cs:0.5)}}, anchor={near ticklabel}, font={{\fontsize{26 pt}{14.3 pt}\selectfont}}, color={rgb,1:red,0.0;green,0.0;blue,0.0}, draw opacity={1.0}, rotate={0.0}}, xmode={log}, log basis x={10}, xmajorgrids={true}, xmin={0.002249003586429591}, xmax={7.882259078854072}, xticklabels={{$10^{-2}$,$10^{-1}$,$10^{0}$}}, xtick={{0.01,0.1,1.0}}, xtick align={inside}, xticklabel style={font={{\fontsize{20 pt}{10.4 pt}\selectfont}}, color={rgb,1:red,0.0;green,0.0;blue,0.0}, draw opacity={1.0}, rotate={0.0}}, x grid style={color={rgb,1:red,0.0;green,0.0;blue,0.0}, draw opacity={0.1}, line width={0.5}, solid}, axis x line*={left}, x axis line style={color={rgb,1:red,0.0;green,0.0;blue,0.0}, draw opacity={1.0}, line width={1}, solid}, scaled y ticks={false}, ylabel={Reconstruction error}, y tick style={color={rgb,1:red,0.0;green,0.0;blue,0.0}, opacity={1.0}}, y tick label style={color={rgb,1:red,0.0;green,0.0;blue,0.0}, opacity={1.0}, rotate={0}}, ylabel style={at={(ticklabel cs:0.5)}, anchor=near ticklabel, at={{(ticklabel cs:0.5)}}, anchor={near ticklabel}, font={{\fontsize{26 pt}{14.3 pt}\selectfont}}, color={rgb,1:red,0.0;green,0.0;blue,0.0}, draw opacity={1.0}, rotate={0.0}}, ymode={log}, log basis y={10}, ymajorgrids={true}, ymin={1.3896656991424517e-5}, ymax={6.696200597330384e-5}, yticklabels={{\(1.78 \times 10^{-5}\),\(3.15 \times 10^{-5}\),\(5.62 \times 10^{-5}\)}}, ytick={{1.778279410038923e-5,3.1622776601683795e-5,5.623413251903491e-5}}, ytick align={inside}, yticklabel style={font={{\fontsize{20 pt}{10.4 pt}\selectfont}}, color={rgb,1:red,0.0;green,0.0;blue,0.0}, draw opacity={1.0}, rotate={0.0}}, y grid style={color={rgb,1:red,0.0;green,0.0;blue,0.0}, draw opacity={0.1}, line width={0.5}, solid}, axis y line*={left}, y axis line style={color={rgb,1:red,0.0;green,0.0;blue,0.0}, draw opacity={1.0}, line width={1}, solid}, colorbar={false}]
    \addplot[color={rgb,1:red,0.0;green,0.6056;blue,0.9787}, name path={24963a4e-2389-43f5-b879-4f74a60ab373}, draw opacity={1.0}, line width={2}, solid, mark={*}, mark size={3.0 pt}, mark repeat={1}, mark options={color={rgb,1:red,0.0;green,0.6056;blue,0.9787}, draw opacity={1.0}, fill={rgb,1:red,0.0;green,0.6056;blue,0.9787}, fill opacity={1.0}, line width={0.75}, rotate={0}, solid}]
        table[row sep={\\}]
        {
            \\
            0.0028334185301999997  2.1040209817350433e-5  \\
            0.0037379942884000006  1.955374993018194e-5  \\
            0.0084193877026  1.9095491235337058e-5  \\
            0.010181012451600002  1.836096878096219e-5  \\
            0.0189446295446  1.8373916093626262e-5  \\
            0.031308930622200006  1.616412356286385e-5  \\
            0.038426807412599996  1.489290726566159e-5  \\
            0.08536121517520003  1.5134227100681767e-5  \\
            0.1287781440518  1.4652880673868332e-5  \\
            0.20646168981120003  1.4529083151318067e-5  \\
        }
        ;
    \addlegendentry {MTSF}
    \addplot[color={rgb,1:red,0.8889;green,0.4356;blue,0.2781}, name path={9b517cf8-ae39-4f91-aa22-6a0c9d3c3303}, draw opacity={1.0}, line width={2}, solid, mark={*}, mark size={3.0 pt}, mark repeat={1}, mark options={color={rgb,1:red,0.8889;green,0.4356;blue,0.2781}, draw opacity={1.0}, fill={rgb,1:red,0.8889;green,0.4356;blue,0.2781}, fill opacity={1.0}, line width={0.75}, rotate={0}, solid}]
        table[row sep={\\}]
        {
            \\
            0.0132382501892  1.4533952116559921e-5  \\
            0.011192393792799998  1.4666388467189694e-5  \\
            0.018433897098799997  1.4742905088188501e-5  \\
            0.021356811400000005  1.4620339303767194e-5  \\
            0.0292819788492  1.4690461868593096e-5  \\
            0.031296631864399994  1.465452016557325e-5  \\
            0.05138992778119999  1.4651529602723103e-5  \\
            0.0704137448438  1.4617802932388759e-5  \\
            0.142327679929  1.4655915604597402e-5  \\
            0.22315145818339996  1.4644127842001571e-5  \\
        }
        ;
    \addlegendentry {MTSF+GS}
    \addplot[color={rgb,1:red,0.7644;green,0.4441;blue,0.8243}, name path={249ee1ab-9245-49e6-893c-c40a02259f75}, draw opacity={1.0}, line width={2}, solid, mark={*}, mark size={3.0 pt}, mark repeat={1}, mark options={color={rgb,1:red,0.7644;green,0.4441;blue,0.8243}, draw opacity={1.0}, fill={rgb,1:red,0.7644;green,0.4441;blue,0.8243}, fill opacity={1.0}, line width={0.75}, rotate={0}, solid}]
        table[row sep={\\}]
        {
            \\
            0.008095425743399998  6.404726428895852e-5  \\
            0.017744916993200004  6.186110463768229e-5  \\
            0.022478603143000004  3.2360352773168454e-5  \\
            0.0398031081262  1.4671243133032395e-5  \\
            0.07672043277299999  1.4635048312711097e-5  \\
            0.10974329127839998  1.463356340219918e-5  \\
            0.09112879555300002  1.4633563382536939e-5  \\
            0.1027907322758  1.4633563382536939e-5  \\
            0.10555339700040003  1.4633563382536939e-5  \\
            0.09911791991000003  1.4633563382536939e-5  \\
        }
        ;
    \addlegendentry {$\mathsf{D}^{-1}$+CG}
    \addplot[color={rgb,1:red,0.0;green,0.6658;blue,0.681}, name path={0ecbcedd-0811-4cea-8b68-d962bffbdf4b}, draw opacity={1.0}, line width={2}, dotted, forget plot]
        table[row sep={\\}]
        {
            \\
            6.256480907625998  2.8839798439476127e-6  \\
            6.256480907625998  0.00032266107213668396  \\
        }
        ;
    \addplot[color={rgb,1:red,0.9308;green,0.3675;blue,0.5758}, name path={5e182299-7e60-4c4e-b5bf-f2d9cfd61646}, draw opacity={1.0}, line width={2}, dotted]
        table[row sep={\\}]
        {
            \\
            6.416963818586512e-7  1.4633563382513295e-5  \\
            27625.570968880504  1.4633563382513295e-5  \\
        }
        ;
    \addlegendentry {EXACT}
\end{axis}
\end{tikzpicture}}
            \caption{DC-SBM with mean degree $\overline{d} \simeq 420$.}
        \end{subfigure}
        
        \caption{Runtime-precision comparisons for the graph Tikhonov smoothing problem, comparing two MTSF-based estimators ({MTSF} and {MTSF+GS}) against a diagonally-preconditioned conjugate-gradient-descent ({$\mathsf{D}^{-1}$+CG}) on degree-corrected stochastic block model graphs (DC-SBM, see Section~\ref{sect:numerical} for a definition) with $10000$ nodes. Results are averaged over $10$ graphs generated from DC-SBM models with two different parametrizations, resulting in different average degrees $\overline{d}$: on the left, a DC-SBM with average degree $\overline{d} \simeq 30$ and, on the right, a DC-SBM with average degree $\overline{d} \simeq 420$.
        $x$-axis: average runtime. $y$-axis: average reconstruction error. 
        Each data-point corresponds to a number $m \in \{1,2,3,5,8,13,22,36,60,100\}$ of MTSFs (or of CG iterations) used in the estimation. The vertical (resp. horizontal) line is the computation-time (resp. error) of an exact Cholesky solver (see Section~\ref{sect:numerical} for details).}
        \label{fig:runtime_teaser}
    \end{figure}

    We propose novel RandNLA estimators for the smoothing and angular synchronization problems. We rely on propagations along branches of \emph{Multi-Type Spanning Forests} (MTSF): spanning subsets of both edges and nodes of a graph, whose connected components are either rooted trees or unicycles (connected subsets of edges containing \emph{exactly one} cycle). Our main theorem concerns our estimator for the graph Tikhonov smoothing problem, and can be roughly described as follows (see Theorem~\ref{th:fermionic} for a formal statement).
    \begin{thm}
        For MTSFs sampled from the correct distribution (the one of Eq.~\eqref{eq:prob_dm}), \nt{propagating within each tree the value of its root to all of its other nodes} yields an unbiased estimation of the solution of the connection-aware Tikhonov smoothing problem of Eq.~\eqref{eq:tikhonov_intro}.
        \label{th:teaser}
    \end{thm}
    See Figure~\ref{fig:ex_th} for an illustration. We include a teaser runtime comparison with a standard deterministic solver in Figure~\ref{fig:runtime_teaser}, which shows that our estimators are \emph{less sensitive} to the density of the graph, and can provide significant speed-ups for equivalent precision.

    \noindent In practice, we obtain a \emph{fast} and \emph{scalable} algorithm. MTSF-sampling is achieved via a Wilson-like random-walk-based algorithm with runtime \emph{linear} in the number of edges of the graph. On the theoretical side, our arguments rely on the theory of \emph{Determinantal Point Processes} (DPPs)~\cite{hough_determinantal_2006,macchi1975coincidence}, and generalize the combinatorial analyses of \cite{pilavci_graph_2021,kenyon_spanning_2011}.
    
    A preliminary version of this work already appeared in~\cite{jaquard_smoothing_2023}, where we presented our smoothing estimator along with two variance-reduction techniques, and an application to a ranking problem for a simple synthetic data model.
    We improve and extend this work in a number of different manners, some of which we list below.

    \begin{paragraph}{Theoretical results}
    In addition to Theorem~\ref{th:fermionic}, we derive a connection-aware Feynman-Kac formula (Proposition~\ref{prop:feynmac_kac}), resulting in a local (node-wise) random-walk based estimator for the smoothing problem, similar to walk-on-spheres algorithms~\cite{muller1956some,sawhney2020monte}. In comparison, our MTSF-based estimators allow for \emph{global} estimation on all the nodes at once.  We further relate these two estimators, which may pave the way to generalizations (see Section~\ref{sect:extension} of the Supplementary Material). We also improve one of our variance-reduction techniques to better handle heterogeneous degree distributions, with significant improvements.
    \end{paragraph}

    \begin{paragraph}{Methodological and experimental contributions}
    We compare our estimators with standard Krylov-subspaces methods on synthetic data, for both the angular synchronization and smoothing problems. For the angular synchronization problem, we leverage our smoothing estimator as an iterative step in existing approaches~\cite{yu_angular_2012,singer2011angular,boumal_nonconvex_2016}.
    We analyze and contrast the behavior of these methods on different graph topologies which, up to our knowledge, has never been studied in the literature before. Our results show that MTSF-based estimations can  outperform standard deterministic methods whenever the graph is not \emph{very sparse} (the gain getting bigger as the density increases). See Figure~\ref{fig:runtime_teaser}.
    The code used for these experiments is publicly available\footnote{\url{https://gricad-gitlab.univ-grenoble-alpes.fr/gaia/synchromtsf}}.
    \end{paragraph}

    \subsection*{Related Work}

    Our main result can be understood as a specialization of the approach in~\cite{derezinsky_reverse_2018,derezinski2021determinantal}, based on DPPs, a class of distributions exhibiting negatively-correlated samples~\cite{hough_determinantal_2006,kulesza_determinantal_2012}, and resulting in unbiased estimators for least-squares problems. \{However, and in contrast to generic DPPs that are in general expensive to sample from, the specific structure of MTSFs allows for very efficient sampling algorithms and, in turn, practical graph RandNLA algorithms that can compete with determinstic state-of-the-art solvers.

    There exists similar applications of random spanning forests distributions in RandNLA for problems on connection-free graphs~\cite{avena2013random,avena2018two}, such as the trace estimation of (regularized) inverse Laplacians~\cite{barthelme_inverse_2019} or, the Tikhonov regularization and interpolation of graphs signals~\cite{pilavci_graph_2021}.
    Multi-Type Spanning Forests have also been used to build spectral sparsifiers for the (regularized) connection Laplacian in \cite{fanuel_sparsification_2022}, resulting in randomized preconditioners. Our Feynman-Kac formula is inspired from a similar result in~\cite{kassel2021covariant} for continuous-time random walks, but is specially tailored to the graph Tikhonov smoothing problem; discrete-time walks also exhibit intriguing links with propagations on MTSFs (see Section~\ref{sect:extension} of the Supplementary Material).
    
    \subsection*{Organisation of the Paper}
    In Section~\ref{sect:graphs}, we introduce preliminary background on graphs and connections used in the rest of the paper. In Section~\ref{sect:smoothing}, we present our Feynman-Kac formula (Section~\ref{sect:smoothing:fk}) and our MTSF-based estimators (Section~\ref{sect:smoothing:mtsf}), analyze an efficient sampling algorithm for MTSFs (Section~\ref{sect:smoothing:mtsf_sampling}), and describe our variance-reduction techniques (Section~\ref{par:var_red}). We analyze the numerical behavior of our estimators in Section~\ref{sect:numerical}, and compare them with standard deterministic algorithms (Sections~\ref{sect:numerical:density} and~\ref{subsect:runtime_precision_smoothing}). In Section~\ref{sect:synchro}, we propose an iterative scheme leveraging our smoothing estimators to solve the angular synchronization problem (Section~\ref{subsect:synchro_algo}), illustrate its application to a denoising problem inspired from \emph{cryogenic electron microscopy} (Section~\ref{sect:illustration}), and compare this scheme with standard deterministic methods (Section~\ref{sect:synchro:eval}). Proofs are deferred to the Supplementary Material.
    
    \section{Background}
    \label{sect:graphs}
    
    A graph $\mathcal{G}$ is defined as a set on nodes $\mathcal{V}$ interconnected by a set of edges $\mathcal{E}$. The edges of $\mathcal{E}$ are \emph{non-oriented} edges. However, for the purpose of this paper, we will also need to consider their \emph{oriented} counterparts, denoted by $\overrightarrow{\mathcal{E}}$. The size of $\overrightarrow{\mathcal{E}}$ is twice the size of $\mathcal{E}$:  each edge $e \in \mathcal{E}$ is associated to two oriented edges in $\overrightarrow{\mathcal{E}}$. Choosing arbitrarily an orientation for $e$ and writing $s_e$ and $t_e$ its \emph{source} and \emph{target}, we have $e = (s_e,t_e) \in \overrightarrow{\mathcal{E}}$, and its reversely-oriented edge $e^* = (t_e,s_e) \in \overrightarrow{\mathcal{E}}$. We state our results for \emph{unweighted} graphs, \emph{but all our results generalize to the weighted case}, as detailed in the Supplementary Material.
    
    \subsection{Combinatorial Laplacian}
    
    An elementary instance of graph-supported data is that of real values attached to the nodes of the graph, usually formalized as a vector $f~\in~\R^\mathcal{V}$, and called a \emph{graph signal}~\cite{shuman2013emerging}. The regularity of such a signal can be quantified using the \emph{combinatorial Laplacian} $\mathsf{L}~:~\R^\mathcal{V}~\rightarrow~\R^\mathcal{V}$, a symmetric semi-definite positive operator, conveniently expressed as $\mathsf{L} = \mathsf{D} - \mathsf{A}$, where $\mathsf{D}$ is the diagonal degree matrix ($\mathsf{D}_{i,i} = d_i$ the degree of node $i$) and $\mathsf{A}$ the adjacency matrix of the graph (with $\mathsf{A}_{(i,j)} = 1$ if $(i,j) \in \overrightarrow{\mathcal{E}}$ and $0$ otherwise)~\cite{chung1997spectral}. The quadratic form associated to $\mathsf{L}$ acts as:
    \begin{equation}
        \langle f, \mathsf{L} f \rangle = \frac{1}{2} \sum_{e \in \overrightarrow{\mathcal{E}}} \vert f(t_e) - f(s_e) \vert^2,
        \label{eq:action_laplacian}
    \end{equation}
    which associates a high value to functions with important local variations over the edges of the graph, and serves as a basis for the notion of \emph{frequency} for graph signals~\cite{shuman2013emerging,tremblay_design_2018}. These operators are also related to combinatorial properties of the graph $\mathcal{G}$~\cite{bollobas1998modern}, for instance:   
    \begin{enumerate}[label={(P\arabic*)}]
        \item $\ker{\mathsf{L}}$ is always non-empty, with dimension the number of connected components of $\mathcal{G}$. \label{enum:1}
        \item $\mathsf{D}^{-1} \mathsf{A}$ is the matrix of the natural \emph{random walk} on $\mathcal{G}$, which transitions from a node to one of its neighbors with uniform probability at each step, with $\left(\left(\mathsf{D}^{-1} \mathsf{A}\right)^l\right)_{i,j}$ the probability that a path of length $l$ starting at $i$ ends up in node $j$. \label{enum:2}
    \end{enumerate}
    
    \subsection{Connection Laplacian: Definition and Basic Properties}
    \label{sect:graphs:connections}

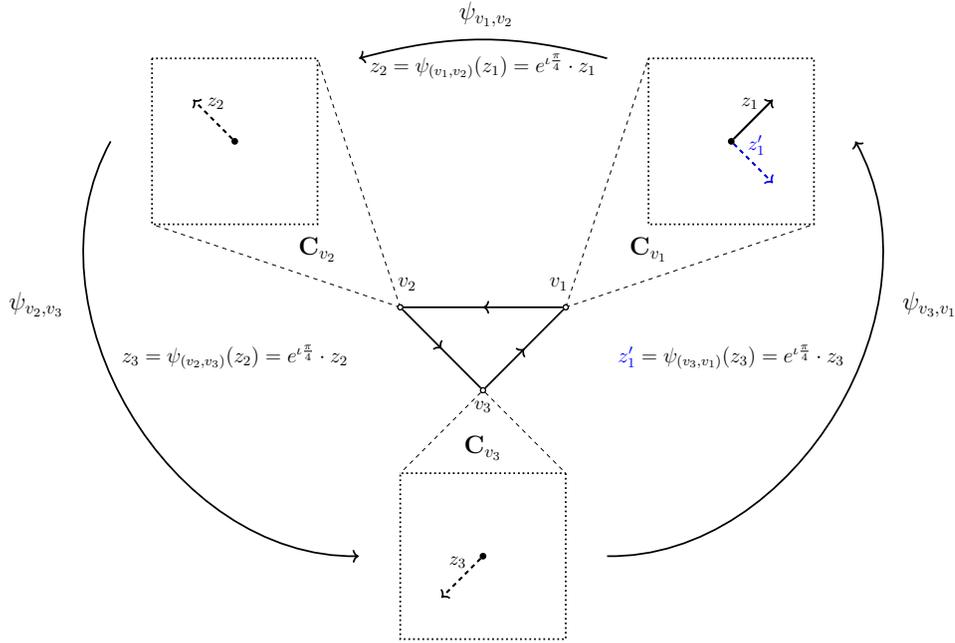
\begin{figure}[ht]
        \centering
        \scalebox{0.55}{\begin{tikzpicture}[every node/.style={circle,inner sep=1pt,thick,draw,scale=1.2},
	                        every edge/.style={very thick,draw},every edge quotes/.append style = {auto, font=\footnotesize, inner sep=2pt}
]

            \node (2) at (0,0){};
            \node (1) at (4,0){};
            \node (3) at (2,-2){};

            \begin{scope}[very thick,decoration={
                            markings,
                            mark=at position 0.5 with {\arrow{>}}}] 
                \draw[postaction={decorate}] (1) -- (2);
                \draw[postaction={decorate}] (2) -- (3);
                \draw[postaction={decorate}] (3) -- (1);
            \end{scope}

            \node[very thick,label={\large $z_2 = \psi_{(v_1,v_2)}(z_1) =  e^{\iota \frac{\pi}{4}} \cdot z_1$},scale=0.001,color=white] (4) at (2,3){};

            \node[very thick,label={\large $z_3 = \psi_{(v_2,v_3)}({z_2}) = e^{\iota \frac{\pi}{4}} \cdot {z_2}$},scale=0.001,color=white] (4) at (-4,-4){};

            \node[very thick,label={\large $\textcolor{blue}{z_1'} = \psi_{(v_3,v_1)}(z_3) = e^{\iota \frac{\pi}{4}} \cdot z_3$},scale=0.001,color=white] (4) at (8,-4){};
            
            \draw[dashed] (4,0) -- (6,6);
            \draw[dashed] (4,0) -- (10,2);
            \draw[very thick,dotted] (6,6) rectangle (10,2);
            \filldraw (6,2) circle (0pt) node[scale=0.001,label={below: \Large $\mathbf{C}_{v_1}$}]{};
            \node[very thick,fill] (c1) at (8,4){};
            \draw[very thick] (c1) -> (9,5) node[midway,label={above:$z_1 \ $},scale=0.001]{};
            \path (c1) edge[->] (9,5);
            \draw[very thick,dashed] (c1) -> (9,3) node[midway,label={above: $\large \ \textcolor{blue}{z_1'}$},scale=0.001]{};
            \path (c1) edge[color=blue,dashed,->] (9,3);

            \draw[dashed] (0,0) -- (-6,2);
            \draw[dashed] (0,0) -- (-2,6);
            \draw[very thick,dotted] (-6,6) rectangle (-2,2);
            \filldraw (-2,2) circle (0pt) node[scale=0.001,label={below: \Large $\mathbf{C}_{v_2}$}]{};
            \node[very thick,fill] (c2) at (-4,4){};
            \draw[very thick,dashed] (c2) -> (-5,5) node[midway,label={above: $
            \ {z_2}$},scale=0.001]{};
            \path (c2) edge[dashed,->] (-5,5);

            \draw[dashed] (3) -- (4,-4);
            \draw[dashed] (3) -- (0,-4);
            \draw[very thick,dotted] (4,-4) rectangle (0,-8);
            \draw[ultra thin,dotted] (0,-4) -- (4,-4) node[scale=0.001,midway,label={above:\Large $\mathbf{C}_{v_3}$}]{};
            \node[very thick,fill] (c3) at (2,-6){};
            \draw[very thick,dashed] (c3) -> (1,-7) node[midway,label={above:$z_3 \ $},scale=0.001]{};
            \path (c3) edge[dashed,->] (1,-7);

            \node[label={above: \large $ \ \ v_2$},fill=white] (2) at (0,0){};
            \node[label={above: \large $v_1 \ \ $},fill=white] (1) at (4,0){};
            \node[label={below: $v_3$},fill=white] (3) at (2,-2){};

            \draw[very thick] (-7,4) edge[->,bend right=60]  (-1,-6);
            \node[label=\Large $\psi_{v_2,v_3}\ \ \ $,scale=0.0001,color=white] (4) at (-8.5,-1){};

            \draw[very thick] (5,-6) edge[->,bend right=60]  (11,4);
            \node[label=\Large $\ \ \ \psi_{v_3,v_1} $,scale=0.0001,color=white] (4) at (12.5,-1){};

            \draw[very thick] (-1,6) edge[bend left=15,<-] (5,6);
            \node[label={\vspace{0.1cm} \Large $\psi_{v_1,v_2}$},scale=0.0001,color=white] (4) at (2,6.2){};

    \end{tikzpicture}}
        \caption{A connection $\Psi$ on the triangle graph $K_3$ (in the center). The connection maps $\psi_e$ are associated with angles $\theta_{(v_1,v_2)} = \theta_{(v_2,v_3)} = \theta_{(v_3,v_1)} = \frac{\pi}{4}$, and we represent their action on some $z_1 \in \C_{v_1}$. A cyclic path $C = ((v_1,v_2),(v_2,v_3),(v_3,v_1))$ is depicted using directed arrows, we denote by $\psi_C = \psi_{(v_3,v_1)} \circ \psi_{(v_2,v_3)} \circ \psi_{(v_1,v_2)}$ the composition of the connection maps along this cycle (a rotation of $\theta_C = \frac{3 \pi}{4})$.
        Top right: $z_1 \in \C_{v_1}$ (bold arrow) and $z_1' = \psi_C(z_1) \in \C_{v_1}$ (blue dotted arrow). Top left, bottom: $z_2 = \psi_{(v_1,v_2)}(z_2) \in \C_{v_2}$ and $z_3 = \psi_{(v_2,v_3)}(z_2) \in \C_{v_3}$.}
        \label{fig:connection}
    \end{figure}

    \emph{Unitary connections} are a practical way to add additional information to the graph's structure\footnote{By describing explicit geometrical transformations between the signal values of the nodes of $\mathcal{G}$.}. We describe in the following how to capture additional rotations associated to edges in the form of a Laplacian-like operator. The main idea is to represent rotations as \emph{unitary complex numbers}, which will naturally lead to consider \emph{complex-valued} functions $f \in \C^\mathcal{V}$ when studying variational properties on $\mathcal{G}$. The entries of $f$ should be understood as belonging to different copies $\C_v$ of the complex plane $\C$ associated to each node $v \in \mathcal{V}$, see Figure~\ref{fig:connection}\footnote{This association of a copy of $\C$ is analogous to a fiber bundle over a manifold (\emph{e.g.} its tangent bundle), and is known as a \emph{complex line bundle} over $\mathcal{G}$.}.

    We will associate to each directed edge in $\overrightarrow{\mathcal{E}}$ a map representing the transformation along that edge, which is known as a \emph{connection}~\cite{kenyon_spanning_2011}. A \emph{unitary connection} $\Psi$ on a graph $\mathcal{G}$ is a collection of unitary linear maps $(\psi_e)_{e \in \overrightarrow{\mathcal{E}}}$, with each map $\psi_e~:~\C_{s_e}\rightarrow\C_{t_e}$ acting as multiplication by a unitary complex number $e^{\iota \theta_e}$, so that $\psi_e(z) = e^{\iota \theta_e} \cdot z$, with $\iota$ the complex imaginary unit (we sometimes abuse notation and write $\psi_e = e^{\iota \theta_e}$ when there is no risk of confusion). In addition, we require that $\psi_{e^*} = \psi_e^*$ (\emph{i.e.} $\psi_{(t_e,s_e)} = \psi_{(s_e,t_e)}^*$), so that the transformation associated to an edge traversed in one direction or the other differs only by conjugation. This notion is illustrated in Figure~\ref{fig:connection}. 
    
    \begin{rem}
    \label{rem:gamma}
    In applications, such as angular synchronization, the $\theta_e$'s often represent \emph{a priori} known rotations, and are frequently modulated by a scale parameter $\gamma \in \R_+$, resulting in connections such that $\psi_e (z) = e^{\iota \gamma \theta_e} z$. While we will forego this additional parameter in our theoretical results, it is easily recovered by replacing each $\theta_e$ with $\gamma \theta_e$.
    \end{rem}

    \begin{paragraph}{The connection Laplacian} To a graph endowed with a unitary connection, we associate the \emph{connection Laplacian} $\mathsf{L}_\theta:~\C^\mathcal{V}\rightarrow~\C^\mathcal{V}$, defined  by $\mathsf{L}_\theta = \mathsf{D} - \mathsf{A}_\theta$, with $\mathsf{A}_\theta$ a connection-aware adjacency matrix such that $(\mathsf{A}_\theta)_{i,j} = e^{-\iota \theta_{(i,j)}}$ if $\{i,j\} \in \mathcal{E}$, and $\mathsf{A}_{i,j} = 0$ otherwise. This operator generalizes the Laplacian to non-trivial connections. \nt{Its quadratic form reads:}
    \begin{equation}
        \langle f, \mathsf{L}_\theta f \rangle = \frac{1}{2} \sum_{e \in \overrightarrow{\mathcal{E}}} \vert f(t_e) - e^{\iota \theta_{e}} f(s_e) \vert^2.
        \label{eq:connection_laplacian_action}
    \end{equation}
    For the trivial connection, with $\psi_e = \id_{\C_{s_e},\C_{t_e}}$ (\emph{i.e.} $\theta_e = 0$) for all edges $e \in \overrightarrow{\mathcal{E}}$, we recover $\mathsf{L}_\theta = \mathsf{L}$ and the expression in Eq.~\eqref{eq:action_laplacian}. Unlike this specific case, $\langle f, \mathsf{L}_\theta f \rangle$ also penalizes functions that are incoherent with respect to the connection, including constant functions.
    \end{paragraph}

    Let us now answer a natural question: when is the kernel $\ker \mathsf{L}_\theta$ non-empty, and what are the functions $f \in \C^\mathcal{V}$ that are not penalized by $\langle f , \mathsf{L}_\theta f \rangle$?

    \begin{paragraph}{($\mathrm{P1'}$) Angular synchronization and $\ker{L_\theta}$}
        When the entries $f_i$ of the vector $f$ are \emph{unitary} complex numbers $f_i = e^{\iota s_i}$, we recover from Eq.~\eqref{eq:connection_laplacian_action} the expression in Eq.~\eqref{eq:cosine_loss}. In fact, one obtains $\vert f(j) - e^{\iota \theta_{(i,j)}} f(i) \vert^2 = 2 - 2 \mathfrak{Re}((e^{\iota s_{j}})^* e^{\iota \theta_{(i,j)}} e^{\iota s_{i}})$, \nt{which yields:}
        \begin{equation}
            \langle f, \mathsf{L}_\theta f \rangle = \sum_{\{i,j\} \in E} 2 - 2 \cos((s_j - s_i) - \theta_{(i,j)}).
        \end{equation}
        Even though this equality only holds for $f \in U(\C)^\mathcal{V}$, note that $\langle f, \mathsf{L}_\theta f \rangle = 0$ if and only if $f_i = e^{\iota \theta_{(j,i)}} f_j$ for all $(i,j) \in \overrightarrow{\mathcal{E}}$, that is, if $s_i = s_j + \theta_{(j,i)}$.
        As a consequence, $\mathsf{L}_\theta$ is not only semi-definite positive (as seen from Eq.~\eqref{eq:connection_laplacian_action}), but \emph{positive} definite, \emph{unless} there is an \emph{exact solution} $x \in U(\C)^\mathcal{V}$ to the associated angular synchronization problem, in which case the kernel $\ker \mathsf{L}_\theta$ is generated by $x$.
    \end{paragraph}

    This behavior of the smallest eigenpair of $\mathsf{L}_\theta$ is in clear contrast to that of $\mathsf{L}$ (which is \emph{always} zero, see property~\ref{enum:1}), and we explain in Section~\ref{sect:smoothing} how a generalization of property~\ref{enum:2} pertaining to \emph{random propagations} can be used to solve connection-Laplacian-based  problems.
    

    \section{Random Estimators for Tikhonov Smoothing under a Connection}
    \label{sect:smoothing}
    
   \nt{ This section shows how to {smooth} a complex signal on a graph via {random propagations}. Let $\mathcal{G} = (\mathcal{V},\mathcal{E})$ be a graph endowed with a connection $\Psi$.  We aim at smoothing a signal $g \in \C^\mathcal{V}$ by solving:}
    \begin{equation}
        f_*=\argmin_{f \in \C^\mathcal{V}} \ q \Vert f - g \Vert_2^2 + \langle f, \mathsf{L}_\theta f \rangle,
        \label{eq:tikhonov}
    \end{equation}
    with $q \in \R_+^*$. This is nothing but a re-writing of Eq.~\eqref{eq:tikhonov_intro}. \nt{Note that $\langle f , \mathsf{L}_\theta f \rangle$ penalizes functions that are not coherent with respect to the connection $\Psi$.} 
    \begin{prop}
        The solution to Problem~\eqref{eq:tikhonov} can be expressed as:
        \begin{equation}
            f_* = q(\mathsf{L}_\theta + q \mathsf{I})^{-1} g.
            \label{eq:solution_tikhonov}
        \end{equation}
        \label{prop:tikhonov_solution}
    \end{prop}
    \nt{The proof, based on standard $\mathbf{CR}$-calculus, is in Section~\ref{supmat:feynman_kac} of the Supplementary Material.} 
    Let us now describe how \nt{to estimate $f_*$} via local propagation over random paths on $\mathcal{G}$.

    \subsection{Local Estimation: a Feynman-Kac Formula}
    \label{sect:smoothing:fk}
    
        Feynman-Kac formulas express solutions of variational problems as the expectations of stochastic processes: we develop here an instance tailored to the Tikhonov smoothing problem of Eq.~\eqref{eq:tikhonov}. We start with a few definitions.
        A \emph{path} $p$ in $\mathcal{G}$ is an ordered sequence of oriented edges in $\overrightarrow{\mathcal{E}}$. We consider the two operations of concatenation $pq$ of two paths $p$ and $q$ and of orientation-reversal $p^*$ of a path $p$, and denote by $P_i^j$ the set of paths from $i$ to $j$ in $\mathcal{G}$. Connection maps extend to paths inductively as $\psi_{pq} = \psi_q \circ \psi_p$, and we will for instance frequently encounter the map
        $\psi_{p^*} : \C_{j} \rightarrow \C_{i}$, acting by composition of the (inverse) rotations encountered along path $p \in P_i^j$ from $i$ to $j$.
        
        The stochastic process we consider is a random walk with a non-zero probability of being interrupted at any node. It is conveniently defined on an extended graph $\mathcal{G}^\Gamma$, with nodes $\mathcal{V}~\cup~\{\Gamma\}$ and edges $\mathcal{E} \cup \bigcup_{v \in \mathcal{V}} \{ \{ v,\Gamma \} \}$ (so that all the nodes of $\mathcal{V}$ are connected to $\Gamma$), where $\Gamma$ is an additional \emph{boundary} node. We then define the random walk $(u_t)_{t \geq 0}$, which starts at some fixed node $u_0 = i \in V$ and transitions to some other node at time $t$ as follows: 
        \begin{itemize}
            \item $u_{t+1} = \Gamma$ with probability $\frac{q}{d_{u_t} + q}$ \nt{(recall that $d_{u_t}$ stands for the degree of node $u_t$)}, 
            \item $u_{t+1} = v \in \mathcal{V}$ with probability $\frac{\mathsf{A}_{u_t,v}}{d_{u_t} + q}$.
        \end{itemize}
        The process ends upon reaching node $\Gamma$. In other words, this is the natural random walk on $\mathcal{G}_\Gamma$, with a stopping criterion at the boundary $\Gamma$. Writing $u_l = \Gamma$ for some time $l$, this process results in a random path $((u_0,u_1),...,(u_{l-1},u_l))$. We denote by $P_i^\Gamma$ the set of all possible paths obtainable via this process, and by $\nu_i$ the resulting probability measure on $P_i^\Gamma$. 
        \\
        We can then show that propagating along paths sampled according to $\nu_i$ converges in expectation to the solution of Problem~\eqref{eq:tikhonov} on node $i$. More precisely:
        
        \begin{prop}
            Denoting by $j = u_{l-1}$ the last node reached before $p$ reaches $\Gamma$, we have:
            \begin{equation}
            f_*(i) = \E_{p \sim \nu_i}\left( \psi_{p^{\;*}_\Gamma}(g_j) \right),
            \label{eq:feynmac_kac}
            \end{equation}
            where $p = p_\Gamma e_\Gamma$ with $e_\Gamma$ the last edge in $p$.
            \label{prop:feynmac_kac}
        \end{prop}
        
        \noindent In other words, an unbiased estimate of $f_*(i)$ is obtained by drawing a random path from $i$ to $\Gamma$, and retropropagating the value of $g$ at the node just before $\Gamma$ ($j$ in Proposition~\ref{prop:feynmac_kac}), taking into account the rotations along the path. See Alg.~\ref{alg:feynman-kac}. This formula is a discrete-time analog of a Theorem in~\cite{kassel2021covariant}, and we generalize it to weighted graphs in the Supplementary Material (Section~\ref{supmat:feynman_kac}, Proposition~\ref{prop:feynman_kac_bis}), building upon a connection-aware version of property~\ref{enum:2}.
        
        \begin{algorithm}[t]
        \caption{Feynman-Kac estimator (Proposition~\ref{prop:feynmac_kac}).}
        \begin{algorithmic}[1]
            \State $f_i \leftarrow 0$ 
            \State{\textbf{Repeat} $m$ times}
            \State{\hspace{\algorithmicindent}Sample a path $p \in P^\Gamma_i$ from $\nu_i$} \Comment{By running an interrupted random walk on $\mathcal{G}$}
            \State{\hspace{\algorithmicindent}$f_i \leftarrow f_i + \psi_{p^*_\Gamma}(g_j)$} \Comment{$j$ the last node before interruption of $p$}
            \State \textbf{Output} $\frac{1}{m} f_i$ \Comment{\nt{An unbiased estimator of $f_*(i)$}}
        \end{algorithmic} 
        \label{alg:feynman-kac}
        \end{algorithm}

    Alg.~\ref{alg:feynman-kac} is simple but may be computationally expensive: estimation on one node requires sampling $m$ paths, which needs to be repeated \emph{for each} of the $\vert \mathcal{V} \vert$ nodes in $\mathcal{G}$ (a locality issue inherent to all  Walk-on-Spheres-type algorithms \cite{muller1956some}).
    We address these limitations in the remainder of this section and propose novel propagation-based estimators of $f_*$, \nt{allowing to update the Monte-Carlo estimates on \emph{all the nodes at once}, by i)~sampling and ii)~propagating the signal, on specific subgraphs: the so-called MTSFs.}
    
    \begin{rem}
        \label{rem:higher}
        The restriction to complex signals on $\mathcal{G}$ and unitary connections (encoding 2D rotations) is mostly artificial, and Proposition~\ref{prop:feynmac_kac} generalizes to transformations in \emph{e.g.} $O(\R^d)$. This setting is encountered for instance for 3D molecule alignement in cryogenic electron microscopy (cryo-EM)~\cite{singer2012vector}, or in structure-from-motion problems~\cite{arie2012global}, where the $\psi_e$'s act as 3D rotations (\emph{i.e.} $\psi_e \in SO(\R^3))$. 
    \end{rem}

    \subsection{Multi-Type Spanning Forests: a first global estimator}
    \label{sect:smoothing:mtsf}

\begin{figure}[t]
    \centering
        
        \begin{subfigure}{0.49\textwidth}
            \centering
            \scalebox{0.98}{
\begin{tikzpicture}[every node/.style={circle,inner sep=1pt,thick,draw}, every edge/.style={very thick,draw}]
	        \node (1) at (0,0){};
            \node (2) [right of=1]{};
            \node (3) [right of=2]{};
            \node (4) [right of=3]{};
            \node (5) [below of=1]{};
            \node (6) [right of=5]{};
            \node (7) [right of=6]{};
            \node (8) [right of=7,color=purple,fill,very thick,inner sep=2pt,label={right:$r$}]{};
            \node (9) [below of=5]{};
            \node (10) [right of=9]{};
            \node (11) [right of=10]{};
            \node (12) [right of=11]{};
            \node (13) [below of=9]{};
            \node (14) [right of=13]{};
            \node (15) [right of=14]{};
            \node (16) [right of=15]{};
            \path [-,color=purple] (1) edge (2);
            \path [-,color=black] (2) edge[dashed] (3);
            \path [-,color=purple] (3) edge (4);
            \path [-,color=purple] (5) edge (6);
            \path [-,color=black] (6) edge[dashed] (7);
            \path [-,color=purple] (7) edge (8);
            \path [-,color=black] (9) edge[dashed] (10);
            \path [-,color=purple] (10) edge (11);
            \path [-,color=black] (11) edge[dashed] (12);
            \path [-,color=purple] (13) edge (14);
            \path [-,color=purple] (14) edge (15);
            \path [-,color=black] (15) edge[dashed] (16);
            \path [-,color=purple] (1) edge (5);
            \path [-,color=purple] (2) edge (6);
            \path [-,color=black] (3) edge[dashed] (7);
            \path [-,color=purple] (4) edge (8);
            \path [-,color=purple] (5) edge (9);
            \path [-,color=black] (6) edge[dashed] (10);
            \path [-,color=black] (7) edge[dashed] (11);
            \path [-,color=purple] (8) edge (12);
            \path [-,color=black] (9) edge[dashed] (13);
            \path [-,color=purple] (10) edge (14);
            \path [-,color=purple] (11) edge (15);
            \path [-,color=purple] (12) edge (16);
            
            \node[color=white,scale=0.0001] (20) at (0,-(3 + 1/4){};
        \end{tikzpicture}}
            \caption{A MTSF $\phi$ in purple, with two unicycles (top left, bottom) and one tree (right) with its root $r$ (large purple node).}
        \end{subfigure}
        \hfill
        \begin{subfigure}{0.49\textwidth}
            \centering
            \scalebox{0.98}{
\begin{tikzpicture}[every node/.style={circle,inner sep=1pt,thick,draw},
	                        every edge/.style={very thick,draw}]
         
	        \node[cross out,very thick,color=teal,line width=2.8pt,inner sep=2pt] (1) at (0,0){};
            \node[color=white,scale=0.0001,label={above:\footnotesize $\tilde{f}(j,\phi,g) = 0$}] (20) at (0,-1/2){};
            
            \node[cross out,very thick,color=teal,line width=2.8pt,inner sep=2pt] (2) [right of=1]{};
            \node (3) [right of=2]{};
            \node[label={right:\footnotesize $\ \ \tilde{f}(i,\phi,g) = \psi_{(r,i)}(g(r))$}] (4) [right of=3]{};
            \node[cross out,very thick,color=teal,line width=2.8pt,inner sep=2pt] (5) [below of=1]{};
            \node[cross out,very thick,color=teal,line width=2.8pt,inner sep=2pt] (6) [right of=5]{};
            \node (7) [right of=6]{};
            \node (8) [right of=7]{};
            \node[cross out,very thick,color=teal,line width=2.8pt,inner sep=2pt] (9) [below of=5]{};
            \node[cross out,very thick,color=teal,line width=2.8pt,inner sep=2pt] (10) [right of=9]{};
            \node[cross out,very thick,color=teal,line width=2.8pt,inner sep=2pt] (11) [right of=10]{};
            \node (12) [right of=11]{};
            \node[cross out,very thick,color=teal,line width=2.8pt,inner sep=2pt] (13) [below of=9]{};
            \node[cross out,very thick,color=teal,line width=2.8pt,inner sep=2pt] (14) [right of=13]{};
            \node[cross out,very thick,color=teal,line width=2.8pt,inner sep=2pt] (15) [right of=14]{};
            \node (16) [right of=15]{};   

            \node (8) [right of=7,color=blue,fill,very thick,inner sep=2pt]{};
            \path [-,color=black] (2) edge[dashed] (3);
            \path [-,color=black] (6) edge[dashed] (7);
            \path [-,color=black] (9) edge[dashed] (10);
            \path [-,color=black] (11) edge[dashed] (12);
            \path [-,color=black] (15) edge[dashed] (16);
            \path [-,color=black] (3) edge[dashed] (7);
            \path [-,color=black] (6) edge[dashed] (10);
            \path [-,color=black] (7) edge[dashed] (11);
            \path [-,color=black] (9) edge[dashed] (13);
            
            \node (8) [right of=7,color=blue,fill,very thick,inner sep=2pt,label={right:\footnotesize $\ \ \tilde{f}(r,\phi,g) = g(r)$}]{};
            \path [-,color=black] (2) edge[dashed] (3);
            \path [-,color=black] (6) edge[dashed] (7);
            \path [-,color=black] (9) edge[dashed] (10);
            \path [-,color=black] (11) edge[dashed] (12);
            \path [-,color=black] (15) edge[dashed] (16);
            \path [-,color=black] (3) edge[dashed] (7);
            \path [-,color=black] (6) edge[dashed] (10);
            \path [-,color=black] (7) edge[dashed] (11);
            \path [-,color=black] (9) edge[dashed] (13);

            \node (8) [right of=7,color=blue,fill,very thick,inner sep=2pt]{};
            \path [<-,color=blue] (4) edge (8);\node (4) [right of=3,color=blue,fill]{};
            \path [->,color=blue] (4) edge (3);\node (3) [right of=2,color=blue,fill]{};
            \path [<-,color=blue] (7) edge (8);\node (7) [right of=6,color=blue,fill]{};
            \path [<-,color=blue] (12) edge (8);\node (12) [right of=11,color=blue,fill]{};
            \path [->,color=blue] (12) edge (16);\node (16) [right of=15,color=blue,fill]{};

            \node (1) at (0,0)[color=teal,fill]{};
            \node (2) [right of=1,color=teal,fill]{};
            \node (5) [below of=1,color=teal,fill]{};
            \node (6) [right of=5,color=teal,fill]{};
            \node (9) [below of=5,color=teal,fill]{};
            \node (10) [right of=9,color=teal,fill]{};
            \node (11) [right of=10,color=teal,fill]{};
            \node (13) [below of=9,color=teal,fill]{};
            \node (14) [right of=13,color=teal,fill]{};
            \node (15) [right of=14,color=teal,fill]{};
            
            \path [-,color=teal] (1) edge (2);
            \path [-,color=teal] (5) edge (6);
            \path [-,color=teal] (10) edge (11);
            \path [-,color=teal] (13) edge (14);
            \path [-,color=teal] (14) edge (15);
            \path [-,color=teal] (1) edge (5);
            \path [-,color=teal] (2) edge (6);
            \path [-,color=teal] (5) edge (9);
            \path [-,color=teal] (10) edge (14);
            \path [-,color=teal] (11) edge (15);
            
        \end{tikzpicture}}
            \caption{Estimation: we set $\widetilde{f}(i,\phi,g) = \psi_{r \xrightarrow[]{\phi} i} (g_r)$ for nodes in the tree (blue arrows), and $\widetilde{f}(j,\phi,g) = 0$ on unicycles (teal crosses).}
        \end{subfigure}
        
        \caption{Theorem~\ref{th:fermionic} illustrated on the $4 \times 4$ grid graph, on one MTSF.}
        \label{fig:ex_th_bis}
    \end{figure}
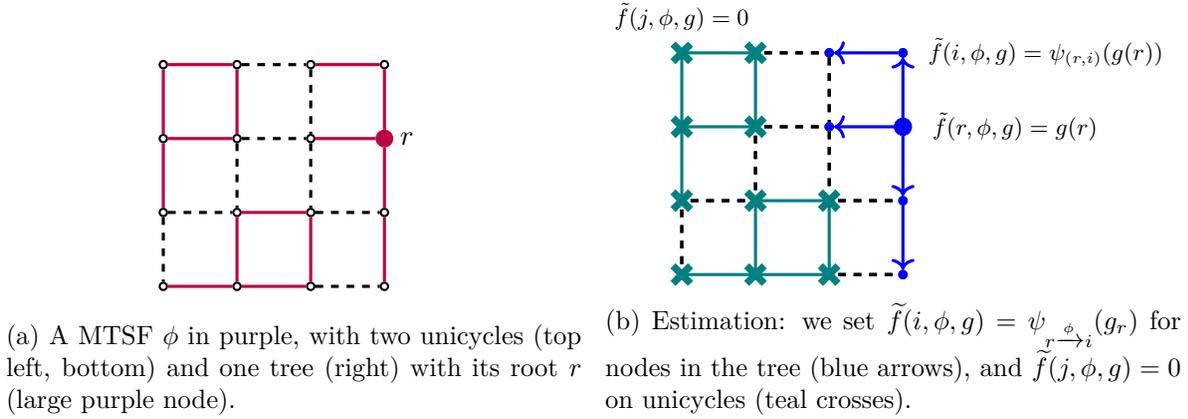

    To build our global estimators, we rely on decompositions of the graph into \emph{Multi-Type Spanning Forests} (MTSFs)~\cite{fanuel_sparsification_2022}. A MTSF $\phi \subseteq \mathcal{V} \cup \mathcal{E}$ decomposes $\mathcal{G}$ into disjoint {components} that are either \emph{rooted trees}, or \emph{unicycles}. More precisely, a MTSF $\phi$ must have fixed cardinality $\vert \phi \vert = \vert \mathcal{V} \vert$, and is divided into (maximal) components $c_\phi \subseteq \phi$, which must be either:
    \begin{itemize}
        \item a rooted tree $c_\phi \subseteq \mathcal{V} \cup \mathcal{E}$, such that $c_\phi \cap \mathcal{E}$ is a connected cycle-free subset of edges, and $c_\phi \cap \mathcal{V}$ is reduced to a single node, \nt{called \emph{root},} connected to $c_\phi \cap \mathcal{E}$;
        \item a unicycle $c_\phi \subseteq \mathcal{E}$, which is a connected subset of edges containing \emph{exactly} one cycle.
    \end{itemize}
    See Figure~\ref{fig:ex_th_bis} for an illustration. 
    Note that in the absence of unicycles, a MTSF is a \emph{rooted spanning forest} \cite{avena2013random,avena2018two}, whereas a tree-free MTSF is a \emph{spanning forest of unicycles} \cite{kenyon_spanning_2011}.

 	\nt{Denote by $\mathcal{M}(\mathcal{G})$ the set of all MTSFs over  $\mathcal{G}$. Consider the distribution $\mathcal{D}_\mathcal{M}$ over $\mathcal{M}(\mathcal{G})$:}
    \begin{equation}
        \prb_{\mathcal{D}_\mathcal{M}}(\phi) \propto q^{\vert \phi \cap \mathcal{V} \vert} \prod_{C \in \mathcal{C}(\phi)} \left( 2 - 2 \cos(\theta_C) \right),
        \label{eq:prob_dm}
    \end{equation}
    where $\mathcal{C}(\phi)$ is the set of cycles belonging to the unicycles of $\phi$, and $\theta_C$ is the argument of the unitary complex number associated to the connection map $\psi_C$ obtained from a path traversing $C$ \emph{one time} (as in Fig.~\ref{fig:connection}). Note that, while $\psi_C$ depends on the path's orientation, $\cos(\theta_C)$ is insensitive to this orientation. Also, a cycle has a non-zero probability of being sampled if and only if it is \emph{inconsistent} (if $C$ is perfectly consistent: $\theta_C = 0$ and $2 - 2\cos(\theta_C) = 0$). 
    $\mathcal{D}_\mathcal{M}$ is a rooted variant of the distribution considered in \cite{fanuel_sparsification_2022}, and is a DPP\footnote{So is the distribution of \cite{fanuel_sparsification_2022}, but their non-combinatorial argument does not carry over to our rooted process. See the Section~\ref{supmat:fermionic} of the Supplementary Material for a proof in this case \nt{*proof of what?*}} \nt{over the set $\mathcal{V}\cup\mathcal{E}$}.
    We can now state our main theorem.
    \begin{thm}
        Let $\phi$ be a MTSF sampled according to $\mathcal{D}_\mathcal{M}$. Denote by $a \xrightarrow[]{\phi} b$ the unique path from $a$ to $b$ in some tree of a MTSF $\phi$ and, \emph{if $v \in \mathcal{V}$} belongs to a tree of $\phi$, by $r_\phi(v)$ the root of this tree.
        Consider the estimator \nt{
        \begin{equation}
        	\widetilde{f}(i,\phi,g) 
        	\begin{cases}
        		= \psi_{r_\phi(i) \xrightarrow[]{\phi} i} \left( g(r_\phi(i)) \right) & \text{if $i$ belongs to a tree} \\
        		= 0 & \text{if $i$ lies in a unicycle.}
        	\end{cases}
        \end{equation}
        When $i$ belongs to a tree, this estimator propagates the value of the signal $g$ from the root $r_\phi(i)$ to $i$, via the connections.} When $i$ belongs to a unicycle, it is simply zero. 
        Then, one has: 
        \begin{equation}
        	f_*(i) =\mathbf{E}_{\mathcal{D}_\mathcal{M}}(\widetilde{f}(i,\phi,g)),
     \label{eq:fermionic}
    \end{equation}
{\emph{i.e.}, the estimator is unbiased.} 
        \color{col}{Moreover, the variance of $\widetilde{f}$ can be characterized as\footnote{Here, the variance of a complex-valued random variable $\widetilde{z}$ is defined as $\mathrm{var}(\widetilde{z}) = \mathbf{E}(\vert \widetilde{z} \vert^2 ) - \left\vert \mathbf{E}(\widetilde{z})^2 \right\vert$, and we compute in Eq.~\eqref{eq:variance_tilde} the expected squared error $\mathbf{E}_{\mathcal{D}_\mathcal{M}}\left( \Vert \widetilde{f}(\phi,g) - f_* \Vert^2 \right) = \sum_{i \in \mathcal{V}} \mathrm{var}\left( \widetilde{f}(i,\phi,g) \right)$.}
        \begin{equation}
            \mathbf{E}_{\mathcal{D}_\mathcal{M}}\left( \Vert \widetilde{f}(\phi,g) - f_* \Vert^2 \right) = \langle g , (\mathsf{I} - \mathsf{K}^2) g   \rangle,
            \label{eq:variance_tilde}
        \end{equation}}
		where $\mathsf{K} = q(\mathsf{L}_\theta + q \mathsf{I})^{-1}$ denotes the connection-aware Tikhonov smoothing operator.
        \label{th:fermionic}
    \end{thm}
    The proof of Theorem~\ref{th:fermionic} relies on the reformulation of $\mathcal{D}_\mathcal{M}$ as a DPP, and an extension of arguments coming from both \cite{kenyon_spanning_2011} and \cite{pilavci_graph_2021}. It is a special case of Theorem~\ref{th:fermionic_bis}, which includes weighted graphs, and is proved in Section~\ref{supmat:fermionic} of the Supplementary Material. 

    In a less technical phrasing, propagating for each tree in a MTSF $\phi$ the value of $g$ at its root $r$ to its other nodes $i$ (along path $r \xrightarrow[]{\phi} i$) results in an unbiased estimator of $f_*$. The estimation on nodes belonging to unicycles is then $0$. See Alg.~\ref{alg:mtsf_propagation} and Figure~\ref{fig:ex_th_bis}.

    \begin{algorithm}[ht]
        \caption{MTSF-based estimator (Theorem~\ref{th:fermionic}).}
        \begin{algorithmic}[1]
            \State $f_i \leftarrow 0 \ \forall i \in \mathcal{V}$
            \State{\textbf{Repeat} $m$ times}
            \Indent
            \State{Sample $\phi \in \mathcal{M}(\mathcal{G})$ from $\mathcal{D}_\mathcal{M}$} \Comment{Alg.~\ref{alg:sampling_mtsf} can be used}
            \For{$i \in \mathcal{V}$}
            \If{$i$ belongs to a rooted tree of $\phi$}
            \State $f_i \leftarrow f_i + \psi_{r_\phi(i) \xrightarrow[]{\phi} i} \left( g(r_\phi(i)) \right)$
            \EndIf
            \EndFor
            \EndIndent
            \State \textbf{Output} $\frac{1}{m} (f_1,...,f_{\vert \mathcal{V} \vert})$
        \end{algorithmic} 
        \label{alg:mtsf_propagation}
    \end{algorithm}
    
    \begin{rem}
        If the connection $\Psi$ is trivial (\emph{i.e.}, $\psi_e = \id_{\C_{s_e},\C_{t_e}}$ for all edges), we recover the result from \cite{pilavci_graph_2021}. There are no inconsistent cycles in this case, and $\mathcal{D}_\mathcal{M}$ reduces to another determinantal distribution $\mathcal{D}_\mathcal{F}$ over \emph{rooted spanning forests} $\phi \in \mathcal{V} \cup \mathcal{E}$ \cite{avena2013random,avena2018two}.
        \label{rem:yigit}
    \end{rem}

    {Unlike Alg.~\ref{alg:feynman-kac}, Alg.~\ref{alg:mtsf_propagation}  allows to update the estimates on \emph{all the nodes} for each sampled MTSF, in $\mathcal{O}(\vert \mathcal{V} \vert)$ time. On the other hand, it becomes necessary to sample a MTSF $\phi$ according to $\mathcal{D}_\mathcal{M}$. A first naive strategy is to use algorithms designed to sample generic DPPs. They however rely on the eigendecomposition of a kernel matrix of size $(\vert V \vert + \vert E \vert) \times (\vert V \vert + \vert E \vert)$, which is prohibitive. In the next section, we show how to take advantage of a very efficient random-walk-based algorithm (from~\cite{fanuel_sparsification_2022}), and analyze its expected computational cost.}
    
    \subsection{Sampling Multi-Type Spanning Forests}
    \label{sect:smoothing:mtsf_sampling}

    Efficient sampling from $\mathcal{D}_\mathcal{M}$ is possible via \emph{loop-erased random walks}, traditionally used in Wilson's algorithm to sample uniform spanning trees~\cite{wilson_random_1996}. \nt{The simplest scenario arises when the following condition is verified}\footnote{Condition~\ref{cond:weak_inconsistency} has already been identified as a technical sampling condition for Wilson-like algorithms (in \cite{kassel_random_2017} and \cite{fanuel_sparsification_2022}). It also appears in a different guise in \cite{kassel2021covariant}, under the name of \emph{trace-positivity}, where it allows to simplify some technical statements.}.
    \begin{cond}[Weak-Inconsistency]
        For all cycles $C$ in $\mathcal{G}$, $\cos(\theta_C) \geq 0$.
        \label{cond:weak_inconsistency}
    \end{cond}
    We stress that this \nt{condition constrains} \emph{all} the cycles of the graph, and not just the cycles $\mathcal{C}(\phi)$ in some MTSF $\phi \in \mathcal{M}(G)$.
    If Condition~\ref{cond:weak_inconsistency} is satisfied, $1 -\cos(\theta_C) \leq 1$ defines a probability measure over the oriented cycles in $G$, which can be leveraged to efficiently sample an MTSF $\phi$, as we detail in the following. If this condition does not hold, one can still efficiently sample MTSFs, by relying on importance sampling, as explained later in this section.


    \begin{paragraph}{MTSF sampling algorithm} 

        The algorithm of \cite{fanuel_sparsification_2022}, recalled here as Alg.~\ref{alg:sampling_mtsf}, samples a MTSF $\phi$ according to $\mathcal{D}_\mathcal{M}$ by simulating multiple random walks on $\mathcal{G}^\Gamma$, and constructs $\phi$ iteratively from these random paths. The precise way in which a path is turned into a part of $\phi$ relies on a \emph{loop-erasure} procedure.
        The full sampling scheme is detailed in Alg.~\ref{alg:sampling_mtsf}, performing random walks that can be stopped in a number of ways: by being interrupted (\emph{i.e.} reaching $\Gamma$), by building a cycle $\theta_C$ (kept with probability $1 - \cos(\theta_C)$), or by reaching a node already spanned by $\phi$\footnote{Some implementation details of Alg.~\ref{alg:sampling_mtsf} are not made completely explicit here, such as the precise way in which we detect cycles (instruction at line 11), or how to track the cycle-acceptance probabilities (\emph{i.e.} $\theta_C$). Please refer to Section~\ref{supmat:complexity} of the Supplementary Material and our Julia implementation for more details}. We suppose that the nodes of $\mathcal{G}$ are arbitrarily ordered in a queue, and say that a node $i$ is spanned by $\phi$ if $i \in \phi$ or if $i$ is an endpoint of some edge $e \in \phi$. We also use a function $\texttt{random$\_$successor}(u)$ that, at each call, randomly outputs either $\Gamma$ with probability $\frac{q}{d_u + q}$, or some node $v$ with probability  $\frac{\mathsf{A}_{u,v}}{d_u + q}$.
        
        \begin{algorithm}[t]
        \caption{MTSF sampling algorithm~\cite{fanuel_sparsification_2022}.}
        \begin{algorithmic}[1]
            \State $\phi \leftarrow \emptyset$
            \While{$\phi$ not spanning}
                \State Let $i \in \mathcal{V}$ be the first node in the queue not spanned by $\phi$
                \State $u \leftarrow i$ \Comment{$u$ is the current node of the random walk}
                \State $p \leftarrow \epsilon$
                \Comment{$\epsilon$ the empty path}
                \While{($u \neq \Gamma$) and ($p$ does not intersect $\phi$ or contain a cycle)}
                    \State $u' \leftarrow \texttt{random$\_$successor}(u)$ \Comment{Move to next node}
                    \If{$u' \neq \Gamma$}
                    \State $e \leftarrow (u,u')$, $p \leftarrow p e$ \Comment{Add $e$ to the path $p$}
                    \EndIf
                    \If{${p}$ contains a cycle $C$}
                        \State Remove this cycle from $p$ with probability $\cos(\theta_C)$
                    \EndIf
                    \State $u_\mathrm{old} \leftarrow u$, $u \leftarrow u'$ \Comment{$u_\mathrm{old}$ the previous node}
                \EndWhile
                \If{$u = \Gamma$}
                    \State $\phi \leftarrow \phi \cup p \cup u_\mathrm{old}$ \Comment{Add the sampled path $p$ and the root $u_\mathrm{old}$ to $\phi$\footnotemark}
                \Else
                    \State $\phi \leftarrow \phi \cup p$
                \EndIf
            \EndWhile
            \State \textbf{Output} $\phi$
        \end{algorithmic}        \label{alg:sampling_mtsf}
        \end{algorithm}
    \end{paragraph}
    One then proves the following by applying the arguments of~\cite{fanuel_sparsification_2022} (which describes how to sample \emph{unrooted} MTSFs), and keeping track of the roots.
    \begin{prop}
        Suppose that Condition~\ref{cond:weak_inconsistency} holds. Then, Alg.~\ref{alg:sampling_mtsf} outputs a MTSF $\phi$ distributed according to $\mathcal{D}_\mathcal{M}$.
    \end{prop}
    Let us now discuss the cost of Alg.~\ref{alg:sampling_mtsf}.
        \begin{prop}
            The expected number of calls to \texttt{random$\_$successor} (\emph{i.e.,} the expected number of random walk steps), in Alg.~\ref{alg:sampling_mtsf}, is upper-bounded by
            \begin{equation}
            \mathrm{tr}\left( (\mathsf{L} + q \mathsf{I})^{-1} (\mathsf{D} + q \mathsf{I}) \right).
            \label{eq:complexity}
            \end{equation}
            Also, Alg.~\ref{alg:sampling_mtsf} can be implemented with an expected time complexity upper-bounded by $\mathcal{O}\left( \frac{\vert \mathcal{E} \vert}{q} \right)$.
            \label{prop:bound_complexity}
        \end{prop}
        See Section~\ref{supmat:complexity} of the Supplementary Material for the proof. The bound is obtained by noting that Alg.~\ref{alg:sampling_mtsf} necessitates \emph{fewer} steps than the Wilson-like algorithm in~\cite{pilavci_graph_2021}, used to sample random spanning forests~\cite{avena2013random,avena2018two}, with expected number of steps exactly given by $\eqref{eq:complexity}$.
        Unlike the algorithm in~\cite{pilavci_graph_2021} though, Alg.~\ref{alg:sampling_mtsf} requires both tracking of the angular offsets accumulated along a cycle, and detecting said cycle, at additional computational cost. We discuss two possible implementations in Section~\ref{supmat:complexity} of the Supplementary Material, one of them resulting in the $\mathcal{O}\left( \frac{\vert \mathcal{E} \vert}{q} \right)$ expected \emph{time complexity} mentioned in Proposition~\ref{prop:bound_complexity}. This translates to $\mathcal{O}\left( \frac{\vert \mathcal{E} \vert}{q} + \vert \mathcal{V} \vert \right)$ expected runtime for  Alg.~\ref{alg:mtsf_propagation}.\footnotetext{{We abuse notations and denote here by $p$ the set of \emph{non-oriented} edges in the path.}}
    \begin{paragraph}{When Condition~\ref{cond:weak_inconsistency} does not hold: Importance Sampling}
    \label{rem:importance_sampling}
        \color{col}{In case Condition~\ref{cond:weak_inconsistency} is not satisfied, one can still estimate $f_*$ using an importance sampling strategy\footnote{The importance weights we consider have been proposed in~\cite{fanuel_sparsification_2022}, in the context of spectral sparsification.}. Specifically, one can threshold the incoherence $2 - 2\cos(\theta_C)$ and sample from
        \begin{equation}
            \prb_{IS}(\phi) \propto q^{\vert \phi \cap \mathcal{V} \vert} \prod_{C \in \mathcal{C}(\phi)} \min(2,2 - 2 \cos(\theta_C))
            \label{eq:importance}
        \end{equation}
        using a variant of Alg.~\ref{alg:sampling_mtsf}, where a cycle is systematically kept, at line 12, if $\cos(\theta_C) \leq 0$. The estimation then uses the importance weights
        \begin{equation}
            w(\phi) = \prod_{C \in \mathcal{C}(\phi)} \max(1,1-\cos(\theta_C)),
        \end{equation}
        which allows to estimate $cf_*$ for some unknown constant $c > 0$. \emph{Note that in the angular-synchronization pipeline we consider in Section~\ref{sect:synchro}, this is sufficient as we only need an estimation of $f_*$ up to a multiplicative constant}. For smoothing purposes, one can further consider a self-normalized importance-sampling estimator~\cite{owen2013monte}: for  $m$ MTSFs $\{\phi_k\}_{k=1}^m$ sampled according to the distribution of Eq.~\eqref{eq:importance}, we set
    \begin{equation}
        \widetilde{f}_{IS} = \frac{1}{\sum_{k=1}^m w(\phi_k)} \sum_{k = 1}^m \widetilde{f}(\phi_k,g), 
    \end{equation}
    which tends to $f_*$ as $m \rightarrow \infty$. This result holds \emph{regardless of the incoherence of the connection}, but the estimator is only guaranteed to be unbiased in the limit.}
    \end{paragraph}

    Both the Feynman-Kac-based Alg.~\ref{alg:feynman-kac} and the MTSF-based Alg.~\ref{alg:mtsf_propagation} can be roughly described as performing a random walk on $\mathcal{G}$ before stopping at some root node, and then retro-propagating the value from this root to the starting point of the random walk. The main difference is that MTSFs allow to update the estimation on all nodes jointly. \nt{We will now see how our estimator $\tilde{f}$ can be improved by implementing two variance reduction techniques.}
    
    \subsection{Variance Reduction for the MTSF-based estimator}
    \label{par:var_red}
    A paramount property of Monte-Carlo estimators is not only their unbiased behavior in expectation, but also their variance: they are significantly improved when used in conjunction with efficient variance reduction techniques.
    We propose two such improvements over $\widetilde{f}$, based on the classical approaches of \emph{Rao-Blackwellization}~\cite{blackwell1947conditional,rao1992information} and \emph{control variates}~\cite{kroese2013handbook} respectively, generalizing the variance-reduction techniques introduced in \cite{pilavci_graph_2021,pilavci2022variance}.

    \begin{paragraph}{Rao-Blackwellization}
        Rao-Blackwellization leverages the two laws of total expectation and variance, roughly stating that conditioning an estimator using another statistic extracted from the same sample still results in an unbiased estimator, with lower variance (see Section~\ref{supmat:variance} of the Supplementary Material for more details in our case). Here, our technique consists in conditioning on the set of edges in $\phi$. This yields the following estimator:
        \begin{equation}
                \overline{f}(i,\phi,g) = \psi_{r_\phi(i) \xrightarrow[]{\phi} i} \left( h_\phi(r_\phi(i),g) \right) \text{ if $i$ belongs to a rooted tree,}
                \label{eq:rao}
        \end{equation}
        where $h_\phi$ is defined as 
        \begin{equation}
            h_\phi(r,g) = \frac{\sum_{j \in c_\phi(r)} \psi_{j \xrightarrow[]{\phi} r} g(j)}{\vert c_\phi(r) \vert},
        \end{equation} and $c_\phi(r)$ is the \emph{set of nodes} spanned by the tree containing $r$. If $i$ belongs to a unicycle, we once again set $\overline{f}(i,\phi,g) = 0$.
        This results in Alg.~\ref{alg:mtsf_propagation_rb}, which amounts to:
        \begin{itemize}
            \item Computing at the root $r$ of each tree the average $h_\phi(r,g)$ of the values of $g$ over the tree (obtained by propagating from each node $i$ in the tree to $r$).
            \item Propagating this average back to the other nodes of the tree.
        \end{itemize}
        As compared with Alg.~\ref{alg:mtsf_propagation}, this procedure can be implemented at little additional cost\footnote{Note also the similarity to the Belief Propagation algorithm over trees~\cite{mezard2009information}.}.

        \begin{algorithm}[t]
        \caption{Rao-Blackwellized MTSF-based estimator (Theorem~\ref{th:fermionic}).}
        \begin{algorithmic}[1]
            \State $f_i = 0 \ \forall i \in \mathcal{V}$, $h_r = 0 \ \forall r \in \mathcal{V}$
            \State{\textbf{Repeat} $m$ times}
            \Indent
            \State Sample $\phi \in \mathcal{M}(\mathcal{G})$ from $\mathcal{D}_\mathcal{M}$ \Comment{Via Alg.~\ref{alg:sampling_mtsf}}
            \For{$r \in \phi \cap \mathcal{V}$}
            \For{$j \in c_\phi(r)$}
            \State $h_r = h_r + \psi_{j \xrightarrow[]{\phi} r} g(j)$ \Comment{Propagate and average at the root}
            \EndFor
            \EndFor
            \For{$i \in \mathcal{V}$ belonging to a rooted tree of $\phi$}
            \State $f_i = f_i + \psi_{r_\phi(i) \xrightarrow[]{\phi} i} \left( h_{r_\phi(i)} \right)$ \Comment{Propagate back}
            \EndFor
            \EndIndent
            \State \textbf{Output} $\frac{1}{m} (f_1,...,f_{\vert \mathcal{V} \vert})$
        \end{algorithmic} 
        \label{alg:mtsf_propagation_rb}
    \end{algorithm}
    \end{paragraph}

    \begin{paragraph}{Control variates}
        Second is the introduction of control variates: an addition of another random quantity to $\overline{f}$, designed to have zero mean (so that the expectation remains unchanged), but resulting in an estimator with lower variance when designed properly. We propose to use a single \emph{gradient-descent} step with parametrized step-size $\alpha$:
        \begin{equation}
           \hat{f}(\phi,g) = \overline{f}(\phi,g) - \alpha \mathsf{P} (q^{-1} (\mathsf{L}_\theta + q \mathsf{I}) \overline{f}(\phi,g) - g),
           \label{eq:gradient_step}
        \end{equation}
        where $\mathsf{P} = \left( q^{-1} \mathsf{D} + \mathsf{I} \right)^{-1}$ is a diagonal preconditioner for the system $q^{-1} (\mathsf{L}_\theta + q \mathsf{I}) f = g$, and $\alpha \in \R^*_+$.
        A good choice of step-size $\alpha$ is crucial in order to obtain a significant reduction of variance. We simply take $\alpha = 1$ in the following (see Section~\ref{supmat:variance} of the Supplementary Material for empirical results backing this choice).
    \end{paragraph}

    We prove in Section~\ref{supmat:variance} of the Supplementary Material that generalizations of both $\overline{f}$ and $\hat{f}$ are unbiased estimators of $f_*$, along with some results on the concentration of $\overline{f}$.
        
    \begin{prop}
        The variance-reduced estimators are unbiased: \begin{equation}
        \mathbf{E}_{\mathcal{D}_\mathcal{M}}(\overline{f}(i,\phi,g)) = \mathbf{E}_{\mathcal{D}_\mathcal{M}}(\hat{f}(i,\phi,g)) = f_*(i).
        \end{equation}
        \color{col}{Further, the variance of the Rao-Blackwellized estimator is given by
        \begin{equation}
            \mathbf{E}_{\phi \sim \mathcal{D}_\mathcal{M}}\left( \Vert \overline{f}(\phi,g) - f_* \Vert^2 \right) = \langle g, (\mathsf{K} - \mathsf{K}^2) g  \rangle,
        \end{equation}
        where $\mathsf{K} = q(\mathsf{L}_\theta + q \mathsf{I})^{-1}$, and we have the following finite-sample concentration bound. Let $\eps,\delta \in (0,1)$, and consider sampling $m$ MTSFs $\{\phi_k\}_{k=1}^m$. Then, provided that\footnote{The factor $6$ in Eq.~\eqref{eq:m_condition} is not optimal: it is obtained from a more explicit bound appearing in the proof.}
        \begin{equation}
            m \geq \frac{6}{\eps^2} \log\left( \frac{\vert \mathcal{V} \vert}{\delta} \right),
        \end{equation}
        it holds that
        \begin{equation}
            \prb\left(\forall g \in \mathcal{C}^\mathcal{V},\quad \left\Vert \left( \frac{1}{m} \sum_{k = 1}^m \overline{f}(\phi_k,g) \right) - f_* \right\Vert_2 \leq \eps \Vert g \Vert_2 \right) \geq 1 - \delta.
            \label{prop:detail_tropp}
        \end{equation}}
        \label{prop:variance_reduction}
    \end{prop}

    Both variance-reduction techniques are easy to implement, and do not incur large additional computational costs: the Rao-Blackwellization $\overline{f}$ comes at $\mathcal{O}(\vert \mathcal{V} \vert)$ additional cost, and the gradient-descent step in $\hat{f}$ entails \emph{a single} matrix-vector multiplication, in $\mathcal{O}(\vert \mathcal{E} \vert )$ time. \textcolor{col}{Further, Eq.~\eqref{eq:concentration} shows that the number of forests required to reach a fixed precision grows \emph{logarithmically} with the size of the graph, and that it is \emph{insensitive to the choice of $q$ and to the incoherence of the graph} (which both affect the conditioning of the system, and the performance of conjugate-gradient solvers\footnote{See the \emph{Discussion} paragraph of Section~\ref{subsect:runtime_precision_smoothing} for a reminder on the expected performance of CG}). Note though that these parameters \emph{do} impact the sampling cost for MTSFs.}

    \section{Numerical Results under Weak-Inconsistency}
    \label{sect:numerical}

    We now analyze the behavior of the estimator proposed in Theorem~\ref{th:fermionic}, along with the improved versions discussed in Section~\ref{par:var_red}, and compare their performance with (deterministic)  conjugate-gradient-based solvers~\cite{saad2003iterative}, on graphs with different topologies\footnote{We perform all our measurements using \emph{a single thread} on a laptop with an intel i7-1185G7 processor.}. We perform experiments on the following connection model.

    \begin{paragraph}{Connection Model}
    \label{par:connection_model}
    For a given graph $\mathcal{G}$ with $n$ nodes, we associate to each node $v$ an angle $\omega_v$ chosen uniformly at random in $[0,2\pi)$, and we set 
    \begin{equation}
        \theta_e = \omega_{t_e} - \omega_{s_e} + \eta \eps_e
        \label{eq:degradation}
    \end{equation}
    for all edges $e$, with $\eps_e$ a perturbation uniformly distributed in $[-1,1]$, and $\eta \in \R_+^*$ a scaling constant. We set  $\eta = \frac{\pi}{2n}$ to ensure that Condition~\ref{cond:weak_inconsistency} is satisfied (so that  Alg.~\ref{alg:sampling_mtsf} in fact samples a MTSF according to $\mathcal{D}_\mathcal{M}$).
    \end{paragraph}

    \subsection{A First Runtime Experiment on Erdös-Rényi Graphs}
    \label{sect:numerical:density}

    We first analyze the behavior of our estimator $\widetilde{f}$ with respect to two parameters: the choice of regularization parameter $q$ and the mean degree $\overline{d}$ of the graph (controlling the graph's density). In our first experiment, we let $\overline{d}$ take values in $\{50,100,150,200\}$, and $q$ take values equal to $q' \overline{d}$, with $q' \in \{10^{-3},10^{-1},1\}$. The complexity bound of Eq.~\eqref{eq:complexity} becomes \emph{linear in $\vert \mathcal{V} \vert$} when using such a parametrization.

\begin{figure}[t]
        \centering
        \scalebox{0.4}{

\begin{tikzpicture}[/tikz/background rectangle/.style={fill={rgb,1:red,1.0;green,1.0;blue,1.0}, fill opacity={1.0}, draw opacity={1.0}}, show background rectangle]
\begin{axis}[point meta max={nan}, point meta min={nan}, legend cell align={left}, legend columns={1}, title={}, title style={at={{(0.5,1)}}, anchor={south}, font={{\fontsize{14 pt}{18.2 pt}\selectfont}}, color={rgb,1:red,0.0;green,0.0;blue,0.0}, draw opacity={1.0}, rotate={0.0}, align={center}}, legend style={color={rgb,1:red,0.0;green,0.0;blue,0.0}, draw opacity={1.0}, line width={1}, solid, fill={rgb,1:red,1.0;green,1.0;blue,1.0}, fill opacity={1.0}, text opacity={1.0}, font={{\fontsize{24 pt}{10.4 pt}\selectfont}}, text={rgb,1:red,0.0;green,0.0;blue,0.0}, cells={anchor={center}}, at={(1.02, 1)}, anchor={north west}}, axis background/.style={fill={rgb,1:red,1.0;green,1.0;blue,1.0}, opacity={1.0}}, anchor={north west}, xshift={1.0mm}, yshift={-1.0mm}, width={145.4mm}, height={99.6mm}, scaled x ticks={false}, xlabel={Mean degree}, x tick style={color={rgb,1:red,0.0;green,0.0;blue,0.0}, opacity={1.0}}, x tick label style={color={rgb,1:red,0.0;green,0.0;blue,0.0}, opacity={1.0}, rotate={0}}, xlabel style={at={(ticklabel cs:0.5)}, anchor=near ticklabel, at={{(ticklabel cs:0.5)}}, anchor={near ticklabel}, font={{\fontsize{24 pt}{14.3 pt}\selectfont}}, color={rgb,1:red,0.0;green,0.0;blue,0.0}, draw opacity={1.0}, rotate={0.0}}, xmajorgrids={true}, xmin={45.5}, xmax={204.5}, xticklabels={{$50$,$100$,$150$,$200$}}, xtick={{50.0,100.0,150.0,200.0}}, xtick align={inside}, xticklabel style={font={{\fontsize{24 pt}{10.4 pt}\selectfont}}, color={rgb,1:red,0.0;green,0.0;blue,0.0}, draw opacity={1.0}, rotate={0.0}}, x grid style={color={rgb,1:red,0.0;green,0.0;blue,0.0}, draw opacity={0.1}, line width={0.5}, solid}, axis x line*={left}, x axis line style={color={rgb,1:red,0.0;green,0.0;blue,0.0}, draw opacity={1.0}, line width={1}, solid}, scaled y ticks={false}, ylabel={Time in ms}, y tick style={color={rgb,1:red,0.0;green,0.0;blue,0.0}, opacity={1.0}}, y tick label style={color={rgb,1:red,0.0;green,0.0;blue,0.0}, opacity={1.0}, rotate={0}}, ylabel style={at={(ticklabel cs:0.5)}, anchor=near ticklabel, at={{(ticklabel cs:0.5)}}, anchor={near ticklabel}, font={{\fontsize{24 pt}{14.3 pt}\selectfont}}, color={rgb,1:red,0.0;green,0.0;blue,0.0}, draw opacity={1.0}, rotate={0.0}}, ymajorgrids={true}, ymin={0.7490744253997002}, ymax={4.135511650350301}, yticklabels={{$1$,$2$,$3$,$4$}}, ytick={{1.0,2.0,3.0,4.0}}, ytick align={inside}, yticklabel style={font={{\fontsize{24 pt}{10.4 pt}\selectfont}}, color={rgb,1:red,0.0;green,0.0;blue,0.0}, draw opacity={1.0}, rotate={0.0}}, y grid style={color={rgb,1:red,0.0;green,0.0;blue,0.0}, draw opacity={0.1}, line width={0.5}, solid}, axis y line*={left}, y axis line style={color={rgb,1:red,0.0;green,0.0;blue,0.0}, draw opacity={1.0}, line width={1}, solid}, colorbar={false}]
    \addplot[color={rgb,1:red,0.0;green,0.6056;blue,0.9787}, name path={9786c281-1979-44be-9ce8-9ed3198ad1df}, draw opacity={1.0}, line width={2}, dotted, mark={*}, mark size={3.0 pt}, mark repeat={1}, mark options={color={rgb,1:red,0.0;green,0.6056;blue,0.9787}, draw opacity={1.0}, fill={rgb,1:red,0.0;green,0.6056;blue,0.9787}, fill opacity={1.0}, line width={0.75}, rotate={0}, solid}]
        table[row sep={\\}]
        {
            \\
            50.0  2.4319999259599996  \\
            100.0  2.5600396339000002  \\
            150.0  2.61139626255  \\
            200.0  2.91035164684  \\
        }
        ;
    \addlegendentry {MTSF, $q' = 10^{-3}$}
    \addplot[color={rgb,1:red,0.0;green,0.6056;blue,0.9787}, name path={8f7fd9d9-74bc-4805-b319-f12409d82fd5}, draw opacity={1.0}, line width={2}, dashed, mark={*}, mark size={3.0 pt}, mark repeat={1}, mark options={color={rgb,1:red,0.0;green,0.6056;blue,0.9787}, draw opacity={1.0}, fill={rgb,1:red,0.0;green,0.6056;blue,0.9787}, fill opacity={1.0}, line width={0.75}, rotate={0}, solid}]
        table[row sep={\\}]
        {
            \\
            50.0  2.25428107808  \\
            100.0  2.34928887878  \\
            150.0  2.6937036674100003  \\
            200.0  2.5295756079500005  \\
        }
        ;
    \addlegendentry {MTSF, $q' = 10^{-1}$}
    \addplot[color={rgb,1:red,0.0;green,0.6056;blue,0.9787}, name path={2aae0a22-45f3-43bf-b818-f9a80093a32a}, draw opacity={1.0}, line width={2}, solid, mark={*}, mark size={3.0 pt}, mark repeat={1}, mark options={color={rgb,1:red,0.0;green,0.6056;blue,0.9787}, draw opacity={1.0}, fill={rgb,1:red,0.0;green,0.6056;blue,0.9787}, fill opacity={1.0}, line width={0.75}, rotate={0}, solid}]
        table[row sep={\\}]
        {
            \\
            50.0  1.73808417128  \\
            100.0  1.82200154652  \\
            150.0  1.9336663219999999  \\
            200.0  1.6747914454000001  \\
        }
        ;
    \addlegendentry {MTSF, $q' = 1$}
    \addplot[color={rgb,1:red,0.2422;green,0.6433;blue,0.3044}, name path={c409a1b5-2f5e-470f-a454-b61c2fcfb49e}, draw opacity={1.0}, line width={2}, dotted, mark={diamond*}, mark size={3.0 pt}, mark repeat={1}, mark options={color={rgb,1:red,0.2422;green,0.6433;blue,0.3044}, draw opacity={1.0}, fill={rgb,1:red,0.2422;green,0.6433;blue,0.3044}, fill opacity={1.0}, line width={0.0}, rotate={0}, solid}]
        table[row sep={\\}]
        {
            \\
            50.0  0.8965733432899999  \\
            100.0  1.80619606815  \\
            150.0  2.66841359052  \\
            200.0  3.9768822379900004  \\
        }
        ;
    \addlegendentry {$\mathsf{L}_\theta f$, $q' = 10^{-3}$}
    \addplot[color={rgb,1:red,0.2422;green,0.6433;blue,0.3044}, name path={c94a4cc6-9123-4adb-8570-53fd3a887cd4}, draw opacity={1.0}, line width={2}, dashed, mark={diamond*}, mark size={3.0 pt}, mark repeat={1}, mark options={color={rgb,1:red,0.2422;green,0.6433;blue,0.3044}, draw opacity={1.0}, fill={rgb,1:red,0.2422;green,0.6433;blue,0.3044}, fill opacity={1.0}, line width={0.75}, rotate={0}, solid}]
        table[row sep={\\}]
        {
            \\
            50.0  0.8449169883700001  \\
            100.0  1.75095113286  \\
            150.0  2.6715416945999997  \\
            200.0  3.95557313557  \\
        }
        ;
    \addlegendentry {$\mathsf{L}_\theta f$, $q' = 10^{-1}$}
    \addplot[color={rgb,1:red,0.2422;green,0.6433;blue,0.3044}, name path={1542edaa-72f9-4d97-9fbb-a3fa9d399e4a}, draw opacity={1.0}, line width={2}, solid, mark={diamond*}, mark size={3.0 pt}, mark repeat={1}, mark options={color={rgb,1:red,0.2422;green,0.6433;blue,0.3044}, draw opacity={1.0}, fill={rgb,1:red,0.2422;green,0.6433;blue,0.3044}, fill opacity={1.0}, line width={0.75}, rotate={0}, solid}]
        table[row sep={\\}]
        {
            \\
            50.0  0.8864710292  \\
            100.0  1.84298632954  \\
            150.0  3.01535530646  \\
            200.0  4.03966908738  \\
        }
        ;
    \addlegendentry {$\mathsf{L}_\theta f$, $q' = 1$}
\end{axis}
\end{tikzpicture}}
        \caption{Runtime when varying $\overline{d}$, for different values of $q$.}
        \label{fig:var_deg}
    \end{figure}
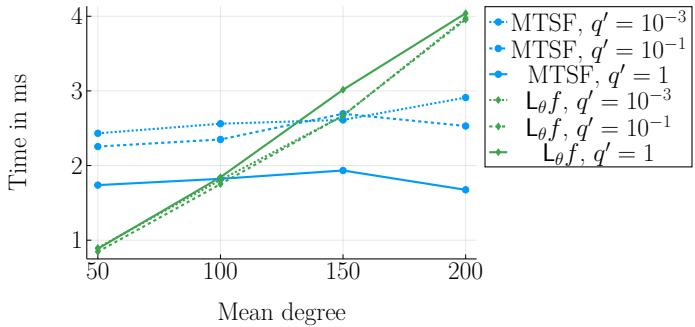

    \begin{paragraph}{Setup}
    We generate Erdös-Rényi random graphs $\mathcal{G} \sim \mathrm{ER}(n,p)$ of size $n = 10^4$ for varying $p \in [0,1]$, so that each edge $e$ independently appears with probability $p$. 
    To control the density, we set $p = \frac{\overline{d}}{n-1}$ so that the expected mean degree of these random graphs is $\overline{d}$.
    For each such graph, we generate a random signal $f \in \C^\mathcal{V}$ with independent complex Gaussian entries $f_v \sim \mathcal{N}_\C(0,1)$.
    \\
    We then measure the running time of sampling \emph{one} MTSF $\phi$ and applying the estimator $\widetilde{f}$. 
    As a reference, we also measure in the same manner the runtime of the matrix-vector multiplication $\mathsf{L}_\theta f$, where $\mathsf{L}_\theta$ is implemented as a sparse Hermitian matrix in CSC format. Note that this operation is the most expensive part of an iteration of the Conjugate-Gradient algorithm, and serves as a simple baseline to compare computation costs.
    \\
    The results, depicted in Figure~\ref{fig:var_deg}, are the average over $10^4$ time measurements, themselves averaged over $10$ realizations of the graph $G$.
    \end{paragraph}

    \begin{paragraph}{Analysis}
        Two trends emerge in Figure~\ref{fig:var_deg}. First, as $q'$ (hence $q$) increases, computing $\widetilde{f}$ becomes less expensive, since the probability of the random walk stopping at some root node increases. Second, in contrast to the matrix-vector product $\mathsf{L}_\theta f$, our estimator is less sensitive to the density of the graph: computing $\widetilde{f}$ is more expensive than computing $\mathsf{L}_\theta f$ for $\overline{d}=50$, and systematically faster for $\overline{d} = 200$. This reflects the $\mathcal{O}(\vert \mathcal{V} \vert)$ expected complexity for $\widetilde{f}$ (obtained with the parametrization $q = q' \overline{d}$)\footnote{Results are similar when varying the size $n$ of the graph (not shown).}. 
    \end{paragraph}

    \subsection{Runtime-Precision Trade-offs in the context of Complex Graph Signal Denoising}
    \label{subsect:runtime_precision_smoothing}

    We compare the performance of our improved estimators $\overline{f}$ and $\hat{f}$ with  conjugate-gradient methods~\cite{saad2003iterative}. The objective is to recover a signal $f_\top \in \C^\mathcal{V}$ given a noisy degradation $g = f_\top + \eps$.

    In the specific instance where $f_\top \in \R^\mathcal{V}$ and the connection on $\mathcal{G}$ is trivial ($\psi_e = \mathrm{id}_{\C_{s_e},\C_{t_e}}$ for all $e \in \mathcal{E}$), a common assumption in the graph signal processing literature is that $f_\top$ is a \emph{smooth} signal, that is, a linear combination of the first (low-frequency) eigenvectors of $\mathsf{L}$ (associated to the smallest eigenvalues)~\cite{puy2018random,lorenzo2018sampling}.

    We here similarly assume that $f_\top \in \C^\mathcal{V}$ is $B$-bandlimited (\emph{i.e.} $f_\top \in \mathrm{span}(u_1,...,u_B)$, with $u_i$ the $i$-th eigenvector of $\mathsf{L}_\theta$), so that solving the Tikhonov smoothing Problem~\eqref{eq:tikhonov}:
    \begin{equation}
        \argmin_{f \in \C^\mathcal{V}} \ q \Vert f - g \Vert_2^2 + \langle f, \mathsf{L}_\theta f \rangle,
    \end{equation} 
    should allow to faithfully recover $f_\top$ from $g$, by penalizing high-frequency components in the optimal solution $f_* = q (\mathsf{L}_\theta + q \mathsf{I})^{-1} g$. We refer the reader to Section~\ref{supmat:tikhonov_numerical} of the Supplementary Material for additional supporting arguments in this connection-aware setting.
    
    Recall that this matrix inversion can  be computed using a Cholesky decomposition (for an exact solution), or a conjugate-gradient-based iteration (for a high quality approximation). We compare our estimators with these methods, and consider two cost functions:
    \begin{itemize}
        \item The reconstruction error $e_r(f) = \Vert f - f_\top \Vert_2 / n$, measuring the quality of the denoising.
        \item The approximation error $e_a(f) = \Vert f - f_* \Vert_2 / n$, measuring the quality of the approximation of $f_*$.
    \end{itemize}

    Our estimators and the conjugate-gradient algorithms are respectively parametrized by the number of MTSFs used and the number of gradient steps, that we both denote by $m$.

\begin{figure}[]
        \centering
        
        \begin{subfigure}{1\textwidth}
        \scalebox{0.42}{

\begin{tikzpicture}[/tikz/background rectangle/.style={fill={rgb,1:red,1.0;green,1.0;blue,1.0}, fill opacity={1.0}, draw opacity={1.0}}, show background rectangle]
\begin{axis}[point meta max={nan}, point meta min={nan}, legend cell align={left}, legend columns={1}, title={}, title style={at={{(0.5,1)}}, anchor={south}, font={{\fontsize{14 pt}{18.2 pt}\selectfont}}, color={rgb,1:red,0.0;green,0.0;blue,0.0}, draw opacity={1.0}, rotate={0.0}, align={center}}, legend style={color={rgb,1:red,0.0;green,0.0;blue,0.0}, draw opacity={1.0}, line width={1}, solid, fill={rgb,1:red,1.0;green,1.0;blue,1.0}, fill opacity={1.0}, text opacity={1.0}, font={{\fontsize{18 pt}{10.4 pt}\selectfont}}, text={rgb,1:red,0.0;green,0.0;blue,0.0}, cells={anchor={center}}, at={(0.98, 0.02)}, anchor={south east}}, axis background/.style={fill={rgb,1:red,1.0;green,1.0;blue,1.0}, opacity={1.0}}, anchor={north west}, xshift={1.0mm}, yshift={-1.0mm}, width={150.4mm}, height={99.6mm}, scaled x ticks={false}, xlabel={Time in s}, x tick style={color={rgb,1:red,0.0;green,0.0;blue,0.0}, opacity={1.0}}, x tick label style={color={rgb,1:red,0.0;green,0.0;blue,0.0}, opacity={1.0}, rotate={0}}, xlabel style={at={(ticklabel cs:0.5)}, anchor=near ticklabel, at={{(ticklabel cs:0.5)}}, anchor={near ticklabel}, font={{\fontsize{26 pt}{14.3 pt}\selectfont}}, color={rgb,1:red,0.0;green,0.0;blue,0.0}, draw opacity={1.0}, rotate={0.0}}, xmode={log}, log basis x={10}, xmajorgrids={true}, xmin={0.0006409864611518314}, xmax={0.767198326372942}, xticklabels={{$10^{-3}$,$10^{-2}$,$10^{-1}$}}, xtick={{0.001,0.01,0.1}}, xtick align={inside}, xticklabel style={font={{\fontsize{18 pt}{10.4 pt}\selectfont}}, color={rgb,1:red,0.0;green,0.0;blue,0.0}, draw opacity={1.0}, rotate={0.0}}, x grid style={color={rgb,1:red,0.0;green,0.0;blue,0.0}, draw opacity={0.1}, line width={0.5}, solid}, axis x line*={left}, x axis line style={color={rgb,1:red,0.0;green,0.0;blue,0.0}, draw opacity={1.0}, line width={1}, solid}, scaled y ticks={false}, ylabel={Approximation error}, y tick style={color={rgb,1:red,0.0;green,0.0;blue,0.0}, opacity={1.0}}, y tick label style={color={rgb,1:red,0.0;green,0.0;blue,0.0}, opacity={1.0}, rotate={0}}, ylabel style={at={(ticklabel cs:0.5)}, anchor=near ticklabel, at={{(ticklabel cs:0.5)}}, anchor={near ticklabel}, font={{\fontsize{26 pt}{14.3 pt}\selectfont}}, color={rgb,1:red,0.0;green,0.0;blue,0.0}, draw opacity={1.0}, rotate={0.0}}, ymode={log}, log basis y={10}, ymajorgrids={true}, ymin={1.0e-13}, ymax={8.912509381337459e-5}, yticklabels={{$10^{-10}$,$10^{-5}$}}, ytick={{1.0e-10,1.0e-5}}, ytick align={inside}, yticklabel style={font={{\fontsize{18 pt}{10.4 pt}\selectfont}}, color={rgb,1:red,0.0;green,0.0;blue,0.0}, draw opacity={1.0}, rotate={0.0}}, y grid style={color={rgb,1:red,0.0;green,0.0;blue,0.0}, draw opacity={0.1}, line width={0.5}, solid}, axis y line*={left}, y axis line style={color={rgb,1:red,0.0;green,0.0;blue,0.0}, draw opacity={1.0}, line width={1}, solid}, colorbar={false}]
    \addplot[color={rgb,1:red,0.0;green,0.6056;blue,0.9787}, name path={b3c8ed41-5cc1-49b8-a2fc-dae1b71eab76}, draw opacity={1.0}, line width={2}, solid, mark={*}, mark size={3.0 pt}, mark repeat={1}, mark options={color={rgb,1:red,0.0;green,0.6056;blue,0.9787}, draw opacity={1.0}, fill={rgb,1:red,0.0;green,0.6056;blue,0.9787}, fill opacity={1.0}, line width={0.75}, rotate={0}, solid}]
        table[row sep={\\}]
        {
            \\
            0.0024888612440000005  3.5771575130541316e-5  \\
            0.003938334589999999  2.5236174683473e-5  \\
            0.0055863218514  2.068343433673873e-5  \\
            0.010645435071400002  1.5974659227380132e-5  \\
            0.018653469546600004  1.2596916273156303e-5  \\
            0.027633077271000003  9.901723751359812e-6  \\
            0.036698266087000005  7.594811951292838e-6  \\
            0.06554787482959999  5.953844240362148e-6  \\
            0.10967975443320001  4.608828981604624e-6  \\
            0.19018408853840002  3.5729288623049988e-6  \\
        }
        ;
    \addlegendentry {MTSF}
    \addplot[color={rgb,1:red,0.8889;green,0.4356;blue,0.2781}, name path={2e64e27c-e53b-4fbb-aec9-13cd76b3200f}, draw opacity={1.0}, line width={2}, solid, mark={*}, mark size={3.0 pt}, mark repeat={1}, mark options={color={rgb,1:red,0.8889;green,0.4356;blue,0.2781}, draw opacity={1.0}, fill={rgb,1:red,0.8889;green,0.4356;blue,0.2781}, fill opacity={1.0}, line width={0.75}, rotate={0}, solid}]
        table[row sep={\\}]
        {
            \\
            0.0045992621394  3.549923373715207e-5  \\
            0.007307466370799998  2.5244813893101607e-5  \\
            0.008752054223799998  2.056206731583079e-5  \\
            0.014165958974200005  1.5987519673998344e-5  \\
            0.017723993555800004  1.2611327955454784e-5  \\
            0.022050524307800003  9.89654205515569e-6  \\
            0.04358923002739999  7.6132465876118415e-6  \\
            0.0704586143956  5.955537848716836e-6  \\
            0.11603293021380001  4.608168934999005e-6  \\
            0.1870521408996  3.560096025888562e-6  \\
        }
        ;
    \addlegendentry {MTSF+GS}
    \addplot[color={rgb,1:red,0.2422;green,0.6433;blue,0.3044}, name path={3bbe6e3f-5b20-4130-b38d-4e8677fa1d82}, draw opacity={1.0}, line width={2}, solid, mark={*}, mark size={3.0 pt}, mark repeat={1}, mark options={color={rgb,1:red,0.2422;green,0.6433;blue,0.3044}, draw opacity={1.0}, fill={rgb,1:red,0.2422;green,0.6433;blue,0.3044}, fill opacity={1.0}, line width={0.75}, rotate={0}, solid}]
        table[row sep={\\}]
        {
            \\
            0.0007833642137999998  5.4127684491883725e-5  \\
            0.0018375166840000003  1.0104891157628804e-5  \\
            0.0024759155474000003  4.171839967540631e-6  \\
            0.0038349307125999997  1.7829030605063885e-6  \\
            0.005397654371399998  3.2283448545597917e-7  \\
            0.008196981080199998  1.3997177347482448e-8  \\
            0.016770241755399998  7.950451362984964e-11  \\
            0.021239833911199997  6.180589722191574e-13  \\
            0.02133251967839999  6.180589722191574e-13  \\
            0.019867141064999996  6.180589722191574e-13  \\
        }
        ;
    \addlegendentry {CG}
    \addplot[color={rgb,1:red,0.7644;green,0.4441;blue,0.8243}, name path={bf04491b-eb78-40a2-b838-0c18cdbed0a8}, draw opacity={1.0}, line width={2}, solid, mark={*}, mark size={3.0 pt}, mark repeat={1}, mark options={color={rgb,1:red,0.7644;green,0.4441;blue,0.8243}, draw opacity={1.0}, fill={rgb,1:red,0.7644;green,0.4441;blue,0.8243}, fill opacity={1.0}, line width={0.75}, rotate={0}, solid}]
        table[row sep={\\}]
        {
            \\
            0.0011277352828000001  6.4949911195634e-5  \\
            0.0018518709881999997  2.0023804396564604e-5  \\
            0.003110688280200001  9.341681291850752e-6  \\
            0.0052443270414000015  2.465675943722226e-6  \\
            0.0083624365004  3.3093994838188014e-7  \\
            0.010870971756199999  1.0696168302992578e-8  \\
            0.018749775943999996  2.159442269886342e-11  \\
            0.022108081976799993  7.531177721886595e-13  \\
            0.0228177663518  7.531177721886595e-13  \\
            0.023127968319000002  7.531177721886595e-13  \\
        }
        ;
    \addlegendentry {$\mathsf{D}^{-1}$+CG}
    \addplot[color={rgb,1:red,0.6755;green,0.5557;blue,0.0942}, name path={719cabc3-9936-48b8-bd7a-e955d9d961a9}, draw opacity={1.0}, line width={2}, solid, mark={*}, mark size={3.0 pt}, mark repeat={1}, mark options={color={rgb,1:red,0.6755;green,0.5557;blue,0.0942}, draw opacity={1.0}, fill={rgb,1:red,0.6755;green,0.5557;blue,0.0942}, fill opacity={1.0}, line width={0.75}, rotate={0}, solid}]
        table[row sep={\\}]
        {
            \\
            0.06354264776600002  2.1501552698155275e-5  \\
            0.07948069751699999  4.087577155187896e-6  \\
            0.0740346127344  8.366620887789885e-7  \\
            0.0731329089826  4.9788481281972815e-8  \\
            0.08450445685980003  4.866229581621992e-10  \\
            0.10494994346000001  6.036660239555044e-13  \\
            0.09905871235959998  4.775160804470074e-13  \\
            0.097464296797  4.775160804470074e-13  \\
            0.09907990184759999  4.775160804470074e-13  \\
            0.101761656415  4.775160804470074e-13  \\
        }
        ;
    \addlegendentry {ILU+CG}
    \addplot[color={rgb,1:red,0.0;green,0.6658;blue,0.681}, name path={9bc4e3c8-2fd3-410c-99e4-6197b4afa283}, draw opacity={1.0}, line width={2}, dotted, forget plot]
        table[row sep={\\}]
        {
            \\
            0.6277587507322  1.1220184543019654e-22  \\
            0.6277587507322  79432.82347242805  \\
        }
        ;
\end{axis}
\end{tikzpicture}}
        \hfill
        \scalebox{0.42}{\input{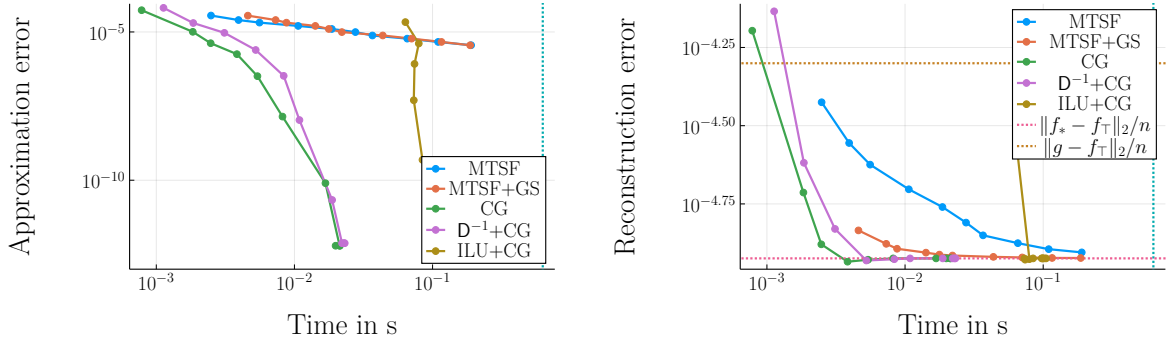}
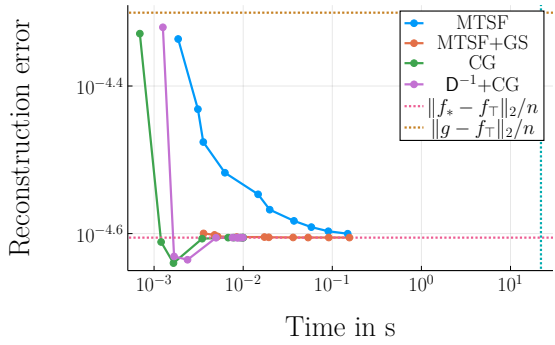
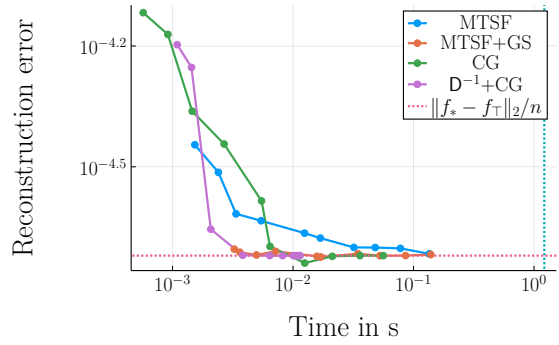
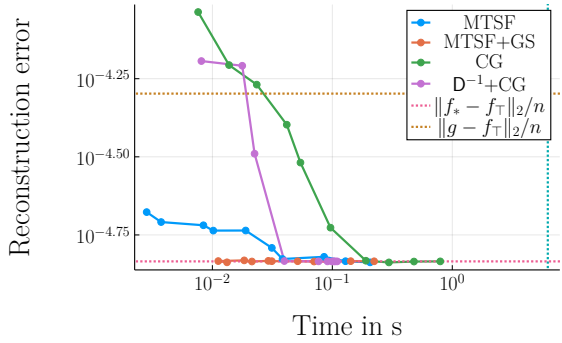
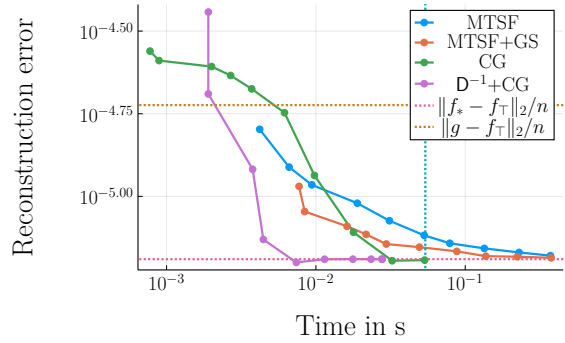}
        \caption{$\eps$-graph. Left: Approximation error. Right: Reconstruction error.}
        \label{subfig:smoothing:eps}
        \end{subfigure}
        
        \begin{subfigure}{0.49\textwidth}
            \scalebox{0.42}{

\begin{tikzpicture}[/tikz/background rectangle/.style={fill={rgb,1:red,1.0;green,1.0;blue,1.0}, fill opacity={1.0}, draw opacity={1.0}}, show background rectangle]
\begin{axis}[point meta max={nan}, point meta min={nan}, legend cell align={left}, legend columns={1}, title={}, title style={at={{(0.5,1)}}, anchor={south}, font={{\fontsize{14 pt}{18.2 pt}\selectfont}}, color={rgb,1:red,0.0;green,0.0;blue,0.0}, draw opacity={1.0}, rotate={0.0}, align={center}}, legend style={color={rgb,1:red,0.0;green,0.0;blue,0.0}, draw opacity={1.0}, line width={1}, solid, fill={rgb,1:red,1.0;green,1.0;blue,1.0}, fill opacity={1.0}, text opacity={1.0}, font={{\fontsize{18 pt}{10.4 pt}\selectfont}}, text={rgb,1:red,0.0;green,0.0;blue,0.0}, cells={anchor={center}}, at={(0.98, 0.98)}, anchor={north east}}, axis background/.style={fill={rgb,1:red,1.0;green,1.0;blue,1.0}, opacity={1.0}}, anchor={north west}, xshift={1.0mm}, yshift={-1.0mm}, width={150.4mm}, height={99.6mm}, scaled x ticks={false}, xlabel={Time in s}, x tick style={color={rgb,1:red,0.0;green,0.0;blue,0.0}, opacity={1.0}}, x tick label style={color={rgb,1:red,0.0;green,0.0;blue,0.0}, opacity={1.0}, rotate={0}}, xlabel style={at={(ticklabel cs:0.5)}, anchor=near ticklabel, at={{(ticklabel cs:0.5)}}, anchor={near ticklabel}, font={{\fontsize{26 pt}{14.3 pt}\selectfont}}, color={rgb,1:red,0.0;green,0.0;blue,0.0}, draw opacity={1.0}, rotate={0.0}}, xmode={log}, log basis x={10}, xmajorgrids={true}, xmin={0.0005072062585816689}, xmax={29.96859881093156}, xticklabels={{$10^{-3}$,$10^{-2}$,$10^{-1}$,$10^{0}$,$10^{1}$}}, xtick={{0.001,0.01,0.1,1.0,10.0}}, xtick align={inside}, xticklabel style={font={{\fontsize{18 pt}{10.4 pt}\selectfont}}, color={rgb,1:red,0.0;green,0.0;blue,0.0}, draw opacity={1.0}, rotate={0.0}}, x grid style={color={rgb,1:red,0.0;green,0.0;blue,0.0}, draw opacity={0.1}, line width={0.5}, solid}, axis x line*={left}, x axis line style={color={rgb,1:red,0.0;green,0.0;blue,0.0}, draw opacity={1.0}, line width={1}, solid}, scaled y ticks={false}, ylabel={Reconstruction error}, y tick style={color={rgb,1:red,0.0;green,0.0;blue,0.0}, opacity={1.0}}, y tick label style={color={rgb,1:red,0.0;green,0.0;blue,0.0}, opacity={1.0}, rotate={0}}, ylabel style={at={(ticklabel cs:0.5)}, anchor=near ticklabel, at={{(ticklabel cs:0.5)}}, anchor={near ticklabel}, font={{\fontsize{26 pt}{14.3 pt}\selectfont}}, color={rgb,1:red,0.0;green,0.0;blue,0.0}, draw opacity={1.0}, rotate={0.0}}, ymode={log}, log basis y={10}, ymajorgrids={true}, ymin={2.2380800501890885e-5}, ymax={5.1227198587161644e-5}, yticklabels={{$10^{-4.6}$,$10^{-4.4}$}}, ytick={{2.511886431509577e-5,3.9810717055349695e-5}}, ytick align={inside}, yticklabel style={font={{\fontsize{18 pt}{10.4 pt}\selectfont}}, color={rgb,1:red,0.0;green,0.0;blue,0.0}, draw opacity={1.0}, rotate={0.0}}, y grid style={color={rgb,1:red,0.0;green,0.0;blue,0.0}, draw opacity={0.1}, line width={0.5}, solid}, axis y line*={left}, y axis line style={color={rgb,1:red,0.0;green,0.0;blue,0.0}, draw opacity={1.0}, line width={1}, solid}, colorbar={false}]
    \addplot[color={rgb,1:red,0.0;green,0.6056;blue,0.9787}, name path={ca16e326-5717-483f-a32b-01a8bfe4d56d}, draw opacity={1.0}, line width={2}, solid, mark={*}, mark size={3.0 pt}, mark repeat={1}, mark options={color={rgb,1:red,0.0;green,0.6056;blue,0.9787}, draw opacity={1.0}, fill={rgb,1:red,0.0;green,0.6056;blue,0.9787}, fill opacity={1.0}, line width={0.75}, rotate={0}, solid}]
        table[row sep={\\}]
        {
            \\
            0.0018632655327999994  4.6094464976668066e-5  \\
            0.0031035786361999994  3.703798719748985e-5  \\
            0.0035661654758000006  3.3427748851838466e-5  \\
            0.006236014886400001  3.037312198118292e-5  \\
            0.014701249190199998  2.839496305657054e-5  \\
            0.0198615820652  2.7063528520967498e-5  \\
            0.037151265315399994  2.6135139695583426e-5  \\
            0.05817498278879999  2.5632358022505907e-5  \\
            0.0902285040462  2.5302460260476415e-5  \\
            0.14888378976519998  2.5109264571537557e-5  \\
        }
        ;
    \addlegendentry {MTSF}
    \addplot[color={rgb,1:red,0.8889;green,0.4356;blue,0.2781}, name path={5db64c47-b5ba-4959-ad96-ebe5d2db0d91}, draw opacity={1.0}, line width={2}, solid, mark={*}, mark size={3.0 pt}, mark repeat={1}, mark options={color={rgb,1:red,0.8889;green,0.4356;blue,0.2781}, draw opacity={1.0}, fill={rgb,1:red,0.8889;green,0.4356;blue,0.2781}, fill opacity={1.0}, line width={0.75}, rotate={0}, solid}]
        table[row sep={\\}]
        {
            \\
            0.0035923954207999994  2.5133541744194787e-5  \\
            0.0047632793134  2.5015963683811668e-5  \\
            0.0051625295425999996  2.4898923068404737e-5  \\
            0.008456681195200003  2.4848430639316365e-5  \\
            0.0173042123938  2.4848612472597792e-5  \\
            0.019383968835199997  2.4829322898094675e-5  \\
            0.0363502558686  2.481909403400939e-5  \\
            0.05355877232620001  2.481131585552164e-5  \\
            0.09120776561720005  2.481391371753439e-5  \\
            0.15565105202340007  2.4809616634883226e-5  \\
        }
        ;
    \addlegendentry {MTSF+GS}
    \addplot[color={rgb,1:red,0.2422;green,0.6433;blue,0.3044}, name path={d49e3ee1-1ca5-482c-86f8-f7284b99881d}, draw opacity={1.0}, line width={2}, solid, mark={*}, mark size={3.0 pt}, mark repeat={1}, mark options={color={rgb,1:red,0.2422;green,0.6433;blue,0.3044}, draw opacity={1.0}, fill={rgb,1:red,0.2422;green,0.6433;blue,0.3044}, fill opacity={1.0}, line width={0.75}, rotate={0}, solid}]
        table[row sep={\\}]
        {
            \\
            0.0006921919472000001  4.6871982105904135e-5  \\
            0.0011996587718  2.4471087108297717e-5  \\
            0.0016507676502  2.2911508373741305e-5  \\
            0.0034622247486  2.4724078285459098e-5  \\
            0.0067802581064  2.480345581010339e-5  \\
            0.007848105523799997  2.4803558772118917e-5  \\
            0.009274831255199996  2.480355877517328e-5  \\
            0.009763153799600001  2.480355877517328e-5  \\
            0.0100140809376  2.480355877517328e-5  \\
            0.009952526642999998  2.480355877517328e-5  \\
        }
        ;
    \addlegendentry {CG}
    \addplot[color={rgb,1:red,0.7644;green,0.4441;blue,0.8243}, name path={f62245d5-5197-4814-9be4-ad5fd3654ca8}, draw opacity={1.0}, line width={2}, solid, mark={*}, mark size={3.0 pt}, mark repeat={1}, mark options={color={rgb,1:red,0.7644;green,0.4441;blue,0.8243}, draw opacity={1.0}, fill={rgb,1:red,0.7644;green,0.4441;blue,0.8243}, fill opacity={1.0}, line width={0.75}, rotate={0}, solid}]
        table[row sep={\\}]
        {
            \\
            0.0012555345316000003  4.781039436785612e-5  \\
            0.0016766699342000005  2.3393586122956167e-5  \\
            0.0023740661493999995  2.315072797187092e-5  \\
            0.004954432405399997  2.4796055914647865e-5  \\
            0.007583537448399998  2.480355868447636e-5  \\
            0.0084280318778  2.4803558775212996e-5  \\
            0.008709491735399998  2.4803558775212996e-5  \\
            0.009477852047200001  2.4803558775212996e-5  \\
            0.009945041298599999  2.4803558775212996e-5  \\
            0.008809355617  2.4803558775212996e-5  \\
        }
        ;
    \addlegendentry {$\mathsf{D}^{-1}$+CG}
    \addplot[color={rgb,1:red,0.0;green,0.6658;blue,0.681}, name path={b6ae2df4-ca8d-4682-b28f-2e37c59d0dd5}, draw opacity={1.0}, line width={2}, dotted, forget plot]
        table[row sep={\\}]
        {
            \\
            21.959603747652  9.778013339011141e-6  \\
            21.959603747652  0.00011725344117458102  \\
        }
        ;
    \addplot[color={rgb,1:red,0.9308;green,0.3675;blue,0.5758}, name path={b883f92b-179c-4747-8d32-6f84d468a8af}, draw opacity={1.0}, line width={2}, dotted]
        table[row sep={\\}]
        {
            \\
            8.584258155258685e-9  2.4803558775335395e-5  \\
            1.7707133922243495e6  2.4803558775335395e-5  \\
        }
        ;
    \addlegendentry {$\Vert f_* - f_\top \Vert_2 / n$}
    \addplot[color={rgb,1:red,0.777;green,0.5097;blue,0.1464}, name path={2c887e9e-4ca6-4c4d-8976-2e7511fc8000}, draw opacity={1.0}, line width={2}, dotted]
        table[row sep={\\}]
        {
            \\
            8.584258155258685e-9  5.004060375021011e-5  \\
            1.7707133922243495e6  5.004060375021011e-5  \\
        }
        ;
    \addlegendentry {$\Vert g - f_\top \Vert_2 / n$}
\end{axis}
\end{tikzpicture}}
            \caption{SBM Reconstruction error.}
        \end{subfigure}
        \hfill
        \begin{subfigure}{0.49\textwidth}
            \scalebox{0.42}{

\begin{tikzpicture}[/tikz/background rectangle/.style={fill={rgb,1:red,1.0;green,1.0;blue,1.0}, fill opacity={1.0}, draw opacity={1.0}}, show background rectangle]
\begin{axis}[point meta max={nan}, point meta min={nan}, legend cell align={left}, legend columns={1}, title={}, title style={at={{(0.5,1)}}, anchor={south}, font={{\fontsize{14 pt}{18.2 pt}\selectfont}}, color={rgb,1:red,0.0;green,0.0;blue,0.0}, draw opacity={1.0}, rotate={0.0}, align={center}}, legend style={color={rgb,1:red,0.0;green,0.0;blue,0.0}, draw opacity={1.0}, line width={1}, solid, fill={rgb,1:red,1.0;green,1.0;blue,1.0}, fill opacity={1.0}, text opacity={1.0}, font={{\fontsize{18 pt}{10.4 pt}\selectfont}}, text={rgb,1:red,0.0;green,0.0;blue,0.0}, cells={anchor={center}}, at={(0.98, 0.98)}, anchor={north east}}, axis background/.style={fill={rgb,1:red,1.0;green,1.0;blue,1.0}, opacity={1.0}}, anchor={north west}, xshift={1.0mm}, yshift={-1.0mm}, width={150.4mm}, height={99.6mm}, scaled x ticks={false}, xlabel={Time in s}, x tick style={color={rgb,1:red,0.0;green,0.0;blue,0.0}, opacity={1.0}}, x tick label style={color={rgb,1:red,0.0;green,0.0;blue,0.0}, opacity={1.0}, rotate={0}}, xlabel style={at={(ticklabel cs:0.5)}, anchor=near ticklabel, at={{(ticklabel cs:0.5)}}, anchor={near ticklabel}, font={{\fontsize{26 pt}{14.3 pt}\selectfont}}, color={rgb,1:red,0.0;green,0.0;blue,0.0}, draw opacity={1.0}, rotate={0.0}}, xmode={log}, log basis x={10}, xmajorgrids={true}, xmin={0.00045315390508036106}, xmax={1.5327119748742344}, xticklabels={{$10^{-3}$,$10^{-2}$,$10^{-1}$,$10^{0}$}}, xtick={{0.001,0.01,0.1,1.0}}, xtick align={inside}, xticklabel style={font={{\fontsize{18 pt}{10.4 pt}\selectfont}}, color={rgb,1:red,0.0;green,0.0;blue,0.0}, draw opacity={1.0}, rotate={0.0}}, x grid style={color={rgb,1:red,0.0;green,0.0;blue,0.0}, draw opacity={0.1}, line width={0.5}, solid}, axis x line*={left}, x axis line style={color={rgb,1:red,0.0;green,0.0;blue,0.0}, draw opacity={1.0}, line width={1}, solid}, scaled y ticks={false}, ylabel={Reconstruction error}, y tick style={color={rgb,1:red,0.0;green,0.0;blue,0.0}, opacity={1.0}}, y tick label style={color={rgb,1:red,0.0;green,0.0;blue,0.0}, opacity={1.0}, rotate={0}}, ylabel style={at={(ticklabel cs:0.5)}, anchor=near ticklabel, at={{(ticklabel cs:0.5)}}, anchor={near ticklabel}, font={{\fontsize{26 pt}{14.3 pt}\selectfont}}, color={rgb,1:red,0.0;green,0.0;blue,0.0}, draw opacity={1.0}, rotate={0.0}}, ymode={log}, log basis y={10}, ymajorgrids={true}, ymin={1.747312360427883e-5}, ymax={7.968969908840074e-5}, yticklabels={{$10^{-4.5}$,$10^{-4.2}$}}, ytick={{3.162277660168373e-5,6.309573444801917e-5}}, ytick align={inside}, yticklabel style={font={{\fontsize{18 pt}{10.4 pt}\selectfont}}, color={rgb,1:red,0.0;green,0.0;blue,0.0}, draw opacity={1.0}, rotate={0.0}}, y grid style={color={rgb,1:red,0.0;green,0.0;blue,0.0}, draw opacity={0.1}, line width={0.5}, solid}, axis y line*={left}, y axis line style={color={rgb,1:red,0.0;green,0.0;blue,0.0}, draw opacity={1.0}, line width={1}, solid}, colorbar={false}]
    \addplot[color={rgb,1:red,0.0;green,0.6056;blue,0.9787}, name path={b36fe5a7-f6b9-427f-91a1-442065d72a89}, draw opacity={1.0}, line width={2}, solid, mark={*}, mark size={3.0 pt}, mark repeat={1}, mark options={color={rgb,1:red,0.0;green,0.6056;blue,0.9787}, draw opacity={1.0}, fill={rgb,1:red,0.0;green,0.6056;blue,0.9787}, fill opacity={1.0}, line width={0.75}, rotate={0}, solid}]
        table[row sep={\\}]
        {
            \\
            0.0015308698232000002  3.58683847044972e-5  \\
            0.0024008618785999994  3.064220908553685e-5  \\
            0.0033631265204  2.416880809224033e-5  \\
            0.005431243193599998  2.3228929799882224e-5  \\
            0.012462797532799998  2.1639479813323863e-5  \\
            0.0168686023766  2.1037033945762114e-5  \\
            0.0319691218302  1.994763611480616e-5  \\
            0.04808730418979999  1.99348635284733e-5  \\
            0.07780924720319998  1.98373931909197e-5  \\
            0.1353063975734  1.9229280453678935e-5  \\
        }
        ;
    \addlegendentry {MTSF}
    \addplot[color={rgb,1:red,0.8889;green,0.4356;blue,0.2781}, name path={b7aa09ba-b325-4cb1-9626-3426b18a85d4}, draw opacity={1.0}, line width={2}, solid, mark={*}, mark size={3.0 pt}, mark repeat={1}, mark options={color={rgb,1:red,0.8889;green,0.4356;blue,0.2781}, draw opacity={1.0}, fill={rgb,1:red,0.8889;green,0.4356;blue,0.2781}, fill opacity={1.0}, line width={0.75}, rotate={0}, solid}]
        table[row sep={\\}]
        {
            \\
            0.0032501763943999997  1.9754309854612547e-5  \\
            0.003622600276000001  1.9408258272925898e-5  \\
            0.004974376870400002  1.91048482906192e-5  \\
            0.0071287508732000016  1.950773105041196e-5  \\
            0.015786409657799996  1.895868487034726e-5  \\
            0.016996112962399997  1.888388827908881e-5  \\
            0.034595695165400016  1.9238859772626545e-5  \\
            0.052620975981799996  1.902530208278843e-5  \\
            0.085794776625  1.903927603236167e-5  \\
            0.1391596934482  1.9145526042606387e-5  \\
        }
        ;
    \addlegendentry {MTSF+GS}
    \addplot[color={rgb,1:red,0.2422;green,0.6433;blue,0.3044}, name path={38009aa5-3cc8-4e18-88e2-0a24e13170fa}, draw opacity={1.0}, line width={2}, solid, mark={*}, mark size={3.0 pt}, mark repeat={1}, mark options={color={rgb,1:red,0.2422;green,0.6433;blue,0.3044}, draw opacity={1.0}, fill={rgb,1:red,0.2422;green,0.6433;blue,0.3044}, fill opacity={1.0}, line width={0.75}, rotate={0}, solid}]
        table[row sep={\\}]
        {
            \\
            0.0005703338724000001  7.633968227869563e-5  \\
            0.0009143089180000002  6.740130060592489e-5  \\
            0.0014550329074  4.345869692458088e-5  \\
            0.0026778562118000013  3.6022966656652674e-5  \\
            0.005474371801999999  2.6020978396068473e-5  \\
            0.006442219036199997  2.005992896748265e-5  \\
            0.0125097720214  1.8239897267007573e-5  \\
            0.021182594912399995  1.8957462656396435e-5  \\
            0.035969484678000006  1.9030753219436857e-5  \\
            0.056163568534  1.9031260384333593e-5  \\
        }
        ;
    \addlegendentry {CG}
    \addplot[color={rgb,1:red,0.7644;green,0.4441;blue,0.8243}, name path={08bc8543-cb69-4b84-8922-63f17df3080b}, draw opacity={1.0}, line width={2}, solid, mark={*}, mark size={3.0 pt}, mark repeat={1}, mark options={color={rgb,1:red,0.7644;green,0.4441;blue,0.8243}, draw opacity={1.0}, fill={rgb,1:red,0.7644;green,0.4441;blue,0.8243}, fill opacity={1.0}, line width={0.75}, rotate={0}, solid}]
        table[row sep={\\}]
        {
            \\
            0.0010923790983999998  6.356999907825005e-5  \\
            0.0014361057674000001  5.582493698922853e-5  \\
            0.0020768403252  2.213021217430731e-5  \\
            0.0038239560488000004  1.9067786945840237e-5  \\
            0.0064090582356  1.90311705710683e-5  \\
            0.008190323317200002  1.9031260483443813e-5  \\
            0.010073858436199996  1.9031260423904773e-5  \\
            0.010642968164399998  1.9031260423904773e-5  \\
            0.0112713690074  1.9031260423904773e-5  \\
            0.011425916066200003  1.9031260423904773e-5  \\
        }
        ;
    \addlegendentry {$\mathsf{D}^{-1}$+CG}
    \addplot[color={rgb,1:red,0.0;green,0.6658;blue,0.681}, name path={8076148f-ce5f-4b2a-984b-f3503280c2e6}, draw opacity={1.0}, line width={2}, dotted, forget plot]
        table[row sep={\\}]
        {
            \\
            1.2178032033323996  3.831236056641673e-6  \\
            1.2178032033323996  0.0003634409213041196  \\
        }
        ;
    \addplot[color={rgb,1:red,0.9308;green,0.3675;blue,0.5758}, name path={7defb784-cf8a-4dc8-a5d8-4d2ffea743d0}, draw opacity={1.0}, line width={2}, dotted]
        table[row sep={\\}]
        {
            \\
            1.3397720188519485e-7  1.903126042498225e-5  \\
            5184.123918133888  1.903126042498225e-5  \\
        }
        ;
    \addlegendentry {$\Vert f_* - f_\top \Vert_2 / n$}
\end{axis}
\end{tikzpicture}}
            \caption{DC-SBM $1$ Reconstruction error.}
        \end{subfigure}
        
        \begin{subfigure}{0.49\textwidth}
            \scalebox{0.42}{

\begin{tikzpicture}[/tikz/background rectangle/.style={fill={rgb,1:red,1.0;green,1.0;blue,1.0}, fill opacity={1.0}, draw opacity={1.0}}, show background rectangle]
\begin{axis}[point meta max={nan}, point meta min={nan}, legend cell align={left}, legend columns={1}, title={}, title style={at={{(0.5,1)}}, anchor={south}, font={{\fontsize{14 pt}{18.2 pt}\selectfont}}, color={rgb,1:red,0.0;green,0.0;blue,0.0}, draw opacity={1.0}, rotate={0.0}, align={center}}, legend style={color={rgb,1:red,0.0;green,0.0;blue,0.0}, draw opacity={1.0}, line width={1}, solid, fill={rgb,1:red,1.0;green,1.0;blue,1.0}, fill opacity={1.0}, text opacity={1.0}, font={{\fontsize{18 pt}{10.4 pt}\selectfont}}, text={rgb,1:red,0.0;green,0.0;blue,0.0}, cells={anchor={center}}, at={(0.98, 0.98)}, anchor={north east}}, axis background/.style={fill={rgb,1:red,1.0;green,1.0;blue,1.0}, opacity={1.0}}, anchor={north west}, xshift={1.0mm}, yshift={-1.0mm}, width={150.4mm}, height={99.6mm}, scaled x ticks={false}, xlabel={Time in s}, x tick style={color={rgb,1:red,0.0;green,0.0;blue,0.0}, opacity={1.0}}, x tick label style={color={rgb,1:red,0.0;green,0.0;blue,0.0}, opacity={1.0}, rotate={0}}, xlabel style={at={(ticklabel cs:0.5)}, anchor=near ticklabel, at={{(ticklabel cs:0.5)}}, anchor={near ticklabel}, font={{\fontsize{26 pt}{14.3 pt}\selectfont}}, color={rgb,1:red,0.0;green,0.0;blue,0.0}, draw opacity={1.0}, rotate={0.0}}, xmode={log}, log basis x={10}, xmajorgrids={true}, xmin={0.002249003586429591}, xmax={7.882259078854072}, xticklabels={{$10^{-2}$,$10^{-1}$,$10^{0}$}}, xtick={{0.01,0.1,1.0}}, xtick align={inside}, xticklabel style={font={{\fontsize{18 pt}{10.4 pt}\selectfont}}, color={rgb,1:red,0.0;green,0.0;blue,0.0}, draw opacity={1.0}, rotate={0.0}}, x grid style={color={rgb,1:red,0.0;green,0.0;blue,0.0}, draw opacity={0.1}, line width={0.5}, solid}, axis x line*={left}, x axis line style={color={rgb,1:red,0.0;green,0.0;blue,0.0}, draw opacity={1.0}, line width={1}, solid}, scaled y ticks={false}, ylabel={Reconstruction error}, y tick style={color={rgb,1:red,0.0;green,0.0;blue,0.0}, opacity={1.0}}, y tick label style={color={rgb,1:red,0.0;green,0.0;blue,0.0}, opacity={1.0}, rotate={0}}, ylabel style={at={(ticklabel cs:0.5)}, anchor=near ticklabel, at={{(ticklabel cs:0.5)}}, anchor={near ticklabel}, font={{\fontsize{26 pt}{14.3 pt}\selectfont}}, color={rgb,1:red,0.0;green,0.0;blue,0.0}, draw opacity={1.0}, rotate={0.0}}, ymode={log}, log basis y={10}, ymajorgrids={true}, ymin={1.3740623880537794e-5}, ymax={9.713746792460557e-5}, yticklabels={{$10^{-4.75}$,$10^{-4.50}$,$10^{-4.25}$}}, ytick={{1.778279410038923e-5,3.1622776601683795e-5,5.623413251903491e-5}}, ytick align={inside}, yticklabel style={font={{\fontsize{18 pt}{10.4 pt}\selectfont}}, color={rgb,1:red,0.0;green,0.0;blue,0.0}, draw opacity={1.0}, rotate={0.0}}, y grid style={color={rgb,1:red,0.0;green,0.0;blue,0.0}, draw opacity={0.1}, line width={0.5}, solid}, axis y line*={left}, y axis line style={color={rgb,1:red,0.0;green,0.0;blue,0.0}, draw opacity={1.0}, line width={1}, solid}, colorbar={false}]
    \addplot[color={rgb,1:red,0.0;green,0.6056;blue,0.9787}, name path={6d6f85b0-2b43-4f8b-a451-a7f5cdcd9167}, draw opacity={1.0}, line width={2}, solid, mark={*}, mark size={3.0 pt}, mark repeat={1}, mark options={color={rgb,1:red,0.0;green,0.6056;blue,0.9787}, draw opacity={1.0}, fill={rgb,1:red,0.0;green,0.6056;blue,0.9787}, fill opacity={1.0}, line width={0.75}, rotate={0}, solid}]
        table[row sep={\\}]
        {
            \\
            0.0028334185301999997  2.1040209817350433e-5  \\
            0.0037379942884000006  1.955374993018194e-5  \\
            0.0084193877026  1.9095491235337058e-5  \\
            0.010181012451600002  1.836096878096219e-5  \\
            0.0189446295446  1.8373916093626262e-5  \\
            0.031308930622200006  1.616412356286385e-5  \\
            0.038426807412599996  1.489290726566159e-5  \\
            0.08536121517520003  1.5134227100681767e-5  \\
            0.1287781440518  1.4652880673868332e-5  \\
            0.20646168981120003  1.4529083151318067e-5  \\
        }
        ;
    \addlegendentry {MTSF}
    \addplot[color={rgb,1:red,0.8889;green,0.4356;blue,0.2781}, name path={fcbf4e37-62d6-4d19-9587-edf63c614530}, draw opacity={1.0}, line width={2}, solid, mark={*}, mark size={3.0 pt}, mark repeat={1}, mark options={color={rgb,1:red,0.8889;green,0.4356;blue,0.2781}, draw opacity={1.0}, fill={rgb,1:red,0.8889;green,0.4356;blue,0.2781}, fill opacity={1.0}, line width={0.75}, rotate={0}, solid}]
        table[row sep={\\}]
        {
            \\
            0.0132382501892  1.4533952116559921e-5  \\
            0.011192393792799998  1.4666388467189694e-5  \\
            0.018433897098799997  1.4742905088188501e-5  \\
            0.021356811400000005  1.4620339303767194e-5  \\
            0.0292819788492  1.4690461868593096e-5  \\
            0.031296631864399994  1.465452016557325e-5  \\
            0.05138992778119999  1.4651529602723103e-5  \\
            0.0704137448438  1.4617802932388759e-5  \\
            0.142327679929  1.4655915604597402e-5  \\
            0.22315145818339996  1.4644127842001571e-5  \\
        }
        ;
    \addlegendentry {MTSF+GS}
    \addplot[color={rgb,1:red,0.2422;green,0.6433;blue,0.3044}, name path={98c03bab-1b02-435f-ab2d-ea643b22a8e6}, draw opacity={1.0}, line width={2}, solid, mark={*}, mark size={3.0 pt}, mark repeat={1}, mark options={color={rgb,1:red,0.2422;green,0.6433;blue,0.3044}, draw opacity={1.0}, fill={rgb,1:red,0.2422;green,0.6433;blue,0.3044}, fill opacity={1.0}, line width={0.75}, rotate={0}, solid}]
        table[row sep={\\}]
        {
            \\
            0.007605919363799998  9.190681429718075e-5  \\
            0.013770885196999998  6.218023387907517e-5  \\
            0.023442879375000003  5.3821598225429655e-5  \\
            0.041699596503399994  4.006410566202006e-5  \\
            0.0542357055444  3.0306346950209035e-5  \\
            0.09653530486560001  1.8764434953247497e-5  \\
            0.19053918895259994  1.4703615681614933e-5  \\
            0.29627084183919994  1.452263819246265e-5  \\
            0.4793237683386001  1.4615146871092484e-5  \\
            0.7939107247726  1.4633457729091718e-5  \\
        }
        ;
    \addlegendentry {CG}
    \addplot[color={rgb,1:red,0.7644;green,0.4441;blue,0.8243}, name path={63696de0-6393-4303-b8a1-a24de0d6839f}, draw opacity={1.0}, line width={2}, solid, mark={*}, mark size={3.0 pt}, mark repeat={1}, mark options={color={rgb,1:red,0.7644;green,0.4441;blue,0.8243}, draw opacity={1.0}, fill={rgb,1:red,0.7644;green,0.4441;blue,0.8243}, fill opacity={1.0}, line width={0.75}, rotate={0}, solid}]
        table[row sep={\\}]
        {
            \\
            0.008095425743399998  6.404726428895852e-5  \\
            0.017744916993200004  6.186110463768229e-5  \\
            0.022478603143000004  3.2360352773168454e-5  \\
            0.0398031081262  1.4671243133032395e-5  \\
            0.07672043277299999  1.4635048312711097e-5  \\
            0.10974329127839998  1.463356340219918e-5  \\
            0.09112879555300002  1.4633563382536939e-5  \\
            0.1027907322758  1.4633563382536939e-5  \\
            0.10555339700040003  1.4633563382536939e-5  \\
            0.09911791991000003  1.4633563382536939e-5  \\
        }
        ;
    \addlegendentry {$\mathsf{D}^{-1}$+CG}
    \addplot[color={rgb,1:red,0.0;green,0.6658;blue,0.681}, name path={18d02316-d1b8-4cb0-80e0-7b695f05f35b}, draw opacity={1.0}, line width={2}, dotted, forget plot]
        table[row sep={\\}]
        {
            \\
            6.256480907625998  1.9436860838595116e-6  \\
            6.256480907625998  0.0006867000914105854  \\
        }
        ;
    \addplot[color={rgb,1:red,0.9308;green,0.3675;blue,0.5758}, name path={37f53a3c-d1f9-4c8f-bf98-e5a704e50af4}, draw opacity={1.0}, line width={2}, dotted]
        table[row sep={\\}]
        {
            \\
            6.416963818586512e-7  1.4633563382513295e-5  \\
            27625.570968880504  1.4633563382513295e-5  \\
        }
        ;
    \addlegendentry {$\Vert f_* - f_\top \Vert_2 / n$}
    \addplot[color={rgb,1:red,0.777;green,0.5097;blue,0.1464}, name path={f8b440d8-12d7-4d92-828e-b7d2383081ad}, draw opacity={1.0}, line width={2}, dotted]
        table[row sep={\\}]
        {
            \\
            6.416963818586512e-7  5.0379026128082445e-5  \\
            27625.570968880504  5.0379026128082445e-5  \\
        }
        ;
    \addlegendentry {$\Vert g - f_\top \Vert_2 / n$}
\end{axis}
\end{tikzpicture}}
            \caption{DC-SBM $2$ Reconstruction error.}
        \end{subfigure}
        \hfill
        \begin{subfigure}{0.49\textwidth}
            \scalebox{0.42}{

\begin{tikzpicture}[/tikz/background rectangle/.style={fill={rgb,1:red,1.0;green,1.0;blue,1.0}, fill opacity={1.0}, draw opacity={1.0}}, show background rectangle]
\begin{axis}[point meta max={nan}, point meta min={nan}, legend cell align={left}, legend columns={1}, title={}, title style={at={{(0.5,1)}}, anchor={south}, font={{\fontsize{14 pt}{18.2 pt}\selectfont}}, color={rgb,1:red,0.0;green,0.0;blue,0.0}, draw opacity={1.0}, rotate={0.0}, align={center}}, legend style={color={rgb,1:red,0.0;green,0.0;blue,0.0}, draw opacity={1.0}, line width={1}, solid, fill={rgb,1:red,1.0;green,1.0;blue,1.0}, fill opacity={1.0}, text opacity={1.0}, font={{\fontsize{18 pt}{10.4 pt}\selectfont}}, text={rgb,1:red,0.0;green,0.0;blue,0.0}, cells={anchor={center}}, at={(0.98, 0.98)}, anchor={north east}}, axis background/.style={fill={rgb,1:red,1.0;green,1.0;blue,1.0}, opacity={1.0}}, anchor={north west}, xshift={1.0mm}, yshift={-1.0mm}, width={150.4mm}, height={99.6mm}, scaled x ticks={false}, xlabel={Time in s}, x tick style={color={rgb,1:red,0.0;green,0.0;blue,0.0}, opacity={1.0}}, x tick label style={color={rgb,1:red,0.0;green,0.0;blue,0.0}, opacity={1.0}, rotate={0}}, xlabel style={at={(ticklabel cs:0.5)}, anchor=near ticklabel, at={{(ticklabel cs:0.5)}}, anchor={near ticklabel}, font={{\fontsize{26 pt}{14.3 pt}\selectfont}}, color={rgb,1:red,0.0;green,0.0;blue,0.0}, draw opacity={1.0}, rotate={0.0}}, xmode={log}, log basis x={10}, xmajorgrids={true}, xmin={0.0006436494332993045}, xmax={0.4529198568985317}, xticklabels={{$10^{-3}$,$10^{-2}$,$10^{-1}$}}, xtick={{0.001,0.01,0.1}}, xtick align={inside}, xticklabel style={font={{\fontsize{18 pt}{10.4 pt}\selectfont}}, color={rgb,1:red,0.0;green,0.0;blue,0.0}, draw opacity={1.0}, rotate={0.0}}, x grid style={color={rgb,1:red,0.0;green,0.0;blue,0.0}, draw opacity={0.1}, line width={0.5}, solid}, axis x line*={left}, x axis line style={color={rgb,1:red,0.0;green,0.0;blue,0.0}, draw opacity={1.0}, line width={1}, solid}, scaled y ticks={false}, ylabel={Reconstruction error}, y tick style={color={rgb,1:red,0.0;green,0.0;blue,0.0}, opacity={1.0}}, y tick label style={color={rgb,1:red,0.0;green,0.0;blue,0.0}, opacity={1.0}, rotate={0}}, ylabel style={at={(ticklabel cs:0.5)}, anchor=near ticklabel, at={{(ticklabel cs:0.5)}}, anchor={near ticklabel}, font={{\fontsize{26 pt}{14.3 pt}\selectfont}}, color={rgb,1:red,0.0;green,0.0;blue,0.0}, draw opacity={1.0}, rotate={0.0}}, ymode={log}, log basis y={10}, ymajorgrids={true}, ymin={5.991924962871175e-6}, ymax={3.805410797196153e-5}, yticklabels={{$10^{-5.00}$,$10^{-4.75}$,$10^{-4.50}$}}, ytick={{1.0e-5,1.778279410038923e-5,3.1622776601683795e-5}}, ytick align={inside}, yticklabel style={font={{\fontsize{18 pt}{10.4 pt}\selectfont}}, color={rgb,1:red,0.0;green,0.0;blue,0.0}, draw opacity={1.0}, rotate={0.0}}, y grid style={color={rgb,1:red,0.0;green,0.0;blue,0.0}, draw opacity={0.1}, line width={0.5}, solid}, axis y line*={left}, y axis line style={color={rgb,1:red,0.0;green,0.0;blue,0.0}, draw opacity={1.0}, line width={1}, solid}, colorbar={false}]
    \addplot[color={rgb,1:red,0.0;green,0.6056;blue,0.9787}, name path={6022068b-0f91-4422-b8fc-5669a9df3ab6}, draw opacity={1.0}, line width={2}, solid, mark={*}, mark size={3.0 pt}, mark repeat={1}, mark options={color={rgb,1:red,0.0;green,0.6056;blue,0.9787}, draw opacity={1.0}, fill={rgb,1:red,0.0;green,0.6056;blue,0.9787}, fill opacity={1.0}, line width={0.75}, rotate={0}, solid}]
        table[row sep={\\}]
        {
            \\
            0.004203589204399999  1.5954757916209495e-5  \\
            0.006616798989599999  1.226037902448055e-5  \\
            0.0093964448744  1.0836233810845338e-5  \\
            0.018933994913199997  9.53563851493414e-6  \\
            0.031089204912399995  8.439801113431489e-6  \\
            0.0531018398106  7.6255700998452935e-6  \\
            0.07883723285800001  7.21422473626598e-6  \\
            0.13423450764400002  6.961961374835717e-6  \\
            0.2294380340380001  6.7690067164657675e-6  \\
            0.3722465075194  6.618327530894972e-6  \\
        }
        ;
    \addlegendentry {MTSF}
    \addplot[color={rgb,1:red,0.8889;green,0.4356;blue,0.2781}, name path={019bf06c-312d-4a2f-a1b6-a649db9e3fdd}, draw opacity={1.0}, line width={2}, solid, mark={*}, mark size={3.0 pt}, mark repeat={1}, mark options={color={rgb,1:red,0.8889;green,0.4356;blue,0.2781}, draw opacity={1.0}, fill={rgb,1:red,0.8889;green,0.4356;blue,0.2781}, fill opacity={1.0}, line width={0.75}, rotate={0}, solid}]
        table[row sep={\\}]
        {
            \\
            0.007717091006599999  1.0724680011677435e-5  \\
            0.008408732203600002  8.998537671719745e-6  \\
            0.016178346834599996  8.119443056139827e-6  \\
            0.0216878288002  7.675397534168341e-6  \\
            0.029592573302199995  7.176255348355279e-6  \\
            0.0493041575274  7.022050258774342e-6  \\
            0.08782824869660001  6.820015292852992e-6  \\
            0.13732254946519998  6.594338255972018e-6  \\
            0.22319866054179996  6.569154995222209e-6  \\
            0.37621435466  6.520628235046238e-6  \\
        }
        ;
    \addlegendentry {MTSF+GS}
    \addplot[color={rgb,1:red,0.2422;green,0.6433;blue,0.3044}, name path={ecf271d8-823f-477e-a0b8-1b738fec5d74}, draw opacity={1.0}, line width={2}, solid, mark={*}, mark size={3.0 pt}, mark repeat={1}, mark options={color={rgb,1:red,0.2422;green,0.6433;blue,0.3044}, draw opacity={1.0}, fill={rgb,1:red,0.2422;green,0.6433;blue,0.3044}, fill opacity={1.0}, line width={0.75}, rotate={0}, solid}]
        table[row sep={\\}]
        {
            \\
            0.0007748816747999998  2.7505812480431975e-5  \\
            0.0008901926748000001  2.574916138824968e-5  \\
            0.0020024384916000008  2.4695727830668818e-5  \\
            0.0026855569585999997  2.320938505421189e-5  \\
            0.0037169597222  2.114258707266317e-5  \\
            0.0061637063457999994  1.791174473555285e-5  \\
            0.0097823743888  1.156642297664856e-5  \\
            0.017846310528400005  7.788650576162001e-6  \\
            0.0324261167112  6.38725930962035e-6  \\
            0.05344942935540001  6.41541807537013e-6  \\
        }
        ;
    \addlegendentry {CG}
    \addplot[color={rgb,1:red,0.7644;green,0.4441;blue,0.8243}, name path={e4da7af6-f540-4efa-9846-ac1c45b129b1}, draw opacity={1.0}, line width={2}, solid, mark={*}, mark size={3.0 pt}, mark repeat={1}, mark options={color={rgb,1:red,0.7644;green,0.4441;blue,0.8243}, draw opacity={1.0}, fill={rgb,1:red,0.7644;green,0.4441;blue,0.8243}, fill opacity={1.0}, line width={0.75}, rotate={0}, solid}]
        table[row sep={\\}]
        {
            \\
            0.0019047500677999997  3.611434942079985e-5  \\
            0.0019027722445999998  2.0414343340576078e-5  \\
            0.0037701460837999998  1.2079388637078168e-5  \\
            0.0044477963648  7.408432617729378e-6  \\
            0.0073962482116  6.313760683881132e-6  \\
            0.011443700883799996  6.454324410388445e-6  \\
            0.017690212085199997  6.457419664766225e-6  \\
            0.0234057576792  6.457396699801601e-6  \\
            0.027562111721599992  6.457396701676616e-6  \\
            0.027844016269199995  6.457396701676616e-6  \\
        }
        ;
    \addlegendentry {$\mathsf{D}^{-1}$+CG}
    \addplot[color={rgb,1:red,0.0;green,0.6658;blue,0.681}, name path={6b1de642-e29e-4aee-ab4d-ee8280aa392c}, draw opacity={1.0}, line width={2}, dotted, forget plot]
        table[row sep={\\}]
        {
            \\
            0.05405087259039999  9.434767144491305e-7  \\
            0.05405087259039999  0.00024167778176711457  \\
        }
        ;
    \addplot[color={rgb,1:red,0.9308;green,0.3675;blue,0.5758}, name path={d10c644f-ce3b-4bdf-8ddb-65118641d385}, draw opacity={1.0}, line width={2}, dotted]
        table[row sep={\\}]
        {
            \\
            9.146973502628493e-7  6.457396700199147e-6  \\
            318.70826906732566  6.457396700199147e-6  \\
        }
        ;
    \addlegendentry {$\Vert f_* - f_\top \Vert_2 / n$}
    \addplot[color={rgb,1:red,0.777;green,0.5097;blue,0.1464}, name path={ac579dfb-dc3e-4008-9590-34d446a0cb99}, draw opacity={1.0}, line width={2}, dotted]
        table[row sep={\\}]
        {
            \\
            9.146973502628493e-7  1.888403144070378e-5  \\
            318.70826906732566  1.888403144070378e-5  \\
        }
        ;
    \addlegendentry {$\Vert g - f_\top \Vert_2 / n$}
\end{axis}
\end{tikzpicture}}
            \caption{AS CAIDA Reconstruction error.}
        \end{subfigure}
        \caption{Runtime-precision trade-offs for Tikhonov smoothing with optimal $q$. Each datapoint corresponds to a value of $m$. The curve labelled MTSF corresponds to the estimator $\overline{f}$, and the curve labelled MTSF+GS to the estimator $\hat{f}$. The vertical blue line records the average runtime of the Cholesky-based solver, the horizontal red line the average reconstruction error of the exact solution to the Tikhonov problem, and the horizontal yellow line the average reconstruction error of the noisy signal. Vertical bars represent standard deviations.}
        \label{fig:smoothing}
    \end{figure}

    \begin{paragraph}{Setup}
        For each graph, we generate a $B$-bandlimited signal $f_\top$, each $u_i$ being weighted by a random complex Gaussian value $a_i \sim \mathcal{N}_\C(0,1)$. We then degrade $f_\top$ with some additive Gaussian noise $\eps \sim \mathcal{N}_\C(0,\sigma^2)$ (independently on each entry), with $\sigma^2$ such that the SNR is equal to $2$, and determine the optimal parameter $q_*$ for which $e_r(f_*) = e_r( q_* (\mathsf{L}_\theta + q_* \mathsf{I})^{-1} g)$ is minimized\footnote{We perform our search in $(0,30)$.}, before measuring the errors and running time associated to different values of $m$ (taken amongst 10 logarithmically-spaced values from $1$ to $100$). The iterative algorithms are initialized at $g$.
        \\
        We compare the following estimation strategies: the estimators $\overline{f}$ and $\hat{f}$, the conjugate-gradient descent with no preconditioner, with a simple diagonal preconditioner $\mathsf{P} = ( q^{-1} \mathsf{D} + \mathsf{I})^{-1}$, and with a CROUT ILU preconditioner\footnote{We use the implementations from the libraries: \url{https://github.com/JuliaLinearAlgebra/IterativeSolvers.jl}, \url{https://github.com/JuliaLinearAlgebra/Preconditioners.jl} and \url{https://github.com/haampie/IncompleteLU.jl}. \\ For the ILU preconditioner, we fix the drop threshold to $0.1$.}. Computing this high-quality ILU preconditioner is too expensive to be competitive with the other methods, and we only include these results on one type of graph, for illustration purposes.
        \\
        Each runtime measurement is averaged over $100$ runs. The results are averaged over $5$ realizations of the noise $\eps$ and $10$ samples for each random graph model. We plot in Figure~\ref{fig:smoothing} the mean results.
    \end{paragraph}

    \begin{paragraph}{Graphs used}
        We use the following graphs:
        \begin{itemize}
            \item An $\eps$-graph obtained by sampling $\vert \mathcal{V} \vert = 10^4$ i.i.d. points $x_i$'s in $[0,1]^3$, with an edge between $x_i$ and $x_j$ whenever $\Vert x_i - x_j \Vert_2 < 0.1$.
            \item A graph generated from a \emph{Stochastic Block Model} (SBM), with $\vert \mathcal{V} \vert = 10^4$ nodes, each belonging to one of two communities $C_1$ and $C_2$ of size $5000$ each. In this model, an edge is drawn randomly between two nodes $i,j$ with probability\footnote{It is usually assumed that $n > c_{k,l}$, so that we always have $\frac{c_{k,l}}{n} \in [0,1]$. Most actual implementations use the probabilities $\min\left(\frac{c_{k,l}}{n},1\right)$ in case this assumption is not satisfied.} $\frac{c_{k,l}}{n}$, where $k,l$ denote respectively the community label of $i$ and $j$. Here, we take $c_{1,1} = c_{2,2} = 36$ and $c_{1,2} = c_{2,1} = 4$, resulting in an average degree $\overline{d} = 40$.
            \item A graph generated from a related \emph{Degree-Corrected Stochastic Block Model} (DC-SBM $1$), with two communities of size $5000$ and an edge between nodes $i \in C_k$ and $j \in C_l$ present with probability $p_i p_j \frac{c_{k,l}}{n}$. $p_i$ is a randomly sampled positive real value, representing the intrinsic \emph{connectivity} of node $i$, with $\mathbf{E}(p_i) = 1$ and finite second moment~\cite{karrer2011stochastic}. Here, we take $c_{1,1} = c_{2,2} = 36$, $c_{1,2} = 4$, and the $p_i$'s are drawn from a normalized mixture of Gaussian distributions\footnote{Specifically, the connectivity parameters are independently sampled from a mixture of $\mathcal{N}(50,20)$, $\mathcal{N}(500,100)$ and $\mathcal{N}(10000,100)$, with weights of $0.59$, $0.4$ and $0.01$ respectively, and then normalized.}, resulting in a graph with mean degree $\overline{d}$ close to  $40$. The objective of adding this model is to illustrate, at constant density, how the degree distribution affects the results.
            \item Another DC-SBM-graph (DC-SBM $2$), with higher density (typically with an average degree more than $10$ times that of the previous DC-SBM model). We take two communities of size $5000$, $c_{1,1} = c_{2,2} = 480$, $c_{1,2} = 20$, and the $p_i$'s are drawn from another mixture of Gaussian distributions\footnote{Here from the mixture of $\mathcal{N}(50,20)$, $\mathcal{N}(1000,50)$, $\mathcal{N}(5000,100)$ and $\mathcal{N}(10000,100)$, with weights of $0.45$, $0.1$, $0.44$ and $0.01$.}. The objective of this second DC-SBM model is to illustrate how density affects the results.
            \item We also illustrate the results on a real-world graph: a relationship graph for internet Autonomous Systems (AS), recorded by the Center for Applied Internet Data Analysis (CAIDA) and provided in the SNAP datasets~\cite{snapnets}\footnote{We use the library available at \url{https://github.com/JuliaGraphs/SNAPDatasets.jl}.}.
        \end{itemize} 
        In the event of a randomly generated graph containing isolated nodes, we remove them so that the graph is connected (to make the interpretations simpler). We summarize in Table~\ref{table:graphs} the information concerning the graphs generated, and the associated optimal $q_*$'s (averaged over all realizations of the noise $\eps)$. All those graphs are endowed with a synthetic connection as specified in Eq.~\eqref{eq:degradation}. The parameters used for the (DC) SBM models ensures that the resulting graphs have a strong community structure.

        \begin{table}[]
        \centering
        \begin{tabular}{m{0.12\textwidth}|m{0.08\textwidth}|m{0.07\textwidth}|m{0.045\textwidth}|m{0.068\textwidth}|m{0.07\textwidth}}
                Graph & $\vert \mathcal{V} \vert$ & $\overline{d}$ & $d_\mathrm{min}$ & $d_\mathrm{max}$ & $q_*$\\
                \hline
                $\eps$-graph & $10000$ & $37.3$ & $5.2$ & $67.5$ & $6.518$ \\
                \hline
                SBM & $10000$ & $39.87$ & $17.7$ & $66.6$ & $21.04$\\
                \hline
                DC-SBM $1$ & $9833.7$ & $33.4$ & $1$ & $906.1$ & $2.634$\\
                \hline
                DC-SBM $2$ & $9932.8$ & $418.8$ & $1$ & $1513.7$ & $4.31$\\
                \hline
                AS CAIDA & $26475$ & $4.032$ & $1$ & $2628$ & $0.4808$ 
                \end{tabular}
            \caption{Experimental parameters associated to each graph. Values are averaged over all $10$ realizations for random graph models.}
            \label{table:graphs}
        \end{table}
        
        \noindent For the $\eps$-graph and the AS-graph, we arbitrarily set $B = 5$ when generating the bandlimited signal $f_\top$. For (DC) SBMs, we take $B = 2$. It is known that the eigenvectors of the \emph{combinatorial} Laplacian $\mathsf{L}$ encode the community structure of these graphs (here,, and taking $B = 2$ here results in signals coherent with the community structure (see Section~\ref{supmat:tikhonov_numerical} of the Supplementary Material for more on the first few eigenvectors of $\mathsf{L}_\theta$).
    \end{paragraph}
        
    \begin{paragraph}{Discussion}
        Let us first comment on the approximation error for the $\eps$-graph (Figure~\ref{subfig:smoothing:eps}, left). Both of our estimators converge linearly in log-log-scale with a mild slope, as expected for Monte-Carlo estimators. Any of the other three conjugate-gradient-based algorithms achieves a better approximation error within $3$ steps than our estimators in $100$ steps (as expected for Monte-Carlo convergence rates), and our 
        methods are not competitive in this setting. However, in this denoising setting, the quality of the denoising is reflected in the reconstruction error (depicted in all the plots of Figure~\ref{fig:smoothing}). Overall, we observe that:
        \begin{enumerate}[label={(\arabic*)}]
            \item $\hat{f}$ consistently performs better than $\overline{f}$, but how much better depends on the graph (see \emph{e.g.} SBM and DC-SBM $1$), and improvements are smaller for graphs with more heterogeneous degree distributions (DC-SBM graphs);
            \item compared with MTSF-based estimators, the CG algorithms perform \emph{worse} on graphs with heterogeneous degree distributions; \label{bp:2}
            \item CG solvers perform worse than our methods when the density increases (MTSF-based methods reach optimal reconstruction error 10 times faster on DC-SBM 2).
        \end{enumerate} 

        \noindent \color{col}{The relatively poor performance of CG on DC-SBM graphs is not unexpected. CG is widely known to be quite sensitive to ill-conditioning~\cite{saad2003iterative}, which is why problem-specific preconditioners are often necessary. For instance, one can show that for a linear system of the form $Ax = b$ with solution $x_*$ (where $A$ is self-adjoint and positive definite), \emph{the time-complexity required by CG} to reach an approximation $\widetilde{x}$ of $x_*$ such that $\frac{\Vert \widetilde{x} - x_* \Vert_A^2}{{\Vert x_* \Vert}_A^2}\leq \eps$ \emph{is bounded in} $\mathcal{O}\left( t_A \sqrt{\kappa} \log\left( \frac{1}{\eps} \right) \right)$, where {$\Vert x \Vert_A^2 = \langle x , A x \rangle$, $t_A$ is the cost of a matrix-vector product involving $A$, and $\kappa$ is the condition number of $A$}~\cite{vishnoi2013lx}. In our case, $\kappa$ is bounded by\footnote{These bounds hold because $1 + d_{\text{max}} \leq \lambda_n \leq 2 d_{\text{max}}$, as we show in the Supplementary Material (Section~\ref{sect:bound_lambda_n}).}
        {\begin{equation}
            \frac{q_* + 1 + d_{\text{max}}}{q_* + \lambda_1} \leq \kappa = \frac{q_* + \lambda_n}{q_* + \lambda_1} \leq \frac{q_* + 2 d_{\text{max}} }{q_* + \lambda_1},
        \end{equation}
        where $\lambda_1 \leq ... \leq \lambda_n$ are the eigenvalues of $\mathsf{L}_\theta$ and $d_{\text{max}}$ is the maximum degree}. Note that $\lambda_1$ is a measure of the quality of the optimal angular assignment (the lower $\lambda_1$ is, the better the quality of the synchronization), which is small for our synthetic connection (good synchronization due to low incoherence).}
\color{col}{In the trivial-connection case, graphs with a highly-heterogenous degree distribution, as is the case here, tend to have poor conditioning; this follows for example from a perturbation argument (see, \emph{e.g.}, Thm.1 in~\cite{zhan2010distributions}).

        \noindent In addition to being sensitive to conditioning, CG is sensitive to graph density: each iteration has cost $t_A$ in $\mathcal{O}(\vert \mathcal{E} \vert)$. On the other hand, we find in practice that our algorithm is
        much less sensitive to conditioning, and less sensitive to density.
        Showing this rigorously is certainly possible but would involves a lengthy exhaustion argument, depending on the comparisons between $q_*$, $\lambda_1$ and $d_{max}$. Note that the concentration bound in Proposition~\ref{prop:variance_reduction} does
        not feature the condition number, which suggests that poorly-conditioned
        graphs do not require more Monte-Carlo replicates. In addition, setting
        $q \approx \overline{d}$ in the runtime bound of Proposition~\ref{prop:bound_complexity} gives $\mathcal{O}(\vert \mathcal{V} \vert)$ expected runtime for sampling MTSFs, which suggests lower sensitivity to graph density. 
        
        In general, we therefore expect MTSF-based solutions to perform better
        than CG on poorly-conditioned or high-density graphs.}
    \end{paragraph}

    Further profiling also shows that the most expensive operations during the computation of $\overline{f}$ (and $\hat{f}$) are the instantiation of the data structures we use, and the successive memory accesses to the different $\theta_e$. These issues will likely benefit from further democratization of randomized computational schemes~\cite{murray2023randomized}, with more efficient optimizations and architectures specialized for these types of computations. Note also that the connection model used here is only representative of a subset of applications (\emph{e.g.} \cite{yu_angular_2012}, \cite{cucuringu_syncrank_2016} or \cite{furutani_gsp_2019}, see also Remark~\ref{rem:gamma}).

    \begin{paragraph}{Choice of $q$}
        \color{col}{In practice, the value of $q_*$ is unknown, and the choice of the regularization parameter $q$ is a classical model-selection issue, that goes beyond the specific case of graphs, and can be approached from different perspectives.
        \begin{itemize}
            \item It may be inferred from an application-specific perspective (this is the case in the vector-field extension problem of~\cite{sharp2019vector}).
            \item From a statistical-signal-processing perspective, it may be estimated by inspecting different criteria, such as (for instance) Akaike's information criterion, the Bayesian information Criterion, leave-one-out cross validation, or Stein's unbiased risk estimator (in the case of graph Tikhonov regularization for trivial connections, this last criterion is discussed in detail in~\cite{pilavci_graph_2021}).
        \end{itemize}
        The methods mentioned in this second point require the computation of the \emph{effective degrees of freedom} of the model~\cite{hastie2009elements}, for different values of $q$. Here, this quantity takes the form
        \begin{equation}
            s(q) = \mathrm{tr}\left( q (\mathsf{L}_\theta + q \mathsf{I})^{-1} \right).
            \label{eq:sq}
        \end{equation}
        \emph{In the trivial-connection setting}, it is known that an unbiased estimator of $s(q)$ is obtained by counting the number of roots of MTSFs: $s(q) = \mathbf{E}_{\phi \sim \mathcal{D}_\mathcal{M}} \left( \vert \phi \cap \mathcal{V} \vert \right)$~\cite{barthelme_inverse_2019}; in conjunction with variance-reduction techniques, this method reaches state-of-the-art performance~\cite{pilavci2022variancebis}. \emph{An important corollary of our work} is the extension of this strategy to \emph{any connection}, as we detail in the Supplementary Material (Remark~\ref{rem:trace_estimation}).}
    \end{paragraph}

    \section{Randomized Angular Synchronization}
    \label{sect:synchro}

    We describe in this section how our novel estimators can be successfully applied to angular synchronization. Recall from Section~\ref{sect:graphs:connections} that, to perform angular synchronization, we aim to minimize the incoherence
    \begin{equation}
        \argmin_{f \in U(\C)^n} \ \langle f , \mathsf{L}_\theta f \rangle,
        \label{eq:synchro}
    \end{equation}
    which is NP-hard in general~\cite{zhang2006complex}. We go over existing approaches before introducing in Section~\ref{subsect:synchro_algo} a randomized schemed based on a {spectral relaxation} of Problem~\eqref{eq:synchro}.

    \begin{paragraph}{Spectral relaxation and other existing approaches}
        Problem~\eqref{eq:synchro} is often relaxed to the following form~\cite{yu_angular_2012,singer2011angular}:
        \begin{equation}
            \argmin_{\Vert f \Vert_2^2 = n} \ \langle f , \mathsf{L}_\theta f \rangle,
            \label{eq:sr}
        \end{equation}
        the solution of which is given by the eigenvector associated to the \emph{smallest} eigenvalue of $\mathsf{L}_\theta$. This solution differs from the exact solution of Problem~\eqref{eq:synchro} (even though guarantees on its quality exist, \emph{e.g.}~\cite{filbir_recovery_2021}), and can be computed using either an inverse power method, a Rayleigh quotient iteration, or a Lanczos iteration~\cite{saad2011numerical}.
        Most theoretical studies focus on Erdös-Rényi graphs, and suggest that spectral relaxations work best in small noise regimes (as compared with other algorithms, regarding the resulting error)~\cite{boumal_nonconvex_2016,perry2018message,lerman2022robust}. 
        Other existing methods include (see \cite{cucuringu2020extension} for a similar discussion):
        \begin{itemize}
            \item Semi-definite relaxations~\cite{singer2011angular}, providing flexible theoretical tools and performance similar to spectral relaxations, but impractical past mid-sized instances ($n \simeq 10^4$).
            \item The \emph{generalized power method} proposed in~\cite{boumal_nonconvex_2016}, provably reaching the optimal solution of Problem~\eqref{eq:synchro} in the presence of (low) Gaussian noise, on complete graphs.
            \item Graph Neural Networks, with state-of-the-art performance for high noise~\cite{he_robust_2023,janco2023unrolled}.
            \item Message-passing algorithms. For Gaussian noise, Approximate Message Passing has been conjectured to be statistically optimal among polynomial time algorithms~\cite{perry2018message}, even at higher noise levels, but is limited to very dense graphs. Cycle-Edge Message-Passing allows exact recovery under a theoretical corruption model (with linear rate), but is more computationally-expensive than spectral relaxations~\cite{lerman2022robust}.
            \item Descent techniques such as~\cite{maunu2023depth}, with exact recovery under a (different) corruption model (see also \cite{liu2023resync}).
        \end{itemize}
        Note that message-passing algorithms and descent techniques (along with others we did not mention as well) can be applied to synchronization problems over more general groups, but are often not competitive on benchmarks for \emph{angular} synchronization~\cite{he_robust_2023}.
    \end{paragraph}

    We will now discuss how to use our randomized estimators for eigenvector-computation.

    \subsection{Proposed Approach}
    \label{subsect:synchro_algo}

    \textcolor{col}{One way to solve the spectral relaxation of Eq.~\eqref{eq:sr} is to perform an inverse power iteration, by setting $f_0 \in \C^\mathcal{V}$ and iterating:
    \begin{equation}
        f_{r+1} = \frac{\mathsf{L}_\theta^{-1} f_r}{\Vert \mathsf{L}_\theta^{-1} f_r \Vert_2}
        \label{eq:ipm}
    \end{equation}
    until convergence. Note that for this iteration to be well-defined, $\mathsf{L}_\theta$ needs to be invertible, \emph{i.e.}, the angular-synchronization problem needs to be non-trivial. Depending on the topology of the graph $\mathsf{L}_\theta$ may be poorly conditioned, yielding a very difficult robust estimation of  $\mathsf{L}_\theta^{-1} f_r$.}

    \textcolor{col}{A classical workaround is to regularize the matrix $\mathsf{L}_\theta$. Noting that the matrices $\mathsf{L}_\theta^{-1}$ and its regularized version $(\mathsf{L}_\theta + q \mathsf{I})^{-1}$ (with a better condition number) share the same eigenvectors, it is a classical result that the power method applied to a matrix of the form $(\mathsf{L}_\theta + q \mathsf{I})^{-1}$ converges to the eigenvalue of $\mathsf{L}_\theta$ closest to $- q$~\cite{saad2011numerical}. In our case $q$ is positive, so that the power method applied to $(\mathsf{L}_\theta + q \mathsf{I})$ converges to the desired spectral solution of angular synchronization, and the iteration reads 
    \begin{equation}
        f_{r+1} = \frac{ (\mathsf{L}_\theta + q \mathsf{I})^{-1} f_r}{\Vert  (\mathsf{L}_\theta + q \mathsf{I})^{-1} f_r \Vert_2}, = \frac{q (\mathsf{L}_\theta + q \mathsf{I})^{-1} f_r}{\Vert q (\mathsf{L}_\theta + q \mathsf{I})^{-1} f_r \Vert_2}.
        \label{eq:power_iteration}
    \end{equation}
    The convergence of this power method iteration is geometric with ratio
    \begin{equation}
        \frac{\mu_1}{\mu_2} = \frac{\lambda_2 + q}{\lambda_1 + q},
    \end{equation}
    where $\mu_1 \leq \mu_2$ (resp. $\lambda_1 \leq \lambda_2$) denote respectively the two largest (resp. smallest) eigenvalues of $q (\mathsf{L}_\theta + q \mathsf{I})^{-1}$ (resp. $\mathsf{L}_\theta$)~\cite{saad2011numerical}. The choice of $q$ is then a classical trade-off between:
    \begin{itemize}
        \item low values of $q$ that induce a smaller $\frac{\lambda_2 + q}{\lambda_1 + q}$ ratio, favoring faster convergence of the power iteration;
        \item large values of $q$ that result in a better conditioning of the system 
    \end{itemize}}

    \begin{rem}
        One could instead maximize $\langle f , \mathsf{A}_\theta f \rangle$ in Eq.~\eqref{eq:synchro}. This formulation is common in theoretical works, that mostly focus on (often dense) Erdös-Rényi graphs, and applying the power method to $\mathsf{A}_\theta$ is efficient  on these graphs. \nt{However, it is often recommended to work with the smallest eigenvalue of $\mathsf{L}_\theta$ rather than the largest of $\mathsf{A}$} when the ratio $\alpha_1/\alpha_2$ of the top eigenvalues of $\mathsf{A}_\theta$ is close to $1$. For trivial connections, such cases include: regular grids (with $\alpha_1/\alpha_2$ growing worse with the dimension), graphs with homogeneous spatial correlations (\emph{e.g.}, $\eps$-graphs, nearest-neighbors-graphs), or graphs with homogeneous degree distributions \emph{and} bottlenecks (poor expansion) such as in SBMs (for regular graphs, this is a classical Cheeger inequality). \color{col}{Our experiments in Section~\ref{sect:A_vs_L} of the Supplementary Material suggest that the same observations can be made for non-trivial connections}\footnote{Further, Cramér-Rao bounds for angular synchronization suggest that angular synchronization is (statistically) difficult on these graphs~\cite{boumal2014cramer}.}. 
        It is also possible to consider a \emph{normalized} Laplacian $\widetilde{\mathsf{L}_\theta} = \mathsf{D}^{-\frac{1}{2}} \mathsf{L}_\theta \mathsf{D}^{-\frac{1}{2}}$ in the spectral relaxation~\eqref{eq:sr}, which was proposed in~\cite{yu_angular_2012,cucuringu2012sensor} and comes with similar guarantees~\cite{bandeira2013cheeger}.
        \label{rem:synchro}
    \end{rem}
    \textcolor{col}{\emph{Our approach} simply consists in estimating each of the successive smoothings {of Eq.~\eqref{eq:power_iteration} with MTSF-based estimators (we in fact recognize the term in $q\left(\mathsf{L}_\theta +q\mathsf{I}\right)^{-1}f_r$ as precisely what we can estimate with our MTSF-based estimators)}. In the remainder of this section, we illustrate the performance of our strategy on a toy problem, and compare it for different graph topologies to a smoothing performed by CG at each step of the power iteration. Further, we also illustrate in Section~\ref{sect:A_vs_L} of the Supplementary Material that, as compared to the power iteration applied to the matrix $\mathsf{A}_\theta$, our strategy can result in very significant speed-ups on graphs for which the top spectral gap of $\mathsf{A}_\theta$ is small (as discussed in Remark~\ref{rem:synchro}).}

    Note that our approach is flexible enough to compute the bottom eigenvector of the \nt{normalized connection Laplacian} $\widetilde{\mathsf{L}_\theta}$ as well (Remark~\ref{rem:synchro}). Indeed, our MTSF-based estimators extend to quantities \nt{of the form} $f_\circ = (\mathsf{L}_\theta + q\mathsf{D})^{-1} (qD) g'$ (see Section~\ref{supmat:fermionic} of the Supplementary Material), which allows to estimate $q (\widetilde{\mathsf{L}_\theta} + q \mathsf{I})^{-1} g$ by setting $g' = \mathsf{D}^{-\frac{1}{2}} g$ and \nt{noting that}:
        \begin{equation}
            q (\widetilde{\mathsf{L}_\theta} + q \mathsf{I})^{-1} g = \mathsf{D}^{\frac{1}{2}} ((\mathsf{L}_\theta + q\mathsf{D})^{-1} (q\mathsf{D})) \ g'.
        \end{equation}

    \subsection{Illustration: Intra-Class Denoising}
    \label{sect:illustration}

    Before delving into precise computation time comparisons with state-of-the-art methods, we first illustrate our proposed method on a toy denoising problem, inspired from an application in cryo-EM~\cite{scheres2012relion,zhao2014rotationally}.

    We consider $n$ copies $(I_i)_{i=1}^n$ of some image $I_*$ that have been rotated and degraded:
    \begin{equation}
        I_i = r_i(I_*) + \eps,
        \label{eq:image_degradation}
    \end{equation}
    with $r_i$ a rotation and $\eps \sim \mathcal{N}(0,\sigma^2 \mathsf{I})$ some Gaussian noise. The rotations $r_i$ are unknown, and the goal is to recover image $I_*$. This is a simplified version of the intra-class denoising problem in cryo-EM, where each (2D) image corresponds to a noisy projection of a (3D) molecule observed in an unknown orientation.

    In the absence of rotations $r_i$, a solution consists in averaging the $I_i$'s, trivially recovering $I_*$ as $n \rightarrow \infty$. This strategy fails in our setting (due to the rotations), and we need to estimate the $r_i$'s before denoising. To do that, we first estimate angles $\theta_{i,j}$ between a subset of pair of images $\mathcal{E}$ (representing the edges of a graph) using image moments~\cite{gonzalez2009digital}\footnote{Specifically, we estimate the orientation of each of the two images using the eigenvectors of the matrix of its covariant moments, and take $\theta_{i,j}$ the angle between these two sets of eigenvectors}, and perform angular synchronization to estimate the rotations $r_i$'s. We then rotate and average the images $I_i$'s. 
    
\begin{figure}[]
	\begin{subfigure}{0.55\linewidth}
		\begin{subfigure}{0.25\linewidth}
			\centering
			\includegraphics[scale=0.2]{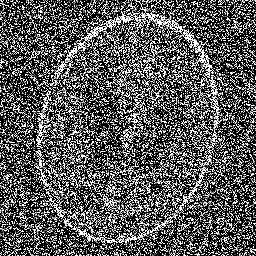}
			\caption{$\sigma^2 = 1$}
		\end{subfigure}
		\hfill
		\begin{subfigure}{0.32\linewidth}
			\centering
			\includegraphics[scale=0.2]{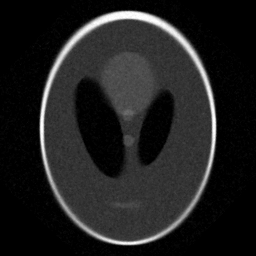}
			\caption{$\sigma^2 = 1$, exact}
		\end{subfigure}
		\hfill
		\begin{subfigure}{0.32\linewidth}
			\centering
			\includegraphics[scale=0.2]{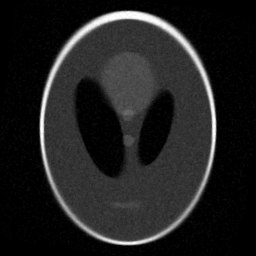}
			\caption{$\sigma^2 = 1$, MTSF}
		\end{subfigure}
		
		\begin{subfigure}{0.25\linewidth}
			\centering
			\includegraphics[scale=0.2]{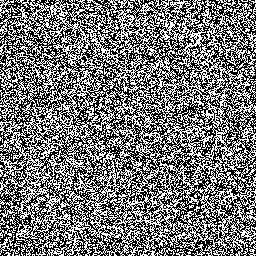}
			\caption{$\sigma^2 = 5$}
		\end{subfigure}
		\hfill
		\begin{subfigure}{0.32\linewidth}
			\centering
			\includegraphics[scale=0.2]{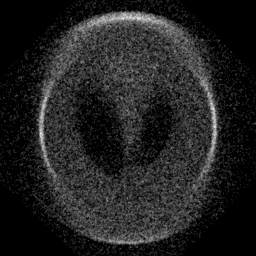}
			\caption{$\sigma^2 = 5$, exact}
		\end{subfigure}
		\hfill
		\begin{subfigure}{0.32\linewidth}
			\centering
			\includegraphics[scale=0.2]{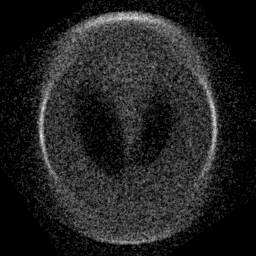}
			\caption{$\sigma^2 = 5$, MTSF}
		\end{subfigure}
	\end{subfigure}
	\begin{subfigure}{0.44\linewidth}
		\centering
		\scalebox{0.35}{

\begin{tikzpicture}[/tikz/background rectangle/.style={fill={rgb,1:red,1.0;green,1.0;blue,1.0}, fill opacity={1.0}, draw opacity={1.0}}, show background rectangle]
\begin{axis}[point meta max={nan}, point meta min={nan}, legend cell align={left}, legend columns={1}, title={}, title style={at={{(0.5,1)}}, anchor={south}, font={{\fontsize{14 pt}{18.2 pt}\selectfont}}, color={rgb,1:red,0.0;green,0.0;blue,0.0}, draw opacity={1.0}, rotate={0.0}, align={center}}, legend style={color={rgb,1:red,0.0;green,0.0;blue,0.0}, draw opacity={1.0}, line width={1}, solid, fill={rgb,1:red,1.0;green,1.0;blue,1.0}, fill opacity={1.0}, text opacity={1.0}, font={{\fontsize{24 pt}{10.4 pt}\selectfont}}, text={rgb,1:red,0.0;green,0.0;blue,0.0}, cells={anchor={center}}, at={(0.98, 0.98)}, anchor={north east}}, axis background/.style={fill={rgb,1:red,1.0;green,1.0;blue,1.0}, opacity={1.0}}, anchor={north west}, xshift={1.0mm}, yshift={-1.0mm}, width={150.4mm}, height={99.6mm}, scaled x ticks={false}, xlabel={m}, x tick style={color={rgb,1:red,0.0;green,0.0;blue,0.0}, opacity={1.0}}, x tick label style={color={rgb,1:red,0.0;green,0.0;blue,0.0}, opacity={1.0}, rotate={0}}, xlabel style={at={(ticklabel cs:0.5)}, anchor=near ticklabel, at={{(ticklabel cs:0.5)}}, anchor={near ticklabel}, font={{\fontsize{26 pt}{14.3 pt}\selectfont}}, color={rgb,1:red,0.0;green,0.0;blue,0.0}, draw opacity={1.0}, rotate={0.0}}, xmajorgrids={true}, xmin={0.7299999999999995}, xmax={10.27}, xticklabels={{$1$,$2$,$3$,$4$,$5$,$6$,$7$,$8$,$9$,$10$}}, xtick={{1.0,2.0,3.0,4.0,5.0,6.0,7.0,8.0,9.0,10.0}}, xtick align={inside}, xticklabel style={font={{\fontsize{24pt}{10.4 pt}\selectfont}}, color={rgb,1:red,0.0;green,0.0;blue,0.0}, draw opacity={1.0}, rotate={0.0}}, x grid style={color={rgb,1:red,0.0;green,0.0;blue,0.0}, draw opacity={0.1}, line width={0.5}, solid}, axis x line*={left}, x axis line style={color={rgb,1:red,0.0;green,0.0;blue,0.0}, draw opacity={1.0}, line width={1}, solid}, scaled y ticks={false}, ylabel={Image Reconstruction Error}, y tick style={color={rgb,1:red,0.0;green,0.0;blue,0.0}, opacity={1.0}}, y tick label style={color={rgb,1:red,0.0;green,0.0;blue,0.0}, opacity={1.0}, rotate={0}}, ylabel style={at={(ticklabel cs:0.5)}, anchor=near ticklabel, at={{(ticklabel cs:0.5)}}, anchor={near ticklabel}, font={{\fontsize{26 pt}{14.3 pt}\selectfont}}, color={rgb,1:red,0.0;green,0.0;blue,0.0}, draw opacity={1.0}, rotate={0.0}}, ymajorgrids={true}, ymin={40.130096301235625}, ymax={51.48631343031269}, yticklabels={{$42$,$44$,$46$,$48$,$50$}}, ytick={{42.0,44.0,46.0,48.0,50.0}}, ytick align={inside}, yticklabel style={font={{\fontsize{24 pt}{10.4 pt}\selectfont}}, color={rgb,1:red,0.0;green,0.0;blue,0.0}, draw opacity={1.0}, rotate={0.0}}, y grid style={color={rgb,1:red,0.0;green,0.0;blue,0.0}, draw opacity={0.1}, line width={0.5}, solid}, axis y line*={left}, y axis line style={color={rgb,1:red,0.0;green,0.0;blue,0.0}, draw opacity={1.0}, line width={1}, solid}, colorbar={false}]
    \addplot[color={rgb,1:red,0.0;green,0.6056;blue,0.9787}, name path={1291ddd6-d984-4a4b-b473-4c24a6273d1b}, draw opacity={1.0}, line width={2}, solid, mark={*}, mark size={3.0 pt}, mark repeat={1}, mark options={color={rgb,1:red,0.0;green,0.6056;blue,0.9787}, draw opacity={1.0}, fill={rgb,1:red,0.0;green,0.6056;blue,0.9787}, fill opacity={1.0}, line width={0.75}, rotate={0}, solid}]
        table[row sep={\\}]
        {
            \\
            1.0  51.164911058735036  \\
            2.0  45.32179551097983  \\
            3.0  42.324470662923886  \\
            4.0  41.507491200527745  \\
            5.0  40.830849233588545  \\
            6.0  40.80938130702225  \\
            7.0  40.622668978395915  \\
            8.0  40.62200570360175  \\
            9.0  40.63847911303608  \\
            10.0  40.703173970567235  \\
        }
        ;
    \addlegendentry {MTSF}
    \addplot[color={rgb,1:red,0.9308;green,0.3675;blue,0.5758}, name path={3a606e8d-2c74-42a7-a3c0-19eeea1c20e7}, draw opacity={1.0}, line width={2}, solid]
        table[row sep={\\}]
        {
            \\
            -8.809999999999999  40.45149867281327  \\
            19.81  40.45149867281327  \\
        }
        ;
    \addlegendentry {Exact}
    \addplot[color={rgb,1:red,0.9308;green,0.3675;blue,0.5758}, name path={3a606e8d-2c74-42a7-a3c0-19eeea1c20e7}, draw opacity={1.0}, line width={2}, solid, forget plot]
        table[row sep={\\}]
        {
            \\
            -8.809999999999999  40.45149867281327  \\
            19.81  40.45149867281327  \\
        }
        ;
    \addplot[color={rgb,1:red,0.9308;green,0.3675;blue,0.5758}, name path={3a606e8d-2c74-42a7-a3c0-19eeea1c20e7}, draw opacity={1.0}, line width={2}, solid, forget plot]
        table[row sep={\\}]
        {
            \\
            -8.809999999999999  40.45149867281327  \\
            19.81  40.45149867281327  \\
        }
        ;
    \addplot[color={rgb,1:red,0.9308;green,0.3675;blue,0.5758}, name path={3a606e8d-2c74-42a7-a3c0-19eeea1c20e7}, draw opacity={1.0}, line width={2}, solid, forget plot]
        table[row sep={\\}]
        {
            \\
            -8.809999999999999  40.45149867281327  \\
            19.81  40.45149867281327  \\
        }
        ;
    \addplot[color={rgb,1:red,0.9308;green,0.3675;blue,0.5758}, name path={3a606e8d-2c74-42a7-a3c0-19eeea1c20e7}, draw opacity={1.0}, line width={2}, solid, forget plot]
        table[row sep={\\}]
        {
            \\
            -8.809999999999999  40.45149867281327  \\
            19.81  40.45149867281327  \\
        }
        ;
    \addplot[color={rgb,1:red,0.9308;green,0.3675;blue,0.5758}, name path={3a606e8d-2c74-42a7-a3c0-19eeea1c20e7}, draw opacity={1.0}, line width={2}, solid, forget plot]
        table[row sep={\\}]
        {
            \\
            -8.809999999999999  40.45149867281327  \\
            19.81  40.45149867281327  \\
        }
        ;
    \addplot[color={rgb,1:red,0.9308;green,0.3675;blue,0.5758}, name path={3a606e8d-2c74-42a7-a3c0-19eeea1c20e7}, draw opacity={1.0}, line width={2}, solid, forget plot]
        table[row sep={\\}]
        {
            \\
            -8.809999999999999  40.45149867281327  \\
            19.81  40.45149867281327  \\
        }
        ;
    \addplot[color={rgb,1:red,0.9308;green,0.3675;blue,0.5758}, name path={3a606e8d-2c74-42a7-a3c0-19eeea1c20e7}, draw opacity={1.0}, line width={2}, solid, forget plot]
        table[row sep={\\}]
        {
            \\
            -8.809999999999999  40.45149867281327  \\
            19.81  40.45149867281327  \\
        }
        ;
    \addplot[color={rgb,1:red,0.9308;green,0.3675;blue,0.5758}, name path={3a606e8d-2c74-42a7-a3c0-19eeea1c20e7}, draw opacity={1.0}, line width={2}, solid, forget plot]
        table[row sep={\\}]
        {
            \\
            -8.809999999999999  40.45149867281327  \\
            19.81  40.45149867281327  \\
        }
        ;
    \addplot[color={rgb,1:red,0.9308;green,0.3675;blue,0.5758}, name path={3a606e8d-2c74-42a7-a3c0-19eeea1c20e7}, draw opacity={1.0}, line width={2}, solid, forget plot]
        table[row sep={\\}]
        {
            \\
            -8.809999999999999  40.45149867281327  \\
            19.81  40.45149867281327  \\
        }
        ;
\end{axis}
\end{tikzpicture}}
		\caption{Error for $\sigma^2 = 5$.}
	\end{subfigure}        
	\caption{Examples of recovery from noisy images (left), with different noise levels, using both the exact solution to the spectral relaxation~\eqref{eq:sr} and the approach from Section~\ref{subsect:synchro_algo}.}
	\label{fig:phantom}
\end{figure}
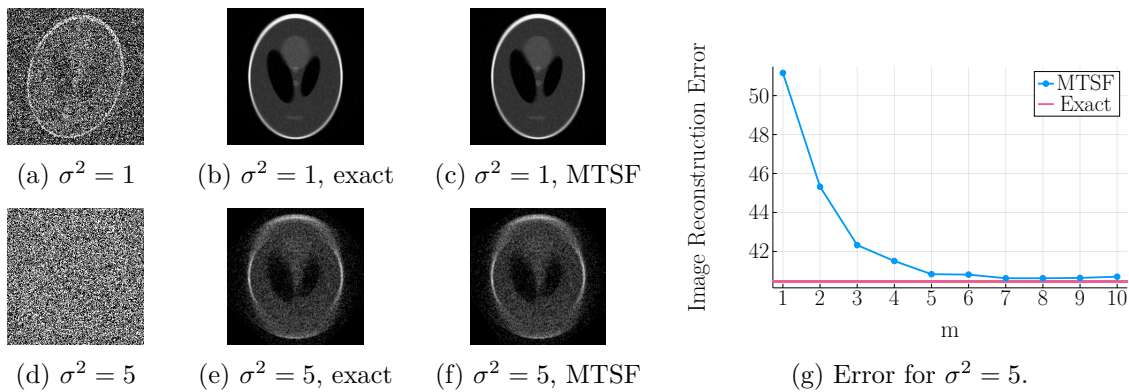

    We take $I_*$ the $256 \times 256$ Shepp-Logan phantom~\cite{shepp1974fourier}, $n = 1000$, $q = 10^{-3}$, $k=20$ iterations of the power method, and randomly choose the underlying graph $\mathcal{G} \sim ER(n,p)$ with $p = \frac{5}{n}$. We uniformly sample rotations $r_i$ with angles in $(-\frac{\pi}{2},\frac{\pi}{2})$\footnote{We restrict the possible rotations so that we can use the simple image-moment-based registration, other methods could be used instead.}. We display the recovered images in Figure~\ref{fig:phantom} (solving the spectral relaxation~\eqref{eq:sr} using both our method, with $m = 10$ MTSFs for the estimator $\overline{f}$ at each step, and an exact solver), for different noise levels $\sigma^2 \in \{1,3,5,10\}$. We also plot the image reconstruction error $e_{I}(x) = \Vert I_* - x \Vert_2$ for $\sigma^2 = 5$ with varying values of $m$. \textcolor{col}{Here, we have no guarantee that the connection resulting from the image-moment-based estimation is weakly-incoherent, and \emph{perform importance sampling} according to the importance distribution of Eq.~\eqref{eq:importance}\footnote{Here, there is no issue with estimating $a f_*$ instead of $f_*$, due to the re-normalization in Eq.~\eqref{eq:power_iteration}.}, based on the estimator $\overline{f}$.}

    Even though the connection may not be weakly-incoherent, the MTSF-based synchronization allows to obtain results of a quality very similar to exact synchronization on this problem, even at higher noise levels. 
    
    \textcolor{col}{We propose in Section~\ref{sect:mnist} of the Supplementary Material a similar application on the MNIST dataset, where we do not control the incoherence of the connection, nor the noise.}
    
    We will now investigate trade-offs between the quality of the synchronization and the runtime of our methods depending on the solver used to compute the power-method iteration in Eq.~\eqref{eq:power_iteration}, and the topology of the graph.
    
    \subsection{Numerical Evaluation}
    \label{sect:synchro:eval}

    We compare the performance of the method from Section~\ref{subsect:synchro_algo} on synthetic graph data using four different iterative solvers: $\overline{f}$, $\hat{f}$ and a conjugate-gradient descent with and without diagonal preconditioning.

\begin{figure}[]
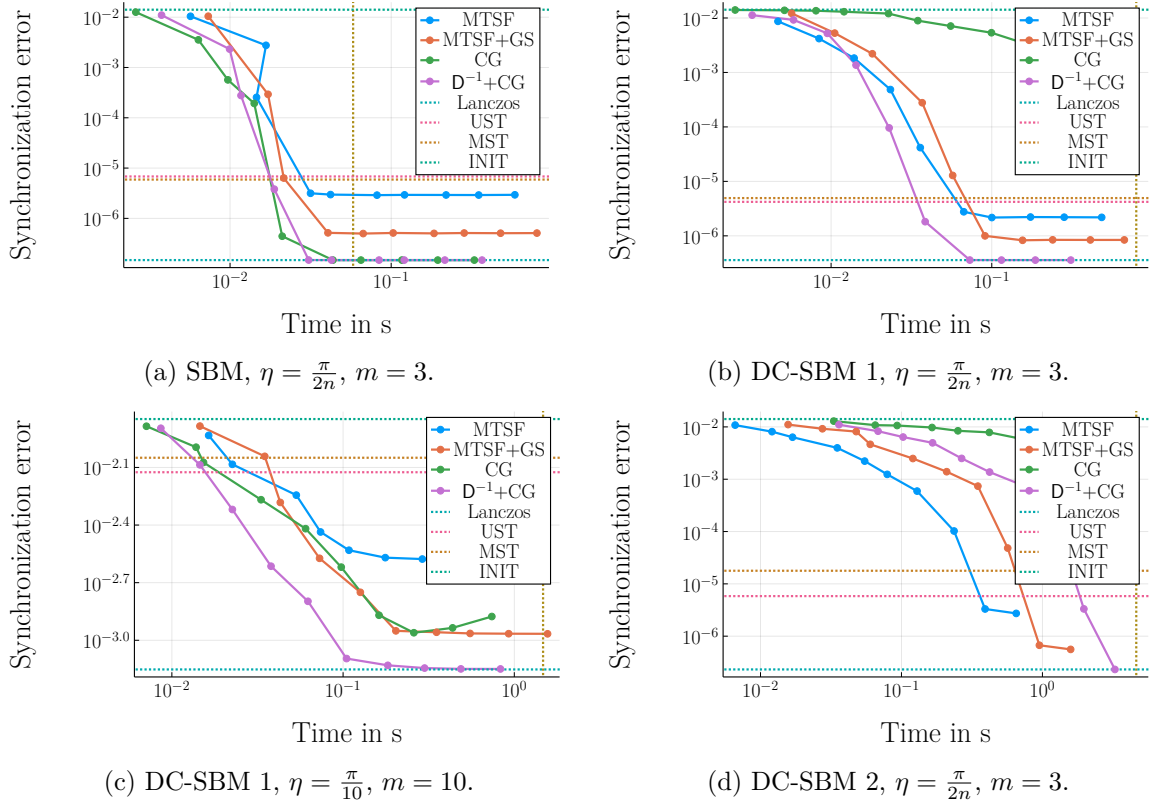

        \centering
        
        \begin{subfigure}{0.49\textwidth}
            \scalebox{0.42}{\input{synchro_sbm_m_3.tikz}}
            \caption{SBM, $\eta = \frac{\pi}{2n}$, $m = 3$.}
        \end{subfigure}
        \hfill
        \begin{subfigure}{0.49\textwidth}
            \scalebox{0.42}{\input{synchro_dc-sbm-1_m_3.tikz}}
            \caption{DC-SBM $1$, $\eta = \frac{\pi}{2n}$, $m = 3$.}
        \end{subfigure}
        
        \begin{subfigure}{0.49\textwidth}
            \scalebox{0.42}{\input{synchro_dc-sbm-1_m_10.tikz}}
            \caption{DC-SBM $1$, $\eta = \frac{\pi}{10}$, $m = 10$.}
        \end{subfigure}
        \hfill
        \begin{subfigure}{0.49\textwidth}
            \scalebox{0.42}{\input{synchro_dc-sbm-2_m_3.tikz}}
            \caption{DC-SBM $2$, $\eta = \frac{\pi}{2n}$, $m = 3$.}
        \end{subfigure}
        
        \caption{Runtime-precision trade-offs for angular synchronization. Each data point corresponds to a value of $k$. The vertical yellow line records the average runtime of the Lanczos-based computation, the horizontal blue line the synchronization error for the Lanczos-based computation. The horizontal red (resp. orange) line is the mean synchronization error for the UST-based (resp. MST-based) synchronization, and the horizontal teal line the error for the random initialization.}
        \label{fig:synchro}
    \end{figure}

    \begin{paragraph}{Setup}
    We work with the SBM and DC-SBM graph models discussed in Section~\ref{subsect:runtime_precision_smoothing}. We endow each graph with a random connection 
    generated according to the model described in Section~\ref{sect:numerical}, for both \nt{weakly-incoherent} ($\eta = \frac{\pi}{2 n}$) and incoherent connections ($\eta =  \frac{\pi}{10}$), and aim to recover the vector $x$ with  $x_i = e^{\iota \omega_i}$.
    Starting from $f_0$ taken uniformly in $U(\C)^n$, we perform $k$ iterations of the power iteration of Eq.~\eqref{eq:power_iteration}, estimating the solution with $m$ MTSFs (resp. conjugate-gradient iterations) for each method. We use $m = 3$ for \nt{weakly-incoherent} connections, and $m = 10$ for incoherent ones. For incoherent connections, we use the importance-sampling distribution from Eq.~\eqref{eq:importance}. We average the synchronization errors
    \begin{equation}
        e_s(f) = \min_{r \in U(\C)} \ \frac{\Vert f - rx \Vert_2}{n}
        \label{eq:sync_error}
    \end{equation}
    over $20$ executions, and measure the mean runtimes over $100$ runs. We set $q = 10^{-2} \times \overline{d}$ and take measurements for different values of $k$ ($10$ logarithmically-spaced values in $\{1,...,100\}$). The results are averaged over $10$ realizations of the random graphs. See Figure~\ref{fig:synchro}.
    \\
    We also take measurements for a Lanczos-iteration-based computation\footnote{We use the one implemented in \url{https://jutho.github.io/KrylovKit.jl/stable/man/eig/\#KrylovKit.eigsolve}.} (for the matrix $\mathsf{L}_\theta$), and for a naive synchronization algorithm going as follows. First, fix a root node $v \in V$, then:
    \begin{itemize}
        \item sample a spanning tree of $G$ uniformly (a UST, using Wilson's algorithm~\cite{wilson_random_1996}),
        \item propagate the value from $v$ to the other nodes (taking into account the offsets).
    \end{itemize}
    \noindent Note that this procedure does \emph{not} depend on $k$ or $m$. Finally, we also evaluate a similar strategy propagating along a \emph{maximum spanning tree} (MST), with respect to the edge-weights $w_{i,j} = \vert \cos(\theta_{i,j}) \vert$ (a strategy inspired from surface reconstruction techniques such as~\cite{hoppe1992surface}).
    \end{paragraph}

    \begin{paragraph}{Comments}
    Results vary greatly depending on the graph and connection. For \nt{weakly-incoherent} connections ($\eta = \frac{\pi}{2 n}$), we observe the following:
    \begin{enumerate}[label={(\arabic*)}]
        \item Krylov-subspaces-based methods are sensitive to the conditioning : all methods perform similarly on the SBM, but not on the (less well-conditioned) DC-SBM graphs.
        \item Higher density results in significant slow-downs for Krylov-subspaces methods, and MTSF-based iterations have a much better performance in this case(DC-SBM 2).
		\item For weakly-inconsistent connections, the tree-based methods also achieve good synchronization quality, but this is no longer the case for incoherent connections.
    \end{enumerate}
	If the first two factors combine (DC-SBM 2), our methods offer a ``10 times'' speed-up.
    \end{paragraph}

    \begin{rem}
        We also experimented with two other estimation strategies (not shown): the spectral relaxation using the normalized Laplacian $\widetilde{\mathsf{L}_\theta}$, and a generalized-power-method-like algorithm, replacing the global normalization in Eq.~\eqref{eq:power_iteration} by a component-wise normalization, so that each $f_k$ belongs to $U(\C)^n$ (like in~\cite{boumal_nonconvex_2016}).
        In our experiments, we did not observe any qualitative differences with the method described in Section~\ref{subsect:synchro_algo}.
    \end{rem}

    For spectral relaxations, our results suggest that MTSF-based solutions perform no-worse than deterministic solutions, and get comparatively better as the density of the graph increases. These techniques are mostly suitable for approximations, especially since previous studies (\emph{e.g.}~\cite{perry2018message}) suggest that spectral-relaxation-based techniques mainly offer good-quality solutions at \emph{low} noise-level, but are otherwise \emph{sub-optimal}. One drawback of our strategy is that it requires to set a value for $q$: further refinements may include adapting our methodology to Rayleigh-quotient iterations~\cite{trefethen2022numerical}, such that $q$ no longer needs to be fixed, and that better convergence rates may be obtained.
    
    \section{Conclusion and Perspectives}

    We propose MTSF-based estimators for connection-aware smoothing and angular synchronization, reaching competitive performance on both problems, especially in low-precision regimes, and even \emph{without} implementing parallelization. Our estimators perform propagation along the branches of MTSFs and, compared to standard deterministic methods, have another noteworthy advantage: they are less sensitive to both the density of the graph \emph{and} the conditioning of the system considered (which is typically bad for graphs with broad degree distributions). 
    Our techniques can apply to a variety of other problems: one result we did not include is the extension of the smoothing estimator to the interpolation problem (\emph{i.e.}, extending a signal that is known only on a few nodes to all the graph), which may find applications to \emph{e.g.} synchronization in presence of anchor nodes. 
    
    \textcolor{col}{Our work opens a number of research directions.  On the theoretical side: i)~{can  the weak-inconsistency condition (Condition~\ref{cond:weak_inconsistency}) be lifted without resorting to importance sampling strategies? We include a few preliminary developments in this direction in  Section~\ref{sect:extension} of the Supplementary Material.} ii)~Can similar estimators be developed for $O(\R^d)$ synchronization? iii)~Instead of using a fixed number $m$ of MTSFs  for the estimation (or CG iterations) of each of $k$ iterations of the power method, how should those $km$ MTSFs be used the to achieve maximal precision? Alternatively, we could re-use the set of $m$ MTSFs for each of the $k$ iterations of the power method.}
    
    From a more applied perspective, many applications of angular synchronization only consider the \emph{exact} solution of the spectral relaxation, but how much loss in precision (and gain in speed) is actually acceptable? Our estimators can be seamlessly implemented on distributed systems (with communication complexity of the order of the number of steps in the sampling algorithm): can this be leveraged in applications? 
    Can random decompositions such as MTSFs be useful in other inference algorithms, \emph{e.g.}, message-passing algorithms~\cite{mezard2009information,lerman2022robust}?
	
	\bibliographystyle{siam}
	\bibliography{ref}

    \newpage

    \begin{center}{\huge \color{header1}\textbf{Supplementary Material}}\end{center}

    \vspace{0.5cm}

    We work with a weighted graph $\mathcal{G}$ in all the following proofs, with $w_e \in \R_+$ the weight of edge $e$. We still denote by $\mathsf{L}_\theta$ the resulting (weighted) connection Laplacian, with off-diagonal entries $(\mathsf{L}_\theta)_{i,j} = - w_{i,j} e^{\iota \theta_{(j,i)}}$ if $\{i,j\}\in \mathcal{E}$ ($(\mathsf{L}_\theta)_{i,j} = w_{(i,j)} =0$ otherwise), and diagonal $(\mathsf{L}_\theta)_{i,i} = d_i = \sum_j w_{i,j}$ the weighted degree of node $i$.

    \section{Proof of Proposition~\ref{prop:tikhonov_solution}}
    \label{supmat:tikhonov_solution}

    Denote by $C: \C^\mathcal{V} \rightarrow \R$ the cost function 
    \begin{equation}
        C(f) = q \Vert f - g \Vert_2^2 + \langle f , \mathsf{L}_\theta f \rangle.
    \end{equation}
    We use a standard argument from $\C \R$-calculus, and seek the zeros of the Fréchet Wirtinger derivatives of $C$ (see \emph{e.g.}~\cite{bouboulis2010wirtinger}). Let us first compute these derivatives: for $f,h \in \C^\mathcal{V}$, we have
    \begin{align}
        C(f + h) & = q \langle (f + h) - g, (f + h) - g \rangle + \langle (f+h), \mathsf{L}_\theta (f+h) \rangle, \\
        & = C(f) + D_f^{W}C(h) + D_f^{W^*}C(h) + \left( q \langle h,h \rangle + \langle h, \mathsf{L}_\theta h \rangle \right),
    \end{align}
    where $(h \mapsto q \langle h,h \rangle + \langle h, \mathsf{L}_\theta h \rangle) \in o(h)$, and the associated (conjugate) Fréchet-Wirtinger derivatives of $C$ (at $f$) $D_f^W C$ and $D_f^{W^*} C: \C^\mathcal{V} \rightarrow \C$ are given by:
    \begin{align}
        D_f^W C(h) = \langle q(f-g) + \mathsf{L}_\theta, h \rangle \\
        D_f^{W^*} C(h) = \langle h, q (f-g) + \mathsf{L}_\theta \rangle.
    \end{align}
    The result follows from the fact that $D_f^W C = 0$ (resp. $D_f^{W^*} C = 0$) identically if and only if $f = q (\mathsf{L}_\theta + q \mathsf{I})^{-1} g$.
    
    \section{A Generalization of Property~\ref{enum:2}}
    \label{supmat:random_walk_connection_proof}
    
    Define $\Delta_\theta = \mathsf{D}^{-1} \mathsf{L}_\theta$ the weighted connection-aware random walk Laplacian $\Delta_\theta$. Proposition~\ref{prop:random_walk_connection_bis} generalizes property~\ref{enum:2}.
    
    \begin{prop}
    For all $l \geq 1$, we have
        \begin{equation}
            (\mathsf{I} - \Delta_\theta)^l_{i,j} = \sum_{\substack{{p \in P_i^j} \\ {l(p) = l}}} \left( \prod_{0 \leq k < l} \frac{w_{(u_k,u_{k+1})}}{d_{u_k}} \right) \psi_{p^*},
            \label{eq:random_walk_laplacian_connection_bis}
        \end{equation}
        with $l(p)$ the length of a path $p = ((u_0,u_1),(u_1,u_2),...,(u_{l(p)-1},u_{l(p)}))$.
        \label{prop:random_walk_connection_bis}
    \end{prop}
    
    For each path $p$, the product in the r.h.s. of Equation~\eqref{eq:random_walk_laplacian_connection_bis} is the probability of observing $p$ when performing a random walk from $i$ to $j$ on $G$. 
    For the trivial connection this expression is equal to $(\mathsf{D}^{-1} \mathsf{A})^l_{i,j}$, the probability of a path going from $i$ to $j$ in $l$ steps, but Proposition~\ref{prop:random_walk_connection_bis} more broadly captures \emph{random propagations} along the corresponding sampled paths. 
    \\
    Observations similar to Proposition~\ref{prop:random_walk_connection_bis} have been laid out in a few works, such as \cite{kassel2021covariant} for continuous time random walk, and in a number of other situations, where the connection often stems from theoretical physics considerations~\cite{brydges1979construction}.
    
    \begin{proof}

    The result is obvious for $l = 1$. We proceed by induction and assume the result true for some $l \in \N$. Then, recalling that $(\mathsf{I} - \Delta_\theta)_{v,j} = \frac{w_{v,j}}{d_v}$, we have:
    \begin{align}
        (\mathsf{I} - \Delta_\theta)^{l+1}_{i,j} & = \left( (\mathsf{I} - \Delta_\theta)^{l} (\mathsf{I} - \Delta_\theta) \right)_{i,j} \\
        & = \sum_{v \in \mathcal{V}} \sum_{\substack{{p \in P_i^v} \\ {l(p) = l}}} \left( \prod_{0 \leq k < l} \frac{w_{(u_k,u_{k+1})}}{d_{u_k}} \right) \frac{w_{v,j}}{d_v} \psi_{p^*}  \psi_{(v,j)^*} \\
        & = \sum_{\substack{{p \in P_i^j} \\ {l(p) = l + 1}}} \left( \prod_{0 \leq k < l + 1} \frac{w_{(u_k,u_{k+1})}}{d_{u_k}} \right) \psi_{p^*}.
    \end{align}
    
    \end{proof}

    \begin{rem}
        For $l = 0$, we have $(\mathsf{I} - \Delta_\theta)^l = \mathsf{I}$.
    \end{rem}

    \section{Proof of Proposition~\ref{prop:feynmac_kac}}
    \label{supmat:feynman_kac}

    We will show an extension of Proposition~\ref{prop:feynmac_kac} to heterogeneous values of $q$. That is, we associate to each node $i$ some non-negative weight $q_i \in \R$ (at least one of them needs to be \emph{positive}), that we record in the diagonal matrix $\mathsf{Q}$ with $\mathsf{Q}_{i,i} = q_i$. These $q_i$'s and can be interpreted as edge weights in the extended graph $\mathcal{G}^\Gamma$, and we extend the definition of the measure $\nu_i$ over $P_i^\Gamma$ to account for these weights, so that
    \begin{equation}
        \prb_{\nu_i}(p) = \mathbf{1}_{\{i\}}(u_0) \mathbf{1}_{\{ \Gamma \}}(u_l) \prod_{0 \leq k < l} \left( \mathbf{1}_{\mathcal{V}\setminus \{\Gamma\}(u_k) } \frac{w_{(u_k,u_{k+1})}}{d_{u_k} + q_{u_k}} \right)
    \end{equation}
    for a path $p = ((u_0,u_1),...,(u_{l-1},u_l))$ of length $l$, and $\mathbf{1}_S(x)$ the indicator that $x \in S$.
    \noindent We show the following Feynman-Kac formula, on weighted graphs and for heterogeneous $q_i$'s.
    
    \vspace{0.2cm}

    \begin{prop} For $j$ the last node reached before absorption of a path $p$, we have
        \begin{equation}
            \left((\mathsf{L}_\theta + \mathsf{Q})^{-1} \mathsf{Q} g \right)(i) = \mathbf{E}_{p \sim \nu_i} \left(\psi_{p^*_\Gamma(g_j)}\right). 
        \end{equation}
        \label{prop:feynman_kac_bis}
    \end{prop}

    Note that the l.h.s. is a more general smoothing operation than the solution to the Tikhonov problem, that can be leveraged in the presence of heteroscedastic noise (\emph{e.g.} for the problem considered in Section~\ref{subsect:runtime_precision_smoothing} with different values of $\sigma^2$ for each node).

    \begin{proof}
    Consider the connection-aware random-walk Laplacian $\Delta_\theta^\mathsf{Q}$ associated to the extended graph $\mathcal{G}^\Gamma$. Its restriction to the rows and columns indexed by $\mathcal{V}$ reads $(\Delta_\theta^q )_\mathcal{V} = (\mathsf{D} + \mathsf{Q})^{-1} (\mathsf{L}_\theta + \mathsf{Q})$. It is clearly invertible, and all its eigenvalues lie in $(0,2)$ (for instance by the Gershgorin circle theorem), hence the eigenvalues of $(\mathsf{I} - (\Delta^\mathsf{Q}_\theta)_\mathcal{V})$ are in $(-1,1)$, ensuring the convergence of the right-hand-side expression in:
        \begin{align}
            \label{eq:proof_feynmac_kac_1}
            \left((\Delta_\theta^\mathsf{Q} )_\mathcal{V}^{-1}\right)_{i,j} & = \left[ \sum_{l \geq 0} \left(\mathsf{I} - (\Delta^q_\theta)_\mathcal{V} \right)^l \right]_{i,j} \\
            & = \sum_{p \in P_i^j} \left( \prod_{0 \leq k < l(p)} \frac{w_{(u_k,u_{k+1})}}{d_{u_{k}} + q_{u_k}} \right) \psi_{p^*}, \label{eq:fk_point}
        \end{align}
        where Equation~\eqref{eq:fk_point} follows from applying Proposition~\ref{prop:random_walk_connection_bis} to $\mathcal{G}^\Gamma$.
        \\
        The result follows by right multiplication with $(\mathsf{D} + \mathsf{Q})^{-1} \mathsf{Q}$, where $\left( (\mathsf{D} + \mathsf{Q})^{-1} \mathsf{Q} \right)_{j,j} = \frac{q_j}{d_{u_{l-1}} + q_j}$ accounts for the probability of the last transition from $u_{l-1} = j$ to $\Gamma$. 
    \end{proof}

    \section{Proof of Theorem~\ref{th:fermionic}}
    \label{supmat:fermionic}

    Before presenting the proof of our main Theorem, we state it for weighted graphs with heterogeneous $q_i$'s.

    \begin{paragraph}{Statement in the presence of weights} 
        Consider the distribution $\mathcal{D}_\mathcal{M}$ over $\mathcal{M}(\mathcal{G})$ such that
        \begin{equation}
            \prb_{\mathcal{D}_\mathcal{M}}(\phi) \propto \prod_{r \in \phi \cap \mathcal{V}} q_r \prod_{e \in \phi \cap \mathcal{E}} w_e \prod_{C \in \mathcal{C}(\phi)} \left((2 - 2 \cos(\theta_C)\right),
        \end{equation}
        and the estimator $\tilde{f}(i,\phi,g) = \psi_{r_\phi(i) \xrightarrow[]{\phi} i} \left( g(r_\phi(i)) \right)$ propagating the value from $r_\phi(i)$ to $i$ if $i$ belongs to a rooted tree, and $\tilde{f}(i,\phi,g) = 0$ if $i$ lies in a unicycle. 
        \begin{thm}
        \textcolor{col}{Letting $\Vert f \Vert^2_\mathsf{Q} = \sum_{i \in \mathcal{V}} q_i \vert f(i) \vert^2$, we have both:}
        \begin{equation}
            \left( (\mathsf{L}_\theta + \mathsf{Q})^{-1} \mathsf{Q} g \right) (i) = \mathbf{E}_{\mathcal{D}_\mathcal{M}}(\tilde{f}(i,\phi,g)),
            \label{eq:heterogeneous_fermionic}
        \end{equation}
        \begin{equation}
            \color{col}{\mathbf{E}_{\phi \sim \mathcal{V}} \left( \Vert \tilde{f}(\phi,g) - f_* \Vert^2_\mathsf{Q} \right) = \Vert g \Vert_\mathsf{Q}^2 - \Vert f_* \Vert^2_\mathsf{Q}.}
        \end{equation}
        \label{th:fermionic_bis}
        \end{thm}
        We now proceed with the proof.
    \end{paragraph}

    \begin{paragraph}{Determinantal Point Processes}
        The first step in our proof consists in re-stating the definition of $\mathcal{D}_\mathcal{M}$ as a Determinantal Point Process (DPP), which will provide us with powerful tools for reasoning about MTSFs. A (discrete) DPP over a finite set $\mathcal{X}$ associates a probability to each subset $X \subseteq \mathcal{X}$, and is defined by its marginal probabilities. It is parametrized by an Hermitian matrix $K \in M_{\vert \mathcal{X} \vert} (\C)$, whose eigenvalues must all lie in $[0,1]$, known as the \emph{marginal kernel} of the DPP, denoted $\mathrm{DPP}(K)$. Precisely, we say that $X \sim \mathrm{DPP}(K)$ if 
        \begin{equation}
            \prb_{\mathrm{DPP}(K)}(A \subseteq X) = \det(K)_{A,A}
        \end{equation}
        for all $A \subseteq \mathcal{X}$, where $\det(K)_{R,C}$ is the minor of $K$ restricted to the rows and columns indexed respectively by $R$ and $C$. We write $\det(K)_{:,C}$ (resp. $\det(K)_{R,:}$) in case $R = \mathcal{X}$ (resp. $C = \mathcal{X}$).
        \end{paragraph}
        
        \begin{paragraph}{Marginal kernel of $\mathcal{D}_\mathcal{M}$}
        Let us now describe the marginal kernel $K$ that we will associate to distribution $\mathcal{D}_\mathcal{M}$. To this end, we consider a \emph{twisted discrete differential} operator $\nabla: \C^\mathcal{V} \rightarrow \C^\mathcal{E}$, mapping complex values defined on the complex planes $\C_v$ (associated to the nodes of $\mathcal{G}$) to copies of the complex planes $\C_e$ associated to the \emph{edges} of $\mathcal{G}$. Its expression relies on a splitting of the connection maps $\psi_e: \C_{s_e} \xrightarrow[]{\psi_{s_e,e}} \C_e \xrightarrow[]{\psi_{e,t_e}} \C_{t_e}$  such that $\psi_e = \psi_{e,t_e} \circ \psi_{s_e, e}$ (this is always possible, \emph{e.g.} $\psi_{e,t_e} = \id_{\C_e,\C_{t_e}}$ and $\psi_{s_e,e}(z) = e^{\iota \theta_e} \cdot z$), and reads~\cite{kenyon_spanning_2011}:
        \begin{equation}
            (\nabla f)(e) = \sqrt{w_e} \psi_{t_e,e}(f(t_e)) - \sqrt{w_e} \psi_{s_e,e}(f(s_e)).
        \end{equation}
        Expliciting the entries of the matrix of $\nabla$, we have:
        \begin{equation}
            \nabla_{e,v} = \begin{cases}
                - \sqrt{w_e} \psi_{s_e,e} & \text{if $v = s_e$} \\
                \sqrt{w_e} \psi_{e,t_e} & \text{if $v = t_e$} \\
                0 & \text{otherwise.}
            \end{cases}
        \end{equation}
        We finally take $\mathcal{X} = \mathcal{V} \cup \mathcal{E}$ and define the marginal kernel $K$:
        \begin{equation}
            K = \nabla_\mathsf{Q} (\mathsf{L}_\theta + \mathsf{Q})^{-1} \nabla_Q^*, \text{ where } \nabla_\mathsf{Q} = \begin{bmatrix}
                \nabla \\ \sqrt{\mathsf{Q}}
            \end{bmatrix}.
            \label{eq:kernel}
        \end{equation}
        \begin{rem}
            For trivial connections, $\nabla$ is the edge-vertex incidence matrix of $\mathcal{G}$, and we have $\mathsf{L} = \nabla^t \nabla$. In a similar manner, one shows that $\mathsf{L}_\theta = \nabla^* \nabla$.
        \end{rem}
        \noindent We will show that $\mathcal{D}_\mathcal{M}$ is a DPP with kernel $K$. First, notice that
        \begin{equation}
            \mathsf{L}_\theta + \mathsf{Q} = \nabla_\mathsf{Q}^* \nabla_\mathsf{Q},
        \end{equation}
        so that $K$ is a projection operator (\emph{i.e.} $K^2 = K$), and all its eigenvalues are in $\{0,1\}$. DPPs associated to such kernels are known as \emph{projection DPPs}, and can be defined without resorting to marginal probabilities, with $\mathbf{P}_{\mathrm{DPP}(K)}(X) = \det(K)_{X,X}$ when $\vert X \vert = \mathrm{rk}(K)$ the rank of $K$ ($\mathrm{rk}(K) = \vert \mathcal{ V } \vert$ in our case), and $\mathbf{P}_{\mathrm{DPP}(K)}(X) = 0$ otherwise~\cite{hough_determinantal_2006}.
        \\
        This will allow to show that the samples of $\mathrm{DPP}(K)$ are MTSFs, distributed according to $\mathcal{D}_\mathcal{M}$ (we generalize here an argument from~\cite{kenyon_spanning_2011}). Let $\phi \subseteq \mathcal{V} \cup \mathcal{E}$ with $\vert \phi \vert = \vert \mathcal{V} \vert$, and remark that:
        \begin{equation}
            \det(K)_{\phi,\phi} = \frac{\det(\nabla_\mathsf{Q}^* \nabla_\mathsf{Q})_{\phi,\phi}}{\det(\mathsf{L}_\theta + \mathsf{Q})}.
        \end{equation}
        The next step then consists in expliciting $\det(\nabla_\mathsf{Q}^* \nabla_\mathsf{Q})_{\phi,\phi} = \det(\nabla_\mathsf{Q})_{\phi,:} \det(\nabla_\mathsf{Q}^*)_{:,\phi}$ depending on $\phi$: we start by inspecting $\det(\nabla_\mathsf{Q}^*)_{:,\phi}$ (similar arguments apply to $\det(\nabla_\mathsf{Q})_{\phi,:} $). In this case, $\phi$ indexes the columns of $\nabla_\mathsf{Q}$, with two columns linearly independent if they do not belong to the same component $c_\phi$. If $c_\phi$ spans $m$ nodes, $\det(\nabla_\mathsf{Q}^*)_{:,c_\phi}$ can only be non-zero if $\vert c_\phi \vert = m$, so that $c_\phi$ must be either a unicycle or a rooted tree. We can then can compute $\det(\nabla_\mathsf{Q}^* \nabla_\mathsf{Q})_{c_\phi,c_\phi}$ from $\det(\nabla_\mathsf{Q}^*)_{:,\phi}$ and $\det(\nabla_\mathsf{Q})_{\phi,:} $ explicitly in those cases.
        \begin{lem}
            \begin{equation}
                \det(\nabla_\mathsf{Q}^* \nabla_\mathsf{Q})_{c_\phi,c_\phi} = \begin{cases}
                    (2 - 2\cos(\theta_C)) \prod_{e \in c_\phi} w_e & \text{if $c_\phi$ is a unicycle with cycle $C$} \\
                    q_r \prod_{e \in c_\phi \cap \mathcal{E}} w_e & \text{if $c_\phi$ is a tree rooted in $r$}
                \end{cases}
        \end{equation}
        \label{lem:calcul}
        \end{lem}
        \begin{proof}
            Suppose first that $c_\phi$ is a unicycle (this is the case treated in~\cite{kenyon_spanning_2011}). We fix an orientation of the edges of $c_\phi$ such that all edges are oriented \emph{towards} the cycle $C$, and that the edges belonging to $C$ form a directed cycle (there are only two such orientations, we choose one arbitrarily). The only non-zero permutations in the determinant
            \begin{equation}
            \det(\Delta_\mathsf{Q}^*)_{:,c_\phi} = \sum_{\sigma \in \mathcal{S}_m} \prod_{1 \leq i \leq m} \mathrm{sgn}(\sigma) (\Delta_\mathsf{Q}^*)_{i,\sigma(i)},
            \end{equation}
            where $\mathrm{sgn}(\sigma)$ denotes the signature of permutation (bijection from nodes to edges) $\sigma$,
            correspond to these orientations, and map vertices $s_e$ to edges $e$ ($\sigma(s_e) = e$). These two bijections, denoted $\sigma_C$ and $\sigma_{C^*}$, differ only on nodes adjacent to edges in $C$, and correspond to the two possible orientations of the cycle. We then obtain:
            \begin{equation}
                \det(\Delta_\mathsf{Q}^*)_{:,c_\phi} = (-1)^m \mathrm{sgn}(\sigma_C) \prod_{e \in c_\phi \setminus C} \sqrt{w_e} \psi_{s_e,e} \left( \prod_{e' \in C} \sqrt{w_{e'}} \psi_{s_{e'},e'} - \prod_{e' \in C} \sqrt{w_{e'}} \psi_{t_{e'},e'} \right).
            \end{equation}
            Performing a similar computation for $\det(\Delta_\mathsf{Q})_{c_\phi,:}$ and multiplying the two resulting expressions then allows to show that (this computation also appears in~\cite{kenyon_spanning_2011}):
            \begin{equation}
            \det(\nabla_\mathsf{Q}^* \nabla_\mathsf{Q})_{c_\phi,c_\phi} = (2 - 2\cos(\theta_C)) \prod_{e \in c_\phi} w_e,
            \end{equation}
            where we used $\psi_e = \psi_{e,t_e} \circ \psi_{s_e, e}$ and $\psi_C + \psi_C^* = 2 \cos(\theta_C)$.
            \\
            Similarly, for $c_\phi$ a rooted tree with root $r \in c_\phi \cap \mathcal{V}$, we have
            \begin{equation}
                \det(\Delta_\mathsf{Q}^*)_{:,c_\phi} = \sqrt{q_r} \times (-1)^m \mathrm{sgn}(\sigma_C) \prod_{e \in (c_\phi \cap \mathcal{e}) \setminus C} \sqrt{w_e} \psi_{s_e,e},
            \end{equation}
            which results in
            \begin{equation}
            \det(\nabla_\mathsf{Q})_{c_\phi,c_\phi} = q_r \prod_{e \in c_\phi \cap \mathcal{E}} w_e.
            \end{equation}
        \end{proof}
        
        \noindent As a corollary, we obtain
        \begin{equation}
            \prb_{\mathrm{DPP}(K)}(\phi) = \frac{\prod_{r \in \phi \cap \mathcal{V}} q_r \prod_{e \in \phi \cap \mathcal{E}} w_e \prod_{C \in \mathcal{C}(\phi)} \left((2 - 2 \cos(\theta_C)\right)}{\det(\mathsf{L}_\theta + \mathsf{Q})},
        \end{equation}
        which is exactly $\prb_{\mathcal{D}_\mathcal{M}}(\phi)$.
        \end{paragraph}
        
        \begin{paragraph}{Unbiased estimator}
        It remains to show that $\tilde{f}(i,\phi,g)$ is an unbiased estimator of the desired quantity $\left( (\mathsf{L}_\theta + \mathsf{Q})^{-1} \mathsf{Q} g \right) (i)$. Our argument is inspired from that of~\cite{pilavci_graph_2021}, and we will require another determinantal tool.
        \begin{prop}[Cauchy-Binet formula]
            For two matrices $A \in M_{m,n}(\C)$ and $B \in M_{n,m}(\C)$ with $m < n$, we have:
            \begin{equation}
                \det(AB) = \sum_{T} \det(A)_{:,T} \det(B)_{T,:},
            \end{equation}
            where $T$ ranges over all subsets of $\{1,...,n\}$ of size $m$.
        \end{prop}
        We begin by rewriting $(\mathsf{L}_\theta + \mathsf{Q})^{-1}_{i,j}$ as an expectation, starting from Cramer's rule
        \begin{equation}
            (\mathsf{L}_\theta + \mathsf{Q})^{-1}_{i,j}  = (-1)^{i+j} \frac{\det(\mathsf{L}_\theta + \mathsf{Q})_{\mathcal{V} \setminus \{j\}, \mathcal{V} \setminus \{i\}}}{\det(\mathsf{L}_\theta + \mathsf{Q})}.
        \end{equation}
        We can rewrite the numerator using the Cauchy-Binet formula:
        \begin{align}
            \det(\mathsf{L}_\theta + \mathsf{Q})_{\mathcal{V} \setminus \{j\}, \mathcal{V} \setminus \{i\}} & = \sum_{\substack{{\phi \subseteq \mathcal{V} \cup \mathcal E} \\ {\vert \phi \vert = \vert \mathcal{V} \vert - 1}}} (-1)^{i+j} \det\left( \begin{bmatrix} (\nabla_\mathsf{Q}^*)_{:,\phi} \ \delta_j \end{bmatrix} \right) \det\left( \begin{bmatrix} (\nabla_\mathsf{Q})_{\phi,:} \\ \delta_i \end{bmatrix} \right),
        \end{align}
        where $\delta_i$ is the $i$-th vector in the usual basis of $\R^\mathcal{V}$. An argument analogous to Lemma~\ref{lem:calcul} then shows that the product of determinants 
        \begin{equation}
            \det\left( \begin{bmatrix} (\nabla_\mathsf{Q}^*)_{:,\phi} \ \delta_j \end{bmatrix} \right) \det\left( \begin{bmatrix} (\nabla_\mathsf{Q})_{\phi,:} \\ \delta_i \end{bmatrix} \right)
        \end{equation}
        can only be non-zero if $\phi$ contains a tree $T_i^j \subseteq \mathcal{E}$ spanning both $i$ and $j$ (with no root), with contribution $\psi_{j \xrightarrow[]{\phi} i} \prod_{e \in T_i^j} w_e$ to the product, and if \emph{all} the other components are rooted trees or unicycles, with associated contributions described in Lemma~\ref{lem:calcul}. 
        We thus obtain
        \begin{align}
            (\mathsf{L}_\theta + \mathsf{Q})^{-1}_{i,j} & = \frac{1}{q_j} \sum_{\phi \in \mathcal{M}(\mathcal{G})} \prb_{\phi \sim \mathcal{D}_\mathcal{M}}(\phi) \mathbf{1}_{c_\phi(j)}(i) \psi_{j \xrightarrow[]{\phi} i} \\
            & = \frac{1}{q_j} \mathbf{E}_{\phi \sim \mathcal{D}_\mathcal{M}} \left(\mathbf{1}_{c_\phi(j)}(i) \psi_{j \xrightarrow[]{\phi} i} \right),
        \end{align}
        where $\mathbf{1}_{c_\phi(j)}(i)$ is the indicator that $i$ belongs to the set of nodes $c_\phi(j)$ spanned by the tree rooted in $r$. Finally, we have
        \begin{align}
            \left( (\mathsf{L}_\theta + \mathsf{Q})^{-1} \mathsf{Q} g \right) (i) & = \langle \delta_i, (\mathsf{L}_\theta + \mathsf{Q})^{-1} \mathsf{Q} g \rangle \\
            & = \sum_{j \in \mathcal{V}} q_j (\mathsf{L}_\theta + \mathsf{Q})^{-1}_{i,j} g(j) \\
            & = \sum_{j \in \mathcal{V}} \mathbf{E}_{\phi \sim \mathcal{D}_\mathcal{M}}\left( \psi_{j \xrightarrow[]{\phi} i} \left( g(j) \right) \mathbf{1}_{c_\phi(j)}(i) \right) \\
            & = \mathbf{E}_{\phi \sim \mathcal{D}_\mathcal{M}}(\tilde{f}(i,\phi,g)).
        \end{align}
        \end{paragraph}
        
        \begin{paragraph}{Variance computation}
        \textcolor{col}{The expression of the variance relies on the observation that $\widetilde{f}(\phi,g)$ is \emph{linear in $g$}, so that we can write $\widetilde{\mathsf{S}}_\phi g = \widetilde{f}(\phi,g)$ for some linear operator/matrix $\widetilde{\mathsf{S}}_\phi$, and otherwise consists in a simple computation. Note that we have $\mathbf{E}_{\phi \sim \mathcal{D}_\mathcal{M}}\left(\widetilde{\mathsf{S}}_\phi g \right) = f_*$ and
    \begin{equation}
    \mathbf{E}_{\phi \sim \mathcal{D}_\mathcal{M}}\left(\widetilde{\mathsf{S}}_\phi\right) = (\mathsf{L}_\theta + \mathsf{Q})^{-1} \mathsf{Q}.
    \end{equation}
    \\
    First, notice that $\widetilde{\mathsf{S}}_\phi^* \mathsf{Q} \widetilde{\mathsf{S}}_\phi = \mathsf{Q}$. We then have
        \begin{align}
            \mathbf{E}_{\phi \sim \mathcal{D}_\mathcal{M}} \left( \Vert \widetilde{S}(\phi,g) - f_* \Vert_\mathsf{Q}^2 \right) & = \sum_{i \in \mathcal{V}} q_i \left(  \mathbf{E}\left(\left\vert \left(\widetilde{S}_\phi g \right)(i) \right\vert^2\right) - \left\vert \mathbf{E}\left(\left(\widetilde{S}_\phi g \right)(i) \right)^2 \right\vert \right) \label{eq:variance1} \\
            & = \left( \sum_{i \in \mathcal{V}} \langle g, \mathbf{E}_{\phi \sim \mathcal{D}_\mathcal{M}} \left( \widetilde{\mathsf{S}}_\phi^* (q_i \delta_i \delta_i^*) \widetilde{\mathsf{S}}_\phi \right) g  \rangle  \right) - \Vert f_* \Vert^2_\mathsf{Q} \label{eq:variance2} \\
            & = \Vert g \Vert^2_\mathsf{Q} - \Vert f_* \Vert^2_\mathsf{Q} \label{eq:variance3},
        \end{align}
        where $\delta_i$ is usual basis vector of $\C^\mathcal{V}$ with non-zero entry $\delta_i(i) = 1$, Equation~\eqref{eq:variance2} follows from the linearity of the expectation, and Equation~\eqref{eq:variance3} from the previous observation that $\widetilde{\mathsf{S}}_\phi^* \mathsf{Q} \widetilde{\mathsf{S}}_\phi = \mathsf{Q}$.}
        \rightline{\color{header1}\qedsymbol{}}
    \end{paragraph}

\begin{rem}
\label{rem:trace_estimation}
\textcolor{col}{We claimed in Section~\ref{subsect:runtime_precision_smoothing} that the quantity $s(q)$ in Equation~\eqref{eq:sq} is given, in expectation, by the number of roots of MTSFs sampled according to $\mathcal{D}_\mathcal{M}$. This fact hinges on two well-known properties of determinantal point processes (see, \emph{e.g.},~\cite{kulesza_determinantal_2012}):}
\begin{enumerate}
    \item \textcolor{col}{if $X$ follows a DPP with over $\mathcal{X}$ with kernel $K$, and if $S \subseteq \mathcal{X}$, then $X \cap S$ also follows a DPP, with kernel $K_S = K_{S,S}$;}
    \item \textcolor{col}{the cardinality $\vert X \vert$ of a DPP-sample is in general a random quantity, and its expectation is given by
    \begin{equation}
        \mathbf{E}_{X \sim \mathrm{DPP}(K)} ( \vert X \vert ) = \mathrm{tr}(K).
        \label{eq:cardinal_dpp}
    \end{equation}}
\end{enumerate}
\textcolor{col}{In particular, the kernel $K_\mathcal{V}$ obtained by restricting the kernel $K$ of $\mathcal{D}_\mathcal{M}$ (Equation~\eqref{eq:kernel}) describes a distribution over the nodes $\mathcal{V}$, and reads
\begin{equation}
    K_\mathcal{V} = q(\mathsf{L}_\theta + q \mathsf{I})^{-1}
\end{equation}
or, more generally for heterogeneous $q_i$'s:
\begin{equation}
    K_\mathcal{V} = \sqrt{\mathsf{Q}} (\mathsf{L}_\theta + q \mathsf{I})^{-1} \sqrt{\mathsf{Q}}.
\end{equation}
The claim follows from Equation~\eqref{eq:cardinal_dpp}:
\begin{equation}
    \mathbf{E}_{\phi \sim \mathcal{V}}(\vert \phi \cap \mathcal{V} \vert) = \mathrm{tr}\left(\sqrt{\mathsf{Q}} (\mathsf{L}_\theta + q \mathsf{I})^{-1} \sqrt{\mathsf{Q}}\right),
\end{equation}
which translates to an efficient estimator for the effective degrees of freedom of the connection-aware smoothing problem.}
\end{rem}

    \section{Complexity Analysis}
    \label{supmat:complexity}
    
    We will describe the complexity of different implementations of the MTSF-sampling algorithm.
    We recall the MTSF-sampling procedure in Algorithm~\ref{alg:sampling_mtsf_bis}. We consider here the generalization to weighted graphs with heterogeneous $q_i$'s, and use a generalized $\mathrm{random\_successor}(u)$ that outputs either $\Gamma$ with probability $\frac{q_u}{d_u + q_u}$, or some node $v$ with probability  $\frac{w_{u,v}}{d_u + q_u}$. This generalized algorithms provably samples MTSFs from the distribution $\mathcal{D}_\mathcal{M}$, with the argument from~\cite{fanuel_sparsification_2022} still applying.
    
    \begin{algorithm}[ht]
        \caption{MTSF sampling algorithm~\cite{fanuel_sparsification_2022}.}
        \begin{algorithmic}[1]
            \State $\phi \leftarrow \emptyset$
            \While{$\phi$ not spanning}
                \State Let $i \in \mathcal{V}$ be the first node in the queue not spanned by $\phi$
                \State $u \leftarrow i$ \Comment{$u$ is the current node of the random walk}
                \State $p \leftarrow \epsilon$
                \Comment{$\epsilon$ the empty path}
                \While{($u \neq \Gamma$) and ($p$ does not intersect $\phi$ or contain a cycle)}
                    \State $u' \leftarrow \mathrm{random\_successor}(u)$ \Comment{Move to next node}
                    \If{$u' \neq \Gamma$}
                    \State $e \leftarrow (u,u')$, $p \leftarrow p e$ \Comment{Add $e$ to the path $p$}
                    \EndIf
                    \If{${p}$ contains a cycle $C$}
                        \State Remove this cycle from $p$ with probability $\cos(\theta_C)$
                    \EndIf
                    \State $u_\mathrm{old} \leftarrow u$, $u \leftarrow u'$ \Comment{$u_\mathrm{old}$ the previous node}
                \EndWhile
                \If{$u = \Gamma$}
                    \State $\phi \leftarrow \phi \cup p \cup u_\mathrm{old}$ \Comment{Add the sampled path $p$ and the root $u_\mathrm{old}$ to $\phi$}
                \Else
                    \State $\phi \leftarrow \phi \cup p$
                \EndIf
            \EndWhile
            \State \textbf{Output} $\phi$
        \end{algorithmic}        \label{alg:sampling_mtsf_bis}
        \end{algorithm}
    
    \subsection{Number of Steps}
    
    Let us re-state the upper bound on the expected number of steps of the random walk in the sampling Algorithm~\ref{alg:sampling_mtsf}. Denote by $T_\phi$ the number of random neighbors sampled to build $\phi \in \mathcal{M}(\mathcal{G})$ during an execution of Algorithm~\ref{alg:sampling_mtsf_bis}. Then, we have:
    \begin{prop}
    \begin{equation}
        \mathbf{E}_{\phi \sim \mathcal{D}_\mathcal{M}}(T_\phi) \leq \mathrm{tr}\left( (\mathsf{L} + \mathsf{Q})^{-1} (\mathsf{D} + \mathsf{Q})  \right),
        \label{eq:bound_mtsf_bis}
    \end{equation}
    \label{prop:sampling_steps}
    \end{prop}
    
    \noindent where we abuse notation and write $\phi \sim \mathcal{D}_\mathcal{M}$ for $\phi$ sampled using Algorithm~\ref{alg:sampling_mtsf_bis}. This bound is easily obtained by considering the special case of the trivial connection with $\psi_e = \mathrm{id}_{\C_{s_e},\C_{t_e}}$ for all $e \in \overrightarrow{\mathcal{E}}$, in which case Algorithm~\ref{alg:sampling_mtsf_bis} samples spanning forests of $\mathcal{G}$ from the distribution mentioned in Remark~\ref{rem:yigit}, and the expected number of steps is known to be~\cite{pilavci_graph_2021}:
    \begin{equation}
        \mathrm{tr}\left( (\mathsf{L} + \mathsf{Q})^{-1} (\mathsf{D} + \mathsf{Q})  \right).
    \end{equation}
    This is also the slowest scenario: in this case, $\theta_C = 0$ for all cycles, and $p$ can never contain a cycle in step $\mathrm{6}$ of Algorithm~\ref{alg:sampling_mtsf_bis}. Hence,  the only way to exit the $\mathrm{\mathbf{while}}$ loop is to reach $\Gamma$ or to intersect $\phi$, and we obtain the bound of Proposition~\ref{prop:sampling_steps} as a consequence.
    
    Let us briefly comment on the value of the trace in Equation~\eqref{eq:bound_mtsf_bis}. We have:
    \begin{equation}
        \mathrm{tr} \left( \mathsf{L} + \mathsf{Q})^{-1} (\mathsf{D} + \mathsf{Q}) \right) = \sum_{i \in \mathcal{V}} (\mathsf{L} + \mathsf{Q})^{-1}_{i,i} (q_i + d_i),
    \end{equation}
    where it is known that $(\mathsf{L} + \mathsf{Q})^{-1}_{i,i}$ equals the probability that $i \in \mathsf{f} \cap \mathcal{V}$ for $\mathsf{f} \sim \mathcal{D}_\mathcal{F}$, and hence lies in $[0,1]$. It follows that $\mathbf{E}_{\phi \sim \mathcal{D}_\mathcal{M}}(T_\phi)$ is in $\mathcal{O}\left( \frac{\vert \mathcal{E} \vert}{q_\mathrm{min}} \right)$.
    
    \begin{rem}
        Algorithm~\ref{alg:sampling_mtsf_bis} can be used in conjunction with Equation~\eqref{eq:heterogeneous_fermionic} to estimate $q(\widetilde{\mathsf{L}_\theta} + q \mathsf{I})^{-1} g$ (see Remark~\ref{rem:synchro}). In this context, we take $\mathsf{Q} = q \mathsf{D}$, and obtain the bound $\mathbf{E}_{\phi \sim \mathcal{D}_\mathcal{M}}(T_\phi) \in \mathcal{O}\left( \left( 1 + \frac{1}{q} \right) \vert \mathcal{V} \vert \right)$.
    \end{rem}
    
    \subsection{Implementation-specific Complexity}
    
    The time-complexity of Algorithm~\ref{alg:sampling_mtsf_bis} does not only depend on the number of steps taken by the random walk, but also on the cost of the cycle-detection (step 10). We propose two counter-based solutions.
    
    \vspace{0.2cm}
    
    \begin{paragraph}{One-counter detection}
    The first strategy consists in dynamically assigning a numerical value $\mathrm{c}_v \in \N$ to each node $v \in \mathcal{V}$ encountered during the random walk (initialized to $\mathrm{c}_v = 0$), based on a global counter $\mathrm{c} \in \N$. This counter is incremented each time a new random walk begins (instruction at line 3), and we set $\mathrm{c}_{u'} = \mathrm{c}$ whenever the random walk reaches $u'$ on the instruction at line 7 (if it not already spanned by MTSF). If a cycle $C$ is created in node $u'$ and discarded, we reset $\mathrm{c}_v = 0$ for all nodes $v$ spanned by $C$ \emph{except} $u'$. A cycle can then be detected in $\mathcal{O}(1)$ time by checking the value of $\mathrm{c}_{u'}$ at each step of the random walk: there is a cycle if $\mathrm{c}_{u'} = \mathrm{c}$, and no cycle if $\mathrm{c}_{u'} < \mathrm{c}$.
    \\
    In a run of Algorithm~\ref{alg:sampling_mtsf_bis}, resetting the values $\mathrm{c}_v$ when cycles are discarded requires going through at most $T_\phi$ nodes, which results in $\mathcal{O}\left( \frac{\vert \mathcal{E} \vert}{q_\mathrm{min}} \right)$ overall expected \emph{time complexity}.
    \end{paragraph}
    
    \vspace{0.2cm}
    
    In practice, we found it faster to use an implementation of the following strategy.
    
    \vspace{0.2cm}
    
    \begin{paragraph}{Multiple-counters detection}
    The $\mathrm{c}_v$'s need to be reset in the one-counter cycle-detection strategy because, whenever a cycle is discarded, nodes spanned by this cycle should no longer be remembered as spanned by the path $p$. This expensive resetting step can be bypassed by considering \emph{multiple} counters, in addition to the $\mathrm{c}$'s (indexed by a global $\mathrm{c}$). We will keep track of couples of values $(\mathrm{id}_v,\mathrm{val}_v) \in \N^2$ for each $v$, initialized to $(0,0)$) and with updates based on two global counters $\mathrm{id}, \mathrm{val} \in \N^2$, and of values $\mathrm{cap}(\mathrm{id})$ associated to each $\mathrm{id}$. $\mathrm{id}_v$ should be thought of as the ID of a counter at $v$ (amongst multiple others), and $\mathrm{val}_v$ as its value (we only need to store the value associated to the largest $\mathrm{id}_v$). These counters are reset to $0$ whenever $\mathrm{c}$ is incremented (and a new random walk is initiated), and are otherwise updated by applying the following rules for all nodes $u'$ reached by the random walk (if not already spanned by the MTSF).
    \begin{itemize}
        \item If no cycle is created and discarded at $u'$, set $\mathrm{id}_{u'} = \mathrm{id}$, $\mathrm{val}_{u'} = \mathrm{val}$, and increment $\mathrm{val}$.
        \item If a cycle is created at node $u'$ and discarded, store the value $\mathrm{cap}(\mathrm{id}) = \mathrm{val}$ (the maximum value for the $\mathrm{id}^{th}$ counter), set $\mathrm{cap}(k) = 0$ for all $k > \mathrm{id}_{u'}$, $\mathrm{val}_{u'} = \mathrm{val}$, and increment $\mathrm{id}$.
    \end{itemize}
    \noindent Using these counters, a cycle is detected at $u'$ if $\mathrm{val}_{u'} \leq \mathrm{cap}(\mathrm{id}_{u'})$.
    \\
    The number $\mathrm{id}$ of counters used is unbounded, and we did not derive an expected time complexity for this multiple-counters strategy, but consistently obtained better performance in our measurements when using it. 
    \end{paragraph}
    
    \begin{rem}
    The $\mathcal{O}\left( \frac{\vert \mathcal{E} \vert}{q_\mathrm{min}} \right)$ complexity can seem daunting for very small values of $q_i$'s, but one should remember that this is only an upper bound, not accounting for the presence of unicycles in the sample. As the $q_i \rightarrow 0$, $\mathrm{DPP}(K)$ approaches a distribution over spanning forests of unicycles, with a Wilson-like sampling algorithm described in~\cite{kassel_random_2017}. Our implementation strategy can also be applied in this case and translates to the bounds on the number of steps of the random walks they derive (see also~\cite{fanuel_sparsification_2022}).
    \end{rem}
    
    Finally, we point out another implementation detail.
    
    \begin{rem}
            \label{rem:sampling_connection}
            The propagation maps used  in Algorithm~\ref{alg:mtsf_propagation} can be computed during the sampling process and, one actually only needs knowledge of: the root of each tree, whether or not a node $i$ belongs to a rooted tree and, if it does, the propagation map from the root $r_\phi(i)$ to $i$. The actual MTSF $\phi$ is never used in the estimation.
    \end{rem}

    \section{Variance Reduction}
    \label{supmat:variance}

    We show that the estimators in Equations~\eqref{eq:rao} and~\eqref{eq:gradient_step} are unbiased (Proposition~\ref{prop:variance_reduction}).
    
    \begin{paragraph}{Rao-Blackwell estimator}
        Let us first re-define $\overline{f}$ for heterogeneous values of $q$. The difference resides in the aggregation function $h_\phi$:
        \begin{equation}
            h_\phi(r,g) = \frac{\sum_{j \in c_\phi(r)} \psi_{j \xrightarrow[]{\phi} r} g(j) } {\sum_{r' \in c_\phi(r)} q_{r'}},
        \end{equation}
        where we recall that $c_\phi(r)$ denotes the set of nodes spanned by the tree containing $r$. We still take $\overline{f}(i,\phi,g) = \psi_{r_\phi(i) \xrightarrow[]{\psi} i} (h_\phi(r_\phi(i),g))$, and show that:
        
        \begin{prop}
        $\mathbf{E}_{\phi \sim \mathcal{D}_\mathcal{M}} (\overline{f}(i,\phi,g)) = f_*(i)$.
        \label{prop:rao}
        \end{prop}
        \begin{proof}
            We will express $\overline{f}$ as a conditional expectation, the result will follow from the law of total expectation. Here, we choose a connected subset of edges $\pi \subseteq \mathcal{E}$, and condition on $\phi \cap \mathcal{E}$ containing $\pi$ as one of its (maximal) components, which we denote by $\pi \sqsubseteq_\mathcal{E} \phi$. For $i$ belonging to the connected component spanned by $\pi$, we have:
            \begin{align}
                \mathbf{E}_{\phi \sim \mathcal{D}_\mathcal{M}}(\tilde{f}(i,\phi,g) \ \vert \  \pi \sqsubseteq_\mathcal{E} \phi) & = \sum_{\phi \in \mathcal{M}(\mathcal{G})} \prb_{\mathcal{D}_\mathcal{M}}(\phi \ \vert \ \pi \sqsubseteq_\mathcal{E} \phi) \ \tilde{f}(i,\phi,g) \\
                & = \sum_{\substack{ {\phi \in \mathcal{M}(\mathcal{G})} \\ {\pi \sqsubseteq_\mathcal{E} \phi}} } \frac{q_r \mathbf{1}_{\phi \cap c_\phi(i)}(r)}{\sum_{r' \in c_\phi(i)} q_{r'}} \ \tilde{f}(i,\phi,g) \\
                & = \overline{f}(i,\phi,g).
            \end{align}
        \end{proof}
        
        \textcolor{col}{The variance of this improved estimator is derived similarly to that of $\tilde{f}$. In particular, $\overline{f}$ is linear in $g$ and, denoting by $\overline{\mathsf{S}}_\phi$ the operator such that $\overline{\mathsf{S}}_\phi g = \overline{f}(\phi,g)$, which is such that $\overline{\mathsf{S}}_\phi^* \mathsf{Q} \overline{\mathsf{S}}_\phi = \mathsf{Q} \overline{\mathsf{S}}_\phi$ and $\mathbf{E}_{\phi \sim \mathcal{D}_\mathcal{M}}(\overline{\mathsf{S}}_\phi) = (\mathsf{L}_\theta + \mathsf{Q})^{-1} \mathsf{Q}$, we obtain from the same argument that. 
        \begin{prop}
        $\mathbf{E}_{\phi \sim \mathcal{D}_\mathcal{M}}\left( \Vert \overline{f}(\phi,g) - f_* \Vert^2_\mathsf{Q} \right) = \langle g, \mathsf{Q} f_* \rangle - \Vert f_* \Vert^2$.
        \end{prop}}
        
        \vspace{0.2cm}
        
        \textcolor{col}{Finally, we show the following concentration bound.
        \begin{prop}
            Let $\eps,\delta \in (0,1)$, and consider sampling $m$ MTSFs $\{\phi_k\}_{k=1}^m$. Then, if
        \begin{equation}
            m \geq \frac{6}{\eps^2} \log\left( \frac{\vert \mathcal{V} \vert}{\delta} \right),
            \label{eq:m_condition}
        \end{equation}
        it holds for any signal $g \in \mathcal{C}^\mathcal{V}$ that
        \begin{equation}
            \prb\left( \left\Vert \left( \frac{1}{m} \sum_{k = 1}^m \overline{f}(\phi_k,g) \right) - f_* \right\Vert_2 \leq \eps \Vert g \Vert_2 \right) \geq 1 - \delta.
            \label{eq:concentration}
        \end{equation}
        \label{prop:detail_tropp_supmat}
        \end{prop}}

\color{col}{\noindent Let us begin by stating the Matrix-Bernstein concentration bound from from~\cite{tropp2012user}.

\begin{thm}[1.4 from~\cite{tropp2012user}]
Consider a set of $m$ \emph{self-adjoint} matrices $\{\mathsf{X}_k\}_{k=1}^m$ in $M_{n}(\C)$, sampled independently. Assume that, almost surely and for each such matrix, both the conditions
\begin{equation}
    \mathbf{E}(\mathsf{X}_k) = 0
\end{equation}
\begin{equation}
    \Vert \mathsf{X}_k \Vert \leq R
\end{equation}
are satisfied, where $\Vert . \Vert$ denotes the operator norm\footnote{Associated to the usual $l_2$-norm on $\C^n$.}.
\\
Then, for all $\eps \geq 0$,
\begin{equation}
\prb\left( \left\Vert \sum_{k=1}^m \mathsf{X}_k \right\Vert \geq \eps \right) \leq n e^{\frac{\frac{-\eps^2}{2}}{\sigma^2 + R \frac{\eps}{3}}},
\end{equation}
where we define
\begin{equation}
    \sigma^2 = \left\Vert \sum_{k=1}^m \mathbf{E}\left(\mathsf{X}_k^2\right) \right\Vert.
\end{equation}
\end{thm}

\vspace{0.2cm}

\noindent Note that this bound is quite powerful and, in particular, that the matrices $\mathsf{X}_k$ need not be sampled from the same distribution. Our setting is quite simple in comparison, as all the matrices we consider \emph{are} sampled from the same distribution.

\vspace{0.2cm}

\noindent Specifically, we consider the matrices $\mathsf{X}_k = \frac{1}{m} (\overline{\mathsf{S}}_{\phi_k} - \mathsf{K})$ associated to a set of $m$ multi-type spanning forests $\{\phi_k\}_{k=1}^m$, with $\mathsf{K} = (\mathsf{L}_\theta + \mathsf{Q})^{-1} \mathsf{Q}$ denoting the ``matrix'' solution to the connection-aware Tikhonov smoothing problem. We assume that all $m$ multi-type spanning forests are sampled according to $\mathcal{D}_\mathcal{M}$, and forget the subscripts in the following. Note that we have
\begin{equation}
    \mathbf{E} (\mathsf{X}_k ) = 0
\end{equation}
for all $k$.

\vspace{0.2cm}

\noindent In order to determine $R$, we bound the spectrum of each $\mathsf{X}_k$:
\begin{equation}
    \Vert \mathsf{X}_k \Vert \leq \frac{1}{m} \left( \Vert \overline{\mathsf{S}}_{\phi_k} \Vert + \Vert \mathsf{K} \Vert \right).
\end{equation}
Since we have both $\Vert \mathsf{K} \Vert = \frac{q}{q + \lambda_1}$ (where $\lambda_1$ is the smallest eigenvalue of $\mathsf{L}_\theta$) and $\Vert \mathsf{\overline{S}}_{\phi_k} \Vert \leq 1$ (by, \emph{e.g.}, Gershgorin's circle theorem), we can take
\begin{equation}
    R = \frac{1}{m} \left( 1 + \frac{q}{q + \lambda_1} \right).
\end{equation}
In order to simplify our bounds later-on, we can notice that
\begin{equation}
    R \leq \frac{3}{m}.
\end{equation}

\vspace{0.2cm}

\noindent As for the value of $\sigma^2$, we obtain from the triangle inequality that
\begin{align}
    \sigma^2 & = \frac{1}{m} \Vert \mathsf{K} - \mathsf{K}^2 \Vert \\
    & \leq \frac{1}{m} \left(\frac{q}{q + \lambda_1} + \frac{q^2}{(q + \lambda_n)^2} \right).
\end{align}
This expression can be made more interpretable by lower-bounding $\lambda_n$, and we find that:
\begin{equation}
    \sigma^2 \leq \frac{1}{m} \left(  \frac{q}{q + \lambda_1} + \frac{q^2}{(1 + q + \max_{v \in \mathcal{V}} d_v)^2} \right),
\end{equation} 
which is proved by leveraging the bounds on $\lambda_n$ in Section~\ref{sect:bound_lambda_n} below.

\vspace{0.2cm}

\noindent It follows from our previous bounds on $R$ and $\sigma^2$ that 
\begin{equation}
    \left( \sigma^2 + R \frac{\eps}{3} \right) \leq \frac{1}{m} \left( \eps + \frac{q}{q + \lambda_1} + \frac{q^2}{(1 + q + \max_{v \in \mathcal{V}} d_v)^2} \right).
\end{equation}

\vspace{0.2cm}

\noindent Finally, $\Vert \sum_{k=1}^m \mathsf{X}_k \Vert \leq \eps$ implies by definition that
\begin{equation}
    \left\Vert \sum_{k=1}^m \mathsf{X}_k g \right\Vert_2 \leq \eps \Vert g \Vert_2
\end{equation}
for any signal $g \in \C^\mathcal{V}$, so that putting all the previous results together results in:
\begin{equation}
    \prb\left( \left\Vert \left( \frac{1}{m} \sum_{k = 1}^m \overline{\mathsf{S}}_{\phi_k} \right)g - f_* \right\Vert_2 \leq \eps \Vert g \Vert_2 \right) \geq 1 - \delta
\end{equation}
whenever we have
\begin{equation}
    m \geq 2 \frac{\eps + \frac{q}{q + \lambda_1} + \frac{q^2}{(1 + q + \max_{v \in \mathcal{V}} d_v)^2}}{\eps^2} \log\left( \frac{\vert \mathcal{V} \vert}{\delta} \right).
\end{equation}

\noindent The bound from Proposition~\ref{prop:detail_tropp_supmat} follows when taking $\eps < 1$.}
\rightline{\color{header1}\qedsymbol{}}
\end{paragraph}

\vspace{0.2cm}
    
    \begin{paragraph}{Gradient step as control variates}
    We show that $\mathbf{E}_{\phi \sim \mathcal{D}_\mathcal{M}}(\hat{f}(\phi,g)) = f_*$. First recall that $\mathbf{E}_{\mathcal{D}_\mathcal{M}}(\overline{f})~=~f_*$ (and $f_* = q (\mathsf{L}_\theta + q \mathsf{I})^{-1} g$) so that, by linearity of the expectation:
    \begin{equation}
        \mathbf{E}_{\phi \sim \mathcal{D}_\mathcal{M}} (\hat{f}(\phi,g)) = f_* - \alpha \mathsf{P} (g - g) = f_*
    \end{equation}
    where, for heterogeneous $q_i$'s, we set $\mathsf{P} = (\mathsf{Q}^{-1} \mathsf{D} + \mathsf{I})^{-1}$.
    
    This shows that $\hat{f}$ is unbiased. Let us now discuss the choice of $\alpha$. We plot in Figure~\ref{fig:alpha} the mean errors $\Vert \overline{f} - f_* \Vert$ and $\Vert \hat{f} - f_* \Vert$ over $5$ trials, for different values values of $\alpha$. We always take $m = 20$ MTSFs, and set $q = 10^{-2} \times \overline{d}$. We use the random graph models from Section~\ref{sect:numerical}.
    
    \begin{figure}[ht]
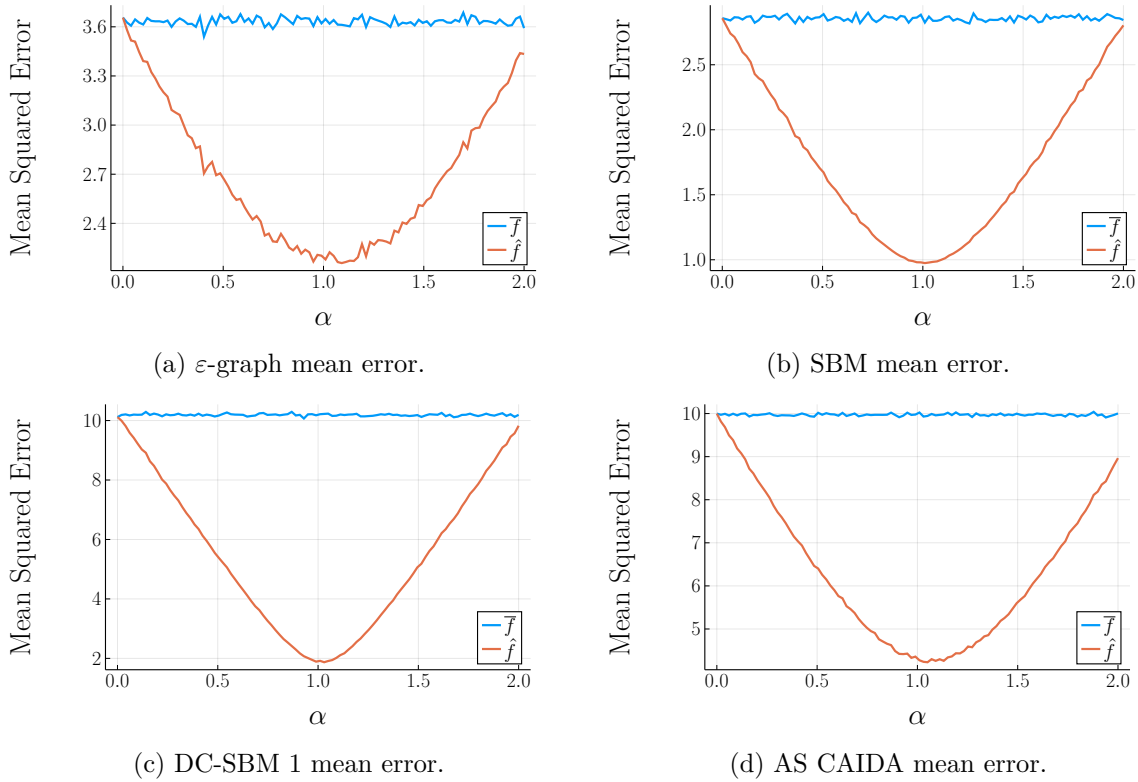

        \centering
        
        \begin{subfigure}{0.49\textwidth}
            \scalebox{0.42}{\input{eps_var.tikz}}
            \caption{$\eps$-graph mean error.}
        \end{subfigure}
        \hfill
        \begin{subfigure}{0.49\textwidth}
            \scalebox{0.42}{\input{sbm_var.tikz}}
            \caption{SBM mean error.}
        \end{subfigure}
        
        \begin{subfigure}{0.49\textwidth}
            \scalebox{0.42}{\input{dc-sbm_mixt_var.tikz}}
            \caption{DC-SBM $1$ mean error.}
        \end{subfigure}
        \hfill
        \begin{subfigure}{0.49\textwidth}
            \scalebox{0.42}{\input{as_caida_var.tikz}}
            \caption{AS CAIDA mean error.}
        \end{subfigure}

        \caption{Mean error as a function of $\alpha$.}
        \label{fig:alpha}
    \end{figure}
    
    The extent of the error reduction depends on the graph, but the optimal choice of $\alpha$ seems to always be close to $1$. This behavior can be observed across order-of-magnitude variations of $q$ (not shown), and setting $\alpha = 1$ resulted in meaningful variance reduction for our other experiments (\emph{e.g.} Section~\ref{sect:numerical}). This contrasts with the (similar) variance reduction technique proposed in~\cite{pilavci2022variance}, developed for connection-free graphs, for which a good choice of step-size was less straightforward.
    \end{paragraph}
    
    \section{Bounds on $\lambda_n$}
    \label{sect:bound_lambda_n}
    
    \textcolor{col}{We show that the largest eigenvalue $\lambda_n$ of $\mathsf{L}_\theta$ is such that
    \begin{equation}
        1 + \max_{v \in \mathcal{V}} d_v \leq \lambda_n \leq 2 \max_{v \in \mathcal{V} d_v},
        \label{eq:bound_lambda_n}
    \end{equation}
    which is nothing but a straightforward generalization of the same bound for the largest eigenvalue of the combinatorial Laplacian.}
    
    \vspace{0.2cm}
    
    \textcolor{col}{The right-hand-side inequality is immediately obtained from Gershgorin circle's theorem. To obtain the left-hand-side inequality, let $v_\mathrm{max} \in \mathcal{V}$ denote one vertex of $\mathcal{G}$ with maximum degree $d_{v_\mathrm{max}} = \max_{v \in \mathcal{V}} d_v$, and consider the signal $f \in \C^\mathcal{V}$ defined by
\begin{equation}
    f(i) = \begin{cases}
        d_{v_\mathrm{max}} \ \text{if $i = v_\mathrm{max}$} \\
        - e^{\iota \theta_{v,i}} \ \text{if $i$ is adjacent to $v_\mathrm{max}$} \\
        0 \ \text{otherwise}.
    \end{cases}
\end{equation}
Clearly, $\frac{\langle f, \mathsf{L}_\theta f \rangle}{\langle f, f \rangle} = 1 + \max_{v \in \mathcal{V}} d_v$, from which the lower-bound on $\lambda_\mathrm{n}$ ensues.}

    \section{Supporting Arguments for Tikhonov Smoothing}
    \label{supmat:tikhonov_numerical}
    
    The goal of this Section is to argue that, when the connection is sufficiently consistent, eigenvectors of $\mathsf{L}$ and $\mathsf{L}_\theta$ can be used similarly in graph-signal-processing applications, such as the smoothing experiment in Section~\ref{sect:numerical}.
    
    \vspace{0.2cm}
    
    \begin{paragraph}{Coherent connection}
        Let us first consider a perfectly coherent connection, such that $\theta_{i,j} = \omega_j - \omega_i$ for all edges $\{i,j\}$. Consider the eigendecomposition $\mathsf{L} = \mathsf{U} \mathsf{\Lambda} \mathsf{U}^*$ of $\mathsf{L}$, and denote respectively by $x$ the vector with entries $x_i = e^{\iota \omega_i}$, and $\mathsf{D}_x$ the diagonal matrix with $(\mathsf{D}_x)_{i,i} = x_i$. In this situation, we have
        \begin{equation}
            \mathsf{L}_\theta = \mathsf{D}_x \mathsf{L} \mathsf{D}_x^* = (\mathsf{D}_x \mathsf{U}) \mathsf{\Lambda} (\mathsf{D}_x \mathsf{U})^*,
            \label{eq:eigenvectors}
        \end{equation}
        and the eigenvectors $v_i$ of $\mathsf{L}_\theta$ are given by $\mathsf{D}_x^* u_i$, with $u_i$ the eigenvectors of $\mathsf{L}$.
    \end{paragraph}
    
    \vspace{0.2cm}
    
    \begin{paragraph}{Noisy connection}
        Equation~\eqref{eq:eigenvectors} no longer holds for incoherent connections. We provide illustrations of the resulting eigenvectors of $\mathsf{L}_\theta$ for connections of the form:
        \begin{equation}
            \theta_{i,j} = \omega_j - \omega_i + \eta \eps_{i,j},
        \end{equation}
        for different values of $\eta$, and $\eps_{i,j}$ uniformly distributed in $[-1,1]$ (like in Section~\ref{sect:numerical}). Specifically, we represent $u_i$, $v_i$, and $\mathsf{D}_{x}^* v_i$ (which should approximate $u_i$). We also measure the difference between eigenvectors $u_i$ and $v_i$ as a function of $\eta$ using the normalized error
        \begin{equation}
            \min_{r \in U(\C)} \ \frac{\Vert u_i - r \left(D_x^* v_i\right) \Vert_2}{n},
        \end{equation} and use the following graphs:
        \begin{itemize}
            \item An $\eps$-graph built from $n = 100$ points uniformly sampled in $[0,1]^3$, with $\eps = 0.3$ (Figure~\ref{fig:gsp1}).
            \item A SBM with two communities on size $50$ ($n = 100$), with $c_{1,1} = c_{2,2} = 19$ and $c_{1,2} = 1$ (Figure~\ref{fig:gsp2}).
        \end{itemize}
        We only use a single realization of these graphs, equipped with a single realization of the connection. The normalized errors are computed across a linear range of values of $\eta$ in $[0,1]$. For each graph, we plot the real and imaginary parts of the entries of the corresponding eigenvectors, for $\eta_1 \simeq 0.11$ and $\eta_2 \simeq 0.75$. The corresponding errors are also highlighted, along with the value $\eta_0 = \frac{\pi}{2n}$ (for which we know that the connection always satisfies the weak-inconsistency Condition~\ref{cond:weak_inconsistency}).
        Recall that the eigenvectors of $\mathsf{L}_\theta$ are only defined up to a global rotation, and this is apparent in our illustrations.

        \begin{figure}[]
        \centering
        
        \begin{subfigure}{\textwidth}
        \begin{subfigure}{0.49\textwidth}
            \scalebox{0.42}{\input{eps_j12_vp2.tikz}}
        \end{subfigure}
        \hfill
        \begin{subfigure}{0.49\textwidth}
            \scalebox{0.42}{\input{eps_j12_vp3.tikz}}
        \end{subfigure}
        \caption{$\eta_1 \simeq 0.11$}
        \end{subfigure}
        
        \begin{subfigure}{\textwidth}
        \begin{subfigure}{0.49\textwidth}
            \scalebox{0.42}{\input{eps_j75_vp2.tikz}}
        \end{subfigure}
        \hfill
        \begin{subfigure}{0.49\textwidth}
            \scalebox{0.42}{\input{eps_j75_vp3.tikz}}
        \end{subfigure}
        \caption{$\eta_2 \simeq 0.75$}
        \end{subfigure}
        
        \begin{subfigure}{\textwidth}
        \begin{subfigure}{0.49\textwidth}
            \scalebox{0.42}{\input{eps_vp2_error.tikz}}
        \end{subfigure}
        \hfill
        \begin{subfigure}{0.49\textwidth}
            \scalebox{0.42}{\input{eps_vp3_error.tikz}}
        \end{subfigure}
        \caption{Normalized error as a function of $\eta$. The colored vertical lines represent the values of $\eta_0$ (orange), $\eta_1$ (green) and $\eta_2$ (purple).}
        \end{subfigure}
        
        \caption{Results for the $\eps$-graph. Left: second eigenvector. Right: third eigenvector.}
        \label{fig:gsp1}
        \end{figure}
    
        \begin{figure}[]
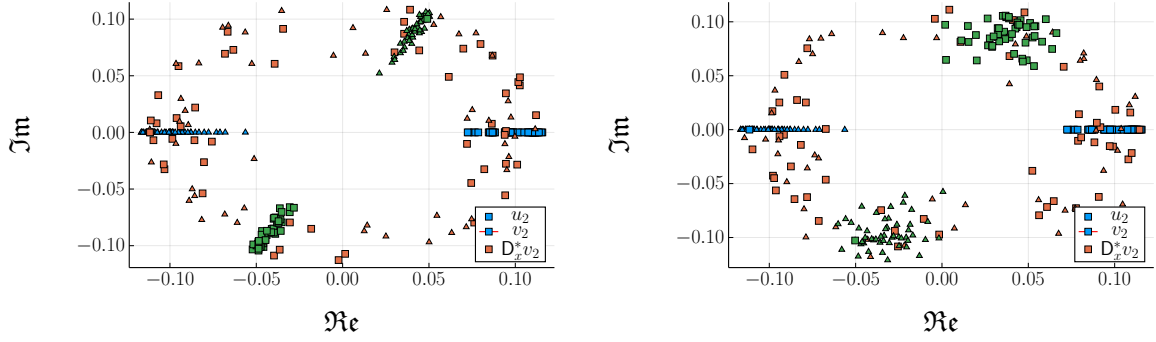

        \centering
            \begin{subfigure}{0.49\textwidth}
                \scalebox{0.42}{\input{sbm_j12_vp2.tikz}}
            \end{subfigure}
            \hfill
            \begin{subfigure}{0.49\textwidth}
                \scalebox{0.42}{\input{sbm_j75_vp2.tikz}}
            \end{subfigure}
            \caption{Second eigenvector for the SBM-graph. The two communities are depicted as triangles and squares respectively. Left: $\eta_1 \simeq 0.11$. Right: $\eta_2 \simeq 0.75$.}
            \label{fig:gsp2}
        \end{figure}

        \noindent Finally, we provide similar illustrations for a cyclic graph on size $n = 100$, in a more controlled setup: the values of $\eps_{i,j}$ are no longer randomly sampled, and we instead set $\eps_e = 1$ for all $e$ along a \emph{coherent} orientation. This way, $\eta_0$ corresponds exactly to the weak-inconsistency threshold in Condition~\ref{cond:weak_inconsistency}. We plot the corresponding normalized errors in Figure~\ref{fig:gsp3}.
        
        \begin{figure}[]
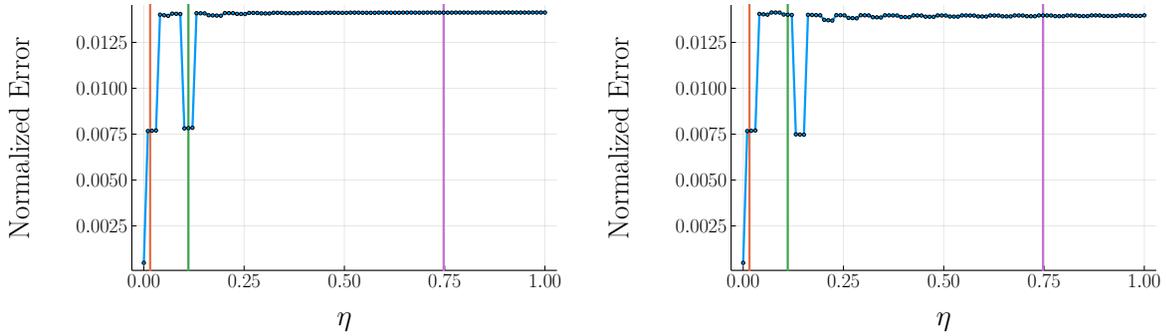

            \begin{subfigure}{0.49\textwidth}
                \scalebox{0.42}{\input{cycle_vp2_error.tikz}}
            \end{subfigure}
            \hfill
            \begin{subfigure}{0.49\textwidth}
                \scalebox{0.42}{\input{cycle_vp3_error.tikz}}
            \end{subfigure}
            \caption{Normalized errors for the cycle-graph. Left: second eigenvector. Right: third eigenvector. Vertical lines represent $\eta_0$ (orange), $\eta_1$ (green) and $\eta_2$ respectively (purple).}
            \label{fig:gsp3}
    \end{figure}
    \end{paragraph}
    
    \begin{paragraph}{Discussion}
        We observe in Figure~\ref{fig:gsp1} that $\mathsf{D}_x^* v_i$ provides a good approximation of $u_i$ for low levels of noise ($\eta = \eta_1$), but this is no longer the case at higher noise-levels ($\eta = \eta_2$): $\mathsf{D}_x^* v_2$ somewhat recovers the shape of $u_2$, but $\mathsf{D}_x ^* v_3$ appears completely unrelated to $u_3$. Results for the SBM (Figure~\ref{fig:gsp2}) suggest that major structural properties of the graph (here, the community structure) are still captured even for high noise: for both $\eta = \eta_1$ and $\eta = \eta_2$, $\mathsf{D}_x^* v_2$ clearly partition the graph into the same communities as $u_i$.
        \\
        These behaviors are reflected in the normalized-error plots (Figure~\ref{fig:gsp1}), with greater error associated to higher values of $\eta$, and smooth decay. Results on the cycle-graph (Figure~\ref{fig:gsp3}) are very different, with three different behaviors (similar for both eigenvectors), occurring at sharp thresholds: perfect correspondence with $u_i$ at $\eta = 0$, two small, low-error, \emph{plateau} (one of them occurring as soon as $\eta > 0$, and containing $\eta = \eta_0$), and a higher, essentially constant, error otherwise. For the SBM (not shown), the error decayed linearly with $\eta$.
    \end{paragraph}
    
    \vspace{0.2cm}
    
    These illustrations provide weak evidence that eigenvectors of $\mathsf{L}$ and $\mathsf{L}_\theta$ behave similarly in low-noise regimes. In this spirit, simple instances (\emph{e.g.} $\omega_v = 0$ for all $v$ but non-trivial noise) may be amenable to perturbative analyses, with $\mathsf{L}_\theta$ understood as an analytic perturbation of $\mathsf{L}$ (this is similar to the setting mentioned in Remark~\ref{rem:gamma}). 
    The situation is much more nuanced for strongly-incoherent connections, and likely requires more involved mathematical descriptions.
    
    \section{Power iteration: $\mathsf{A}_\theta$ v.s. $q(\mathsf{L}_\theta + q \mathsf{I})^{-1}$}
    \label{sect:A_vs_L}

    \textcolor{col}{We compare in the following the performance of our synchronization estimators (Section~\ref{sect:synchro}) with that of the power iteration applied to the matrix $\mathsf{A}_\theta$ (as proposed in~\cite{singer2011angular}).}

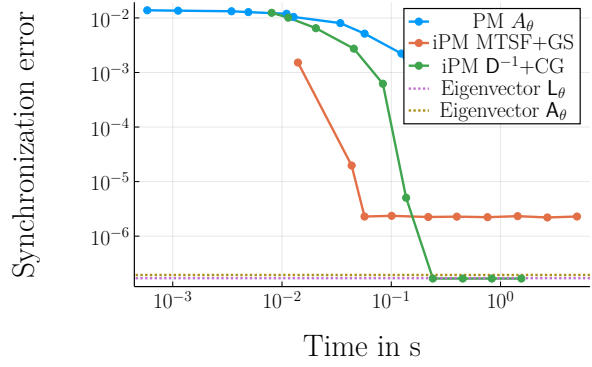
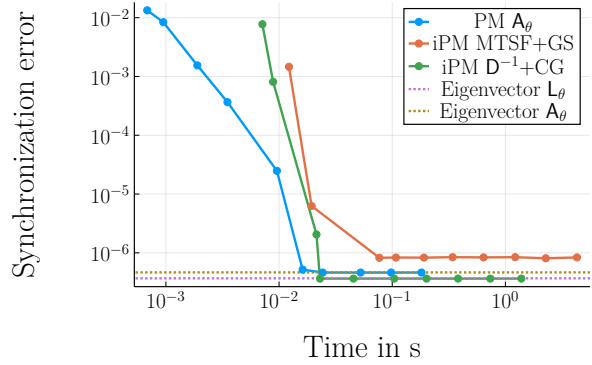
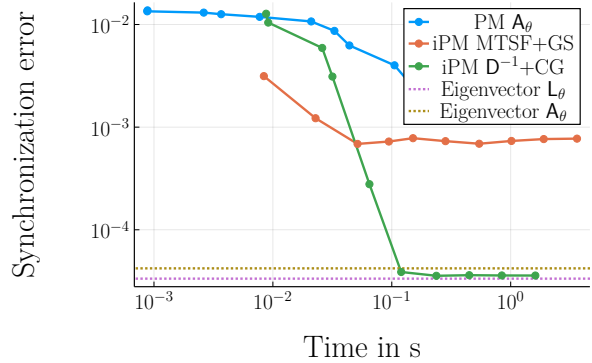
\begin{figure}[h]
        \centering
        \begin{subfigure}{\textwidth}
            \centering
            \scalebox{0.45}{

\begin{tikzpicture}[/tikz/background rectangle/.style={fill={rgb,1:red,1.0;green,1.0;blue,1.0}, fill opacity={1.0}, draw opacity={1.0}}, show background rectangle]
\begin{axis}[point meta max={nan}, point meta min={nan}, legend cell align={left}, legend columns={1}, title={}, title style={at={{(0.5,1)}}, anchor={south}, font={{\fontsize{14 pt}{18.2 pt}\selectfont}}, color={rgb,1:red,0.0;green,0.0;blue,0.0}, draw opacity={1.0}, rotate={0.0}, align={center}}, legend style={color={rgb,1:red,0.0;green,0.0;blue,0.0}, draw opacity={1.0}, line width={1}, solid, fill={rgb,1:red,1.0;green,1.0;blue,1.0}, fill opacity={1.0}, text opacity={1.0}, font={{\fontsize{18 pt}{10.4 pt}\selectfont}}, text={rgb,1:red,0.0;green,0.0;blue,0.0}, cells={anchor={center}}, at={(0.98, 0.98)}, anchor={north east}}, axis background/.style={fill={rgb,1:red,1.0;green,1.0;blue,1.0}, opacity={1.0}}, anchor={north west}, xshift={1.0mm}, yshift={-1.0mm}, width={150.4mm}, height={99.6mm}, scaled x ticks={false}, xlabel={Time in s}, x tick style={color={rgb,1:red,0.0;green,0.0;blue,0.0}, opacity={1.0}}, x tick label style={color={rgb,1:red,0.0;green,0.0;blue,0.0}, opacity={1.0}, rotate={0}}, xlabel style={at={(ticklabel cs:0.5)}, anchor=near ticklabel, at={{(ticklabel cs:0.5)}}, anchor={near ticklabel}, font={{\fontsize{26 pt}{14.3 pt}\selectfont}}, color={rgb,1:red,0.0;green,0.0;blue,0.0}, draw opacity={1.0}, rotate={0.0}}, xmode={log}, log basis x={10}, xmajorgrids={true}, xmin={0.0004451955524441048}, xmax={6.572944792945349}, xticklabels={{$10^{-3}$,$10^{-2}$,$10^{-1}$,$10^{0}$}}, xtick={{0.001,0.01,0.1,1.0}}, xtick align={inside}, xticklabel style={font={{\fontsize{18 pt}{10.4 pt}\selectfont}}, color={rgb,1:red,0.0;green,0.0;blue,0.0}, draw opacity={1.0}, rotate={0.0}}, x grid style={color={rgb,1:red,0.0;green,0.0;blue,0.0}, draw opacity={0.1}, line width={0.5}, solid}, axis x line*={left}, x axis line style={color={rgb,1:red,0.0;green,0.0;blue,0.0}, draw opacity={1.0}, line width={1}, solid}, scaled y ticks={false}, ylabel={Synchronization error}, y tick style={color={rgb,1:red,0.0;green,0.0;blue,0.0}, opacity={1.0}}, y tick label style={color={rgb,1:red,0.0;green,0.0;blue,0.0}, opacity={1.0}, rotate={0}}, ylabel style={at={(ticklabel cs:0.5)}, anchor=near ticklabel, at={{(ticklabel cs:0.5)}}, anchor={near ticklabel}, font={{\fontsize{26 pt}{14.3 pt}\selectfont}}, color={rgb,1:red,0.0;green,0.0;blue,0.0}, draw opacity={1.0}, rotate={0.0}}, ymode={log}, log basis y={10}, ymajorgrids={true}, ymin={1.1890597477676506e-7}, ymax={0.019398701832023676}, yticklabels={{$10^{-6}$,$10^{-5}$,$10^{-4}$,$10^{-3}$,$10^{-2}$}}, ytick={{1.0e-6,1.0e-5,0.0001,0.001,0.01}}, ytick align={inside}, yticklabel style={font={{\fontsize{18 pt}{10.4 pt}\selectfont}}, color={rgb,1:red,0.0;green,0.0;blue,0.0}, draw opacity={1.0}, rotate={0.0}}, y grid style={color={rgb,1:red,0.0;green,0.0;blue,0.0}, draw opacity={0.1}, line width={0.5}, solid}, axis y line*={left}, y axis line style={color={rgb,1:red,0.0;green,0.0;blue,0.0}, draw opacity={1.0}, line width={1}, solid}, colorbar={false}]
    \addplot[color={rgb,1:red,0.0;green,0.6056;blue,0.9787}, name path={39b1d956-bddb-47ec-b155-c4916750e102}, draw opacity={1.0}, line width={2}, solid, mark={*}, mark size={3.0 pt}, mark repeat={1}, mark options={color={rgb,1:red,0.0;green,0.6056;blue,0.9787}, draw opacity={1.0}, fill={rgb,1:red,0.0;green,0.6056;blue,0.9787}, fill opacity={1.0}, line width={0.75}, rotate={0}, solid}]
        table[row sep={\\}]
        {
            \\
            0.0005841808319999999  0.013811699799985179  \\
            0.001119001426  0.013585005520758236  \\
            0.003449721252  0.01325155440993647  \\
            0.004918690322  0.012782396945729349  \\
            0.011006593954000003  0.011917698469847743  \\
            0.012750170967999998  0.010436445679375907  \\
            0.034162859  0.008032591455907736  \\
            0.056956712952  0.0051520550632098506  \\
            0.12378156240400001  0.0022236324668804906  \\
            0.210650504816  0.00030527914541205695  \\
        }
        ;
    \addlegendentry {PM $A_\theta$}
    \addplot[color={rgb,1:red,0.8889;green,0.4356;blue,0.2781}, name path={2916e5cb-2521-405c-8b25-a7009759c30f}, draw opacity={1.0}, line width={2}, solid, mark={*}, mark size={3.0 pt}, mark repeat={1}, mark options={color={rgb,1:red,0.8889;green,0.4356;blue,0.2781}, draw opacity={1.0}, fill={rgb,1:red,0.8889;green,0.4356;blue,0.2781}, fill opacity={1.0}, line width={0.75}, rotate={0}, solid}]
        table[row sep={\\}]
        {
            \\
            0.013986999814  0.0015215276983217947  \\
            0.04336160880400001  1.9722501175265034e-5  \\
            0.056912525623999996  2.2895715380746104e-6  \\
            0.10080563872  2.3573414465627387e-6  \\
            0.21790647754399997  2.242817198801299e-6  \\
            0.39799530254000004  2.2782155758249994e-6  \\
            0.7604354444979999  2.227404889155444e-6  \\
            1.4341331746859998  2.325012243377969e-6  \\
            2.681973575808  2.2004708793286232e-6  \\
            5.009143792448  2.289568521330181e-6  \\
        }
        ;
    \addlegendentry {iPM MTSF+GS}
    \addplot[color={rgb,1:red,0.2422;green,0.6433;blue,0.3044}, name path={3a26cced-aa12-4b7a-a8ec-0258cd5959e7}, draw opacity={1.0}, line width={2}, solid, mark={*}, mark size={3.0 pt}, mark repeat={1}, mark options={color={rgb,1:red,0.2422;green,0.6433;blue,0.3044}, draw opacity={1.0}, fill={rgb,1:red,0.2422;green,0.6433;blue,0.3044}, fill opacity={1.0}, line width={0.75}, rotate={0}, solid}]
        table[row sep={\\}]
        {
            \\
            0.008035773582  0.012418175309561328  \\
            0.01139447855  0.01012130151456096  \\
            0.020399859872000003  0.006489853303253338  \\
            0.0454268197  0.0027352540959325977  \\
            0.08375145960800001  0.0006221393461554647  \\
            0.136734392418  5.053563620170085e-6  \\
            0.24020164487399998  1.6732017253290396e-7  \\
            0.451519064736  1.6700965263436327e-7  \\
            0.8318211904720001  1.6700490049335332e-7  \\
            1.551052027666  1.670064198532511e-7  \\
        }
        ;
    \addlegendentry {iPM $\mathsf{D}^{-1}$+CG}
    \addplot[color={rgb,1:red,0.7644;green,0.4441;blue,0.8243}, name path={9a07a0c9-1110-4172-be63-3e819901264e}, draw opacity={1.0}, line width={2}, dotted]
        table[row sep={\\}]
        {
            \\
            3.015377219183949e-8  1.6941830300859323e-7  \\
            97044.10345952789  1.6941830300859323e-7  \\
        }
        ;
    \addlegendentry {Eigenvector $\mathsf{L}_\theta$}
    \addplot[color={rgb,1:red,0.6755;green,0.5557;blue,0.0942}, name path={e7231d12-a968-4d2a-b78b-6e090a10e28b}, draw opacity={1.0}, line width={2}, dotted]
        table[row sep={\\}]
        {
            \\
            3.015377219183949e-8  1.941183397449594e-7  \\
            97044.10345952789  1.941183397449594e-7  \\
        }
        ;
    \addlegendentry {Eigenvector $\mathsf{A}_\theta$}
\end{axis}
\end{tikzpicture}}
            \caption{\textcolor{col}{$\eps$-graph, $\eta = \frac{\pi}{2n}$}.}
        \end{subfigure}
        
        \begin{subfigure}{\textwidth}
            \centering
            \scalebox{0.45}{

\begin{tikzpicture}[/tikz/background rectangle/.style={fill={rgb,1:red,1.0;green,1.0;blue,1.0}, fill opacity={1.0}, draw opacity={1.0}}, show background rectangle]
\begin{axis}[point meta max={nan}, point meta min={nan}, legend cell align={left}, legend columns={1}, title={}, title style={at={{(0.5,1)}}, anchor={south}, font={{\fontsize{14 pt}{18.2 pt}\selectfont}}, color={rgb,1:red,0.0;green,0.0;blue,0.0}, draw opacity={1.0}, rotate={0.0}, align={center}}, legend style={color={rgb,1:red,0.0;green,0.0;blue,0.0}, draw opacity={1.0}, line width={1}, solid, fill={rgb,1:red,1.0;green,1.0;blue,1.0}, fill opacity={1.0}, text opacity={1.0}, font={{\fontsize{18 pt}{10.4 pt}\selectfont}}, text={rgb,1:red,0.0;green,0.0;blue,0.0}, cells={anchor={center}}, at={(0.98, 0.98)}, anchor={north east}}, axis background/.style={fill={rgb,1:red,1.0;green,1.0;blue,1.0}, opacity={1.0}}, anchor={north west}, xshift={1.0mm}, yshift={-1.0mm}, width={150.4mm}, height={99.6mm}, scaled x ticks={false}, xlabel={Time in s}, x tick style={color={rgb,1:red,0.0;green,0.0;blue,0.0}, opacity={1.0}}, x tick label style={color={rgb,1:red,0.0;green,0.0;blue,0.0}, opacity={1.0}, rotate={0}}, xlabel style={at={(ticklabel cs:0.5)}, anchor=near ticklabel, at={{(ticklabel cs:0.5)}}, anchor={near ticklabel}, font={{\fontsize{26 pt}{14.3 pt}\selectfont}}, color={rgb,1:red,0.0;green,0.0;blue,0.0}, draw opacity={1.0}, rotate={0.0}}, xmode={log}, log basis x={10}, xmajorgrids={true}, xmin={0.0005261183941079966}, xmax={5.581072368566106}, xticklabels={{$10^{-3}$,$10^{-2}$,$10^{-1}$,$10^{0}$}}, xtick={{0.001,0.01,0.1,1.0}}, xtick align={inside}, xticklabel style={font={{\fontsize{18 pt}{10.4 pt}\selectfont}}, color={rgb,1:red,0.0;green,0.0;blue,0.0}, draw opacity={1.0}, rotate={0.0}}, x grid style={color={rgb,1:red,0.0;green,0.0;blue,0.0}, draw opacity={0.1}, line width={0.5}, solid}, axis x line*={left}, x axis line style={color={rgb,1:red,0.0;green,0.0;blue,0.0}, draw opacity={1.0}, line width={1}, solid}, scaled y ticks={false}, ylabel={Synchronization error}, y tick style={color={rgb,1:red,0.0;green,0.0;blue,0.0}, opacity={1.0}}, y tick label style={color={rgb,1:red,0.0;green,0.0;blue,0.0}, opacity={1.0}, rotate={0}}, ylabel style={at={(ticklabel cs:0.5)}, anchor=near ticklabel, at={{(ticklabel cs:0.5)}}, anchor={near ticklabel}, font={{\fontsize{26 pt}{14.3 pt}\selectfont}}, color={rgb,1:red,0.0;green,0.0;blue,0.0}, draw opacity={1.0}, rotate={0.0}}, ymode={log}, log basis y={10}, ymajorgrids={true}, ymin={2.6602981511449e-7}, ymax={0.01827411315299823}, yticklabels={{$10^{-6}$,$10^{-5}$,$10^{-4}$,$10^{-3}$,$10^{-2}$}}, ytick={{1.0e-6,1.0e-5,0.0001,0.001,0.01}}, ytick align={inside}, yticklabel style={font={{\fontsize{18 pt}{10.4 pt}\selectfont}}, color={rgb,1:red,0.0;green,0.0;blue,0.0}, draw opacity={1.0}, rotate={0.0}}, y grid style={color={rgb,1:red,0.0;green,0.0;blue,0.0}, draw opacity={0.1}, line width={0.5}, solid}, axis y line*={left}, y axis line style={color={rgb,1:red,0.0;green,0.0;blue,0.0}, draw opacity={1.0}, line width={1}, solid}, colorbar={false}]
    \addplot[color={rgb,1:red,0.0;green,0.6056;blue,0.9787}, name path={b2ea52a5-246f-4195-bacc-7974160e0028}, draw opacity={1.0}, line width={2}, solid, mark={*}, mark size={3.0 pt}, mark repeat={1}, mark options={color={rgb,1:red,0.0;green,0.6056;blue,0.9787}, draw opacity={1.0}, fill={rgb,1:red,0.0;green,0.6056;blue,0.9787}, fill opacity={1.0}, line width={0.75}, rotate={0}, solid}]
        table[row sep={\\}]
        {
            \\
            0.000683937676  0.013333456358383606  \\
            0.0009482622279999999  0.008376690793220365  \\
            0.001896216926  0.0015429171571035697  \\
            0.0035007607180000005  0.00036587310052986504  \\
            0.009569365488  2.482480211978611e-5  \\
            0.01615524208  5.182390990916096e-7  \\
            0.024301767808000002  4.622877858679628e-7  \\
            0.052601759102000004  4.622878200557793e-7  \\
            0.097706467964  4.622878232746044e-7  \\
            0.18134292787  4.622878184401514e-7  \\
        }
        ;
    \addlegendentry {PM $\mathsf{A}_\theta$}
    \addplot[color={rgb,1:red,0.8889;green,0.4356;blue,0.2781}, name path={dcad3a6e-7f6b-4102-92a5-401f0567ecae}, draw opacity={1.0}, line width={2}, solid, mark={*}, mark size={3.0 pt}, mark repeat={1}, mark options={color={rgb,1:red,0.8889;green,0.4356;blue,0.2781}, draw opacity={1.0}, fill={rgb,1:red,0.8889;green,0.4356;blue,0.2781}, fill opacity={1.0}, line width={0.75}, rotate={0}, solid}]
        table[row sep={\\}]
        {
            \\
            0.012228356232  0.0014601824769562527  \\
            0.019399729868  6.266359185566196e-6  \\
            0.07692768920399999  8.238534296562751e-7  \\
            0.107781970822  8.317879010936095e-7  \\
            0.19080237958399998  8.294846213269638e-7  \\
            0.34189058256  8.431836121656734e-7  \\
            0.64136163272  8.36539161919972e-7  \\
            1.221004951276  8.432002700762088e-7  \\
            2.2756911913780002  8.087705345625025e-7  \\
            4.293234508038  8.357155453763904e-7  \\
        }
        ;
    \addlegendentry {iPM MTSF+GS}
    \addplot[color={rgb,1:red,0.2422;green,0.6433;blue,0.3044}, name path={12b4f842-5f3f-40c2-a876-08014a9c1d8b}, draw opacity={1.0}, line width={2}, solid, mark={*}, mark size={3.0 pt}, mark repeat={1}, mark options={color={rgb,1:red,0.2422;green,0.6433;blue,0.3044}, draw opacity={1.0}, fill={rgb,1:red,0.2422;green,0.6433;blue,0.3044}, fill opacity={1.0}, line width={0.75}, rotate={0}, solid}]
        table[row sep={\\}]
        {
            \\
            0.007120954164000001  0.007733610753607895  \\
            0.008858158226  0.0008110643320500581  \\
            0.021410297765999998  2.0546456658700436e-6  \\
            0.022980126996000005  3.6461520963959334e-7  \\
            0.045399486182  3.646076594232302e-7  \\
            0.10416680424000001  3.646153646549009e-7  \\
            0.20116290549000002  3.646065564916285e-7  \\
            0.38171288291399996  3.6460863773880925e-7  \\
            0.735150622748  3.646099380215362e-7  \\
            1.382698195006  3.646060565846208e-7  \\
        }
        ;
    \addlegendentry {iPM $\mathsf{D}^{-1}$+CG}
    \addplot[color={rgb,1:red,0.7644;green,0.4441;blue,0.8243}, name path={57e6e554-aeaf-4273-9837-4f81d8753e11}, draw opacity={1.0}, line width={2}, dotted]
        table[row sep={\\}]
        {
            \\
            4.959630449835808e-8  3.687433275713261e-7  \\
            59204.1052584416  3.687433275713261e-7  \\
        }
        ;
    \addlegendentry {Eigenvector $\mathsf{L}_\theta$}
    \addplot[color={rgb,1:red,0.6755;green,0.5557;blue,0.0942}, name path={1839d9c9-71e5-4d8b-be87-162acd2db11c}, draw opacity={1.0}, line width={2}, dotted]
        table[row sep={\\}]
        {
            \\
            4.959630449835808e-8  4.6460139159756505e-7  \\
            59204.1052584416  4.6460139159756505e-7  \\
        }
        ;
    \addlegendentry {Eigenvector $\mathsf{A}_\theta$}
\end{axis}
\end{tikzpicture}}
            \caption{\textcolor{col}{DC-SBM 1, $\eta = \frac{\pi}{2n}$}.}
        \end{subfigure}
        
        \begin{subfigure}{\textwidth}
            \centering
            \scalebox{0.45}{

\begin{tikzpicture}[/tikz/background rectangle/.style={fill={rgb,1:red,1.0;green,1.0;blue,1.0}, fill opacity={1.0}, draw opacity={1.0}}, show background rectangle]
\begin{axis}[point meta max={nan}, point meta min={nan}, legend cell align={left}, legend columns={1}, title={}, title style={at={{(0.5,1)}}, anchor={south}, font={{\fontsize{14 pt}{18.2 pt}\selectfont}}, color={rgb,1:red,0.0;green,0.0;blue,0.0}, draw opacity={1.0}, rotate={0.0}, align={center}}, legend style={color={rgb,1:red,0.0;green,0.0;blue,0.0}, draw opacity={1.0}, line width={1}, solid, fill={rgb,1:red,1.0;green,1.0;blue,1.0}, fill opacity={1.0}, text opacity={1.0}, font={{\fontsize{18 pt}{10.4 pt}\selectfont}}, text={rgb,1:red,0.0;green,0.0;blue,0.0}, cells={anchor={center}}, at={(0.98, 0.98)}, anchor={north east}}, axis background/.style={fill={rgb,1:red,1.0;green,1.0;blue,1.0}, opacity={1.0}}, anchor={north west}, xshift={1.0mm}, yshift={-1.0mm}, width={150.4mm}, height={99.6mm}, scaled x ticks={false}, xlabel={Time in s}, x tick style={color={rgb,1:red,0.0;green,0.0;blue,0.0}, opacity={1.0}}, x tick label style={color={rgb,1:red,0.0;green,0.0;blue,0.0}, opacity={1.0}, rotate={0}}, xlabel style={at={(ticklabel cs:0.5)}, anchor=near ticklabel, at={{(ticklabel cs:0.5)}}, anchor={near ticklabel}, font={{\fontsize{26 pt}{14.3 pt}\selectfont}}, color={rgb,1:red,0.0;green,0.0;blue,0.0}, draw opacity={1.0}, rotate={0.0}}, xmode={log}, log basis x={10}, xmajorgrids={true}, xmin={0.0006836902113902431}, xmax={4.6374020645993355}, xticklabels={{$10^{-3}$,$10^{-2}$,$10^{-1}$,$10^{0}$}}, xtick={{0.001,0.01,0.1,1.0}}, xtick align={inside}, xticklabel style={font={{\fontsize{18 pt}{10.4 pt}\selectfont}}, color={rgb,1:red,0.0;green,0.0;blue,0.0}, draw opacity={1.0}, rotate={0.0}}, x grid style={color={rgb,1:red,0.0;green,0.0;blue,0.0}, draw opacity={0.1}, line width={0.5}, solid}, axis x line*={left}, x axis line style={color={rgb,1:red,0.0;green,0.0;blue,0.0}, draw opacity={1.0}, line width={1}, solid}, scaled y ticks={false}, ylabel={Synchronization error}, y tick style={color={rgb,1:red,0.0;green,0.0;blue,0.0}, opacity={1.0}}, y tick label style={color={rgb,1:red,0.0;green,0.0;blue,0.0}, opacity={1.0}, rotate={0}}, ylabel style={at={(ticklabel cs:0.5)}, anchor=near ticklabel, at={{(ticklabel cs:0.5)}}, anchor={near ticklabel}, font={{\fontsize{26 pt}{14.3 pt}\selectfont}}, color={rgb,1:red,0.0;green,0.0;blue,0.0}, draw opacity={1.0}, rotate={0.0}}, ymode={log}, log basis y={10}, ymajorgrids={true}, ymin={2.7890565475211367e-5}, ymax={0.016448554127034443}, yticklabels={{$10^{-4}$,$10^{-3}$,$10^{-2}$}}, ytick={{0.0001,0.001,0.01}}, ytick align={inside}, yticklabel style={font={{\fontsize{18 pt}{10.4 pt}\selectfont}}, color={rgb,1:red,0.0;green,0.0;blue,0.0}, draw opacity={1.0}, rotate={0.0}}, y grid style={color={rgb,1:red,0.0;green,0.0;blue,0.0}, draw opacity={0.1}, line width={0.5}, solid}, axis y line*={left}, y axis line style={color={rgb,1:red,0.0;green,0.0;blue,0.0}, draw opacity={1.0}, line width={1}, solid}, colorbar={false}]
    \addplot[color={rgb,1:red,0.0;green,0.6056;blue,0.9787}, name path={5c2b59ff-a271-4cd8-aea6-78bcd4d3e331}, draw opacity={1.0}, line width={2}, solid, mark={*}, mark size={3.0 pt}, mark repeat={1}, mark options={color={rgb,1:red,0.0;green,0.6056;blue,0.9787}, draw opacity={1.0}, fill={rgb,1:red,0.0;green,0.6056;blue,0.9787}, fill opacity={1.0}, line width={0.75}, rotate={0}, solid}]
        table[row sep={\\}]
        {
            \\
            0.00088255151  0.013731327582234259  \\
            0.0008775980439999999  0.01344589213026655  \\
            0.002621059914  0.013059498346980215  \\
            0.003672605622  0.012601729806830107  \\
            0.007742558200000001  0.011871993532792741  \\
            0.020993942552  0.010713349907524266  \\
            0.032941757015999995  0.008668684555591482  \\
            0.044010551317999996  0.006253576069921389  \\
            0.105442010718  0.004003738722360301  \\
            0.203589082392  0.0016455859091362185  \\
        }
        ;
    \addlegendentry {PM $\mathsf{A}_\theta$}
    \addplot[color={rgb,1:red,0.8889;green,0.4356;blue,0.2781}, name path={56f51497-7007-480f-8506-247fd935c56d}, draw opacity={1.0}, line width={2}, solid, mark={*}, mark size={3.0 pt}, mark repeat={1}, mark options={color={rgb,1:red,0.8889;green,0.4356;blue,0.2781}, draw opacity={1.0}, fill={rgb,1:red,0.8889;green,0.4356;blue,0.2781}, fill opacity={1.0}, line width={0.75}, rotate={0}, solid}]
        table[row sep={\\}]
        {
            \\
            0.008396662632000002  0.0031420763408713023  \\
            0.022852525199999998  0.0012217368360024777  \\
            0.05150562704000001  0.0006864241938183461  \\
            0.094497149658  0.000724159506617086  \\
            0.150814577556  0.0007808140896367496  \\
            0.28294814872  0.0007304589008822358  \\
            0.543404180012  0.0006893775304331516  \\
            1.0152214597999998  0.0007336341151662049  \\
            1.89682974055  0.0007654707218518306  \\
            3.612754631262  0.0007728907960008919  \\
        }
        ;
    \addlegendentry {iPM MTSF+GS}
    \addplot[color={rgb,1:red,0.2422;green,0.6433;blue,0.3044}, name path={63a7b902-2401-4bd6-aadc-99c132f17f46}, draw opacity={1.0}, line width={2}, solid, mark={*}, mark size={3.0 pt}, mark repeat={1}, mark options={color={rgb,1:red,0.2422;green,0.6433;blue,0.3044}, draw opacity={1.0}, fill={rgb,1:red,0.2422;green,0.6433;blue,0.3044}, fill opacity={1.0}, line width={0.75}, rotate={0}, solid}]
        table[row sep={\\}]
        {
            \\
            0.008774046724  0.012733745971592447  \\
            0.009139307372  0.010445899689457241  \\
            0.025991020695999997  0.005911547334840479  \\
            0.031679057184  0.003111346385217115  \\
            0.06487023276400002  0.00027829788155590314  \\
            0.119980483958  3.878350733980825e-5  \\
            0.237084166216  3.561667182460904e-5  \\
            0.448295017338  3.610604248551396e-5  \\
            0.842316663012  3.5881218664010114e-5  \\
            1.61161006758  3.58243505448717e-5  \\
        }
        ;
    \addlegendentry {iPM $\mathsf{D}^{-1}$+CG}
    \addplot[color={rgb,1:red,0.7644;green,0.4441;blue,0.8243}, name path={35a64330-1afa-415f-b3ee-fe4f1158bb17}, draw opacity={1.0}, line width={2}, dotted]
        table[row sep={\\}]
        {
            \\
            1.0079615669279289e-7  3.340969568347928e-5  \\
            31455.032630963135  3.340969568347928e-5  \\
        }
        ;
    \addlegendentry {Eigenvector $\mathsf{L}_\theta$}
    \addplot[color={rgb,1:red,0.6755;green,0.5557;blue,0.0942}, name path={c57a90e4-16fc-4eee-8d68-cd809e809259}, draw opacity={1.0}, line width={2}, dotted]
        table[row sep={\\}]
        {
            \\
            1.0079615669279289e-7  4.207878221814456e-5  \\
            31455.032630963135  4.207878221814456e-5  \\
        }
        ;
    \addlegendentry {Eigenvector $\mathsf{A}_\theta$}
\end{axis}
\end{tikzpicture}}
            \caption{\textcolor{col}{$\eps$-graph, $\eta = \frac{\pi}{100}$}.}
        \end{subfigure}
        \caption{\textcolor{col}{Runtime-precision comparisons for the power iteration based on $\mathsf{A}_\theta$~\cite{singer2011angular} v.s. the approximate (inverse) power methods from Section~\ref{sect:synchro}. The results obtained from the exact eigenvectors of $\mathsf{L}_\theta$ and $\mathsf{A}_\theta$ are represented as horizontal lines.}}
        \label{fig:A_vs_L}
        \end{figure}

    \begin{paragraph}{Setup}
    \textcolor{col}{We perform the experiment on three different graph/connection models, where the connection is generated according to the model of Equation~\eqref{eq:degradation} for different incoherence levels.
    \begin{itemize}
        \item With low-incoherence connections (weak-inconsistency condition satisfied), with $\eta = \frac{\pi}{2n}$, on two different random-graph models.
        \begin{itemize}
            \item On the $\eps$-graph model of Section~\ref{subsect:runtime_precision_smoothing}.
            \item On the DC-SBM 1 model of Section~\ref{subsect:runtime_precision_smoothing}. 
        \end{itemize}
        Note that, for these low-incoherence models, the spectral gaps of the matrices $\mathsf{L}_\theta$ and $\mathsf{A}_\theta$ should be similar to those of the matrices $\mathsf{L}$ and $\mathsf{A}$ (which we analyzed in Remark~\ref{rem:synchro}). Let us stress that low-incoherence connections of this type \emph{are} encountered in practice, in problems such as ranking~\cite{yu_angular_2009,cucuringu_syncrank_2016} for instance.
        \item In higher-incoherence connections (no guarantee that the weak-inconsistency condition is satisfied), with $\eta = \frac{\pi}{100}$, on the same $\eps$-graph model.
    \end{itemize}}
    \noindent \textcolor{col}{For each graph we estimate the solution $x \in U(\C)^\mathcal{V}$ to the angular-synchronization problem by performing $k$ power-method iterations, for $k \in \{ 1, 2, 4, 7, 13, 24, 45, 85, 160, 300 \}$. Each power-method iteration consists in either a multiplication by $\mathsf{A}_\theta$, or in an approximate smoothing by $\hat{f}$ with $m = 5$ MTSFs (resp., by $m = 5$ iterations of diagonally-preconditioned CG). We set $q = 10^{-4} \times \overline{d}$, and measure for each $k$ the resulting synchronization error (Equation~\eqref{eq:sync_error}) and associated runtime. Each runtime measurement is averaged over $100$ runs. The final results are averaged over $10$ uniformly-selected initializations $f_0$ of the power method, themselves averaged over $5$ graph-realizations for each model.}
    \\
    \textcolor{col}{The results are depicted in Figure~\ref{fig:A_vs_L}.}
    \end{paragraph}

    \vspace{0.2cm}

    \begin{paragraph}{Comments} \textcolor{col}{There are two main observations.
        \begin{enumerate}
            \item For low-inconsistence connections, the predonimant parameter is the expansivity of the graph. In particular, the spectral gap of the matrix $q(\mathsf{L}_\theta + q \mathsf{I})^{-1}$ is much larger than that of $\mathsf{A}_\theta$ for the $\eps$-graph, and we obtain a significant gain in performance with both the $\hat{f}$-based and CG-based iterations (for $\hat{f}$, the synchronization error is better in $2$ iterations than in $300$ iterations for $\mathsf{A}_\theta$). The situation is more nuanced for the DC-SBM 1 graph: in that case, the spectral gap of $\mathsf{A}_\theta$ is quite small (due to the clear community structure exhibited for this graph model\footnote{In particular, the parameters of this model ensure that we are well-above the theoretical community-reconstruction threshold for DC-SBMs (see, \emph{e.g.}~\cite{gulikers2016non}).}), so that we could expect a similar speed-up, but not actually small enough for the power-method based on $\mathsf{A}_\theta$ to converge very slowly (as we observe in Figure~\ref{fig:A_vs_L})\footnote{For the DC-SBM 2 model, the community structure is even more clear-cut, and we do observe a similar speed-up.}.
            \item For higher-incoherence connections, the bottom eigenvalue $\lambda_1$ of $\mathsf{L}_\theta$ (resp. the top eigenvalue $\alpha_1$ of $\mathsf{A}_\theta$) gets larger (resp. smaller), which in this example appears to slow-down (resp. speed-up) the convergence of the power-method based on the matrix $q (\mathsf{L}_\theta + q \mathsf{I})^{-1}$ (resp. $\mathsf{A}_\theta$). This translates to a smaller improvement in the performance of our method against the power-method based on $\mathsf{A}_\theta$. Nevertheless, it is hard to derive a general conclusion: the eigenvalues $\lambda_2$ and $\alpha_2$ are also affected by the incoherence of the connection, and the convergence rates of the power-methods depend on the \emph{rates at which the ratios $\frac{\lambda_2 + q}{\lambda_1 + q}$ and $\frac{\alpha_1}{\alpha_2}$ evolve} (when the incoherence is low, this problem should be amenable to a pertubative analysis, as sketched in Section~\ref{supmat:tikhonov_numerical}).
        \end{enumerate}
        Let us note that the choice of $q$ (which has no bearing on the end-result of the iterations, but does impact the runtime of the inverse power methods) is still subject to the trade-off discussed in Section~\ref{sect:synchro}: $q$ should be small enough to allow for a fast convergence of the power iteration, but large enough so that each iteration can be performed efficiently. Here, we purposefully choose a small value of $q$ to showcase the gain in performance of our methods against the power-iteration based on $\mathsf{A}_\theta$, but do not tailor this choice to the exact graph nor to the incoherence level.}
    \end{paragraph}
    
    \section{Illustration on MNIST images: an average digit}
    \label{sect:mnist}
    
    \textcolor{col}{We propose a variant of the illustration from Section~\ref{sect:illustration} in which the incoherence of the connection is not due to the addition of (controlled) noise, but is inherent to the data itself. Our goal is to extract from a set of similar images an ``average representative'' of this dataset. Here, we consider a set of $5000$ ``1'''s in the MNIST dataset~\cite{lecun1998gradient} (a dataset of $28 \times 28$ handwritten digits), that all represent similar images but, due to different calligraphic styles, present some discrepancies: in particular, the images are not aligned properly, and two very similar digits may differ by a simple rotation.}

    \textcolor{col}{For each pair of images $(I,J)$ in the dataset, we compute a rotation $r_{I,J}$ registering the two images (based on image moments), extract for each image $I$ its set of $s = 35$ nearest neighbors (with respect to the distance $\Vert I - r_{J,I}(J) \Vert_2$ in $\R^{28 \times 28}$), and build the corresponding nearest-neighbor graph and connection.}
    
    \textcolor{col}{We then compute a global alignment via angular synchronization, using respectively the exact eigenvector of $\mathsf{L}_\theta$ and a MTSF-based estimate (obtained from $\hat{f}$, with $k = 10$ power-method iterations, $m = 5$ MTSFs sampled at each step, and $q = 10^{-4} \times \overline{d}$), and average the $5000$ aligned images. For comparison, we also include the average obtained without the rotations.}

\begin{figure}[ht]
        \centering
		\begin{subfigure}{\linewidth}
        \begin{subfigure}{0.3\textwidth}
            \centering
            \includegraphics[scale=3]{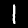}
        \end{subfigure}
        \hfill
        \begin{subfigure}{0.3\textwidth}
            \centering
            \includegraphics[scale=3]{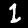}
        \end{subfigure}
        \hfill
        \begin{subfigure}{0.3\textwidth}
            \centering
            \includegraphics[scale=3]{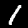}
        \end{subfigure}
		\caption{\textcolor{col}{Some ``1'''s from the MNIST dataset.}}
		\end{subfigure}

		\vspace{0.5cm}

        \begin{subfigure}{0.3\textwidth}
            \centering
            \includegraphics[scale=3]{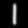}
            \caption{\textcolor{col}{Exact eigenvector.}}
        \end{subfigure}
        \hfill
        \begin{subfigure}{0.3\textwidth}
            \centering
            \includegraphics[scale=3]{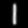}
            \caption{\textcolor{col}{$\overline{f}$-based estimation.}}
        \end{subfigure}
        \hfill
        \begin{subfigure}{0.3\textwidth}
            \centering
            \includegraphics[scale=3]{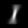}
            \caption{\textcolor{col}{No rotation.}}
        \end{subfigure}
        \caption{\textcolor{col}{Top: some ``1'''s from the MNIST dataset. Bottom: average ``1'' in the dataset, as obtained by three different methods (left: exact synchronization; center: MTSF-based synchronization; right: no synchronization). The results of the synchronization are only given up to a global rotation, which is apparent in the recovered images.}}
        \label{fig:mnist}
    \end{figure}

    \textcolor{col}{The results are depicted in Figure~\ref{fig:mnist}.}
    
    \section{Towards Extensions: an Alternative Construction}
    \label{sect:extension}
    
    We exhibit how the loop-erasure procedure in Algorithm~\ref{alg:sampling_mtsf} can be understood as stemming from the Feynman-Kac formula in Proposition~\ref{prop:feynmac_kac}.
    The aim is twofold: first, to better understand the relation between the Feynman-Kac formula in Proposition~\ref{prop:feynmac_kac} and the MTSF-based estimaton of Theorem~\ref{th:fermionic} and, second, to propose a possible proof outline for generalizations of our work. 
    
    \vspace{0.2cm}
    
    \begin{paragraph}{Loop Erasure from Feynman-Kac}
    \color{col}{Let us begin by defining an equivalence relation on paths in $P_i^\Gamma$ (recall that $P_i^\Gamma$ denotes the set of paths from $i$ to $\Gamma$ in the extended graph $\mathcal{G}^\Gamma$): two paths $p,q \in P_i^\Gamma$ are equivalent, which we denote by $p \simeq q$, if the only difference between the two is the orientation of (at most) \emph{one} cycle; that is, $p \simeq q$ if $p = aCb$ and $q = a C^* b$, with $C$ a (possibly empty) cycle in $\mathcal{G}$.
        For convenience, we denote by $p(C)$ and $p(C^*)$ the two possible representatives of each of the induced equivalence classes, and by $\overline{p(C)}$ the equivalence class itself.
        \\
        For any node $i \in \mathcal{V}$, we can write the expected value obtained at $i$ by performing parallel transport, from which we obtain:
        \begin{align}
            \label{eq:obs}
            \sum_{p \in P_i^\Gamma} \prb_{\nu_i^\Gamma}(p) \psi_{p_\Gamma^*} & = \sum_{\overline{p(C)} \in \faktor{P_i^\Gamma}{\simeq}} \left( \prb_{\nu_i^\Gamma}(p(C)) \psi_{p(C)^*_\Gamma} + \prb_{\nu_i^\Gamma}(p(C^*)) \psi_{p(C^*)^*_\Gamma} \right)  \\
            & = \sum_{\overline{p(C)} \in \faktor{P_i^\Gamma} {\simeq}} \prb_{\nu_i^\Gamma}(\overline{p(C)}) \left( \psi_{a^*} \circ (\psi_C + \psi_{C^*}) \circ \psi_{b^*} \right)\label{eq:point} \\
            & = \sum_{\overline{p(C)} \in \faktor{P_i^\Gamma} {\simeq}} 2 \prb_{\nu_i^\Gamma}(\overline{p(C)}) \left( \left( \psi_{a^*} \circ \left(\frac{\psi_C + \psi_{C^*}}{2} \right) \right) \circ \psi_{b^*} \right), \label{eq:point2}
        \end{align}
        where $p(C) = aCb$, and we denote by $\prb_{\nu_i^\Gamma}(\overline{p(C)}) = \prb_{\nu_i^\Gamma}(p(C))$ in Equation~\ref{eq:point} the probability of sampling either of the two representatives $p(C)$ or $p(C^*)$ (note that $\prb_{\nu_i^\Gamma}(p_C)~=~\prb_{\nu_i^\Gamma}(p_{C^*})$).
        \\
        Let us re-phrase this observation in a more qualitative language: whereas each estimation in the Feynman-Kac estimator proceeds by parallel transport along a path $p(C)$, with the transport application \emph{depending on the orientation of the cycle $C$}, we observe from Equation~\eqref{eq:point2} that we can instead perform an \emph{orientation-agnostic estimation} resulting in the same expectation, by applying the map $\frac{\psi_C + \psi_{C^*}}{2}$ instead of the connection map $\psi_C$; the factor $2$ in Equation~\eqref{eq:point2} in turn accounts for the possibility of sampling either $p(C)$ or $p(C^*)$.
        \\
        \emph{Under the weak-inconsistency condition}, this allows to probabilistically interpret the maps $\frac{\psi_C + \psi_{C^*}}{2}$ in Equation~\eqref{eq:point2}, which act as multiplications by the factors $\cos(\theta_C)$:
        \begin{equation}
            \left( \frac{\psi_C + \psi_{C^*}}{2} \right)(z) = \cos(\theta_C) \cdot z.
        \end{equation}
        We can thus understand Equation~\eqref{eq:point2} as a variant of the Feynman-Kac estimator involving a Bernoulli trial (with success probability $\cos(\theta_C)$), going as follows. Suppose that the path $p$ contains a unique cycle $C$, and takes the form $p = a C b \ (j,\Gamma)$, then:
        \begin{enumerate}
            \item if the trial succeeds, set $\psi_{(ab)^*} (g(j))$
            as the estimation on node $i$.
            \item If the trial fails, set $0$ as the estimation.
        \end{enumerate}
        Let us stress that, in expectation, this procedure still results in the smoothed signal $f_*$.
        \\
    This argument extends to any number $k$ of cycles, and translates to a significant reduction of variance for the Feynman-Kac estimator\footnote{Note that this requires to be careful in the definition of the successive cycles, as one cycle may be contained within another. Here, a convenient choice of definition is to define the set of cycles according to the following procedure: whenever a cycle is observed in the path $p$ constructed in the random walk, add it to the set of cycles and remove it from $p$.}: any path $p$ with $k$ cycles, sampled according to $\nu_i^\Gamma$, belongs to an equivalence class of $2^k$ paths, for each of which the Feynman-Kac estimator results in a different value at node $i$, while the modified estimator assigns the same value to all $2^k$ paths.
    \\
    \emph{One can go a step further, and get rid of the weak-inconsistency condition}. Indeed, for any cycle $C$ such that $\cos(\theta_C) < 0$, one can perform a Bernoulli trial with success probability $-\cos(\theta_C) = (-1) \times \cos(\theta_C)$, and multiply the resulting estimate $-1$ each time such a cycle is constructed in the path $p$, thereby \emph{extending the estimation to any $U(\C)$-connection}.
    \\
    The crucial observation follows, as one can identify:
    \begin{itemize}
        \item the modified Feynman-Kac estimator, and in particular the Bernoulli trial for each cycle;
        \item the cycle-acceptance trial of Algorithm~\ref{alg:sampling_mtsf_bis}.
    \end{itemize}
    More specifically, consider the first node $i$ in the queue of the MTSF-sampling algorithm, at which the first loop-erased random walk starts: this walk constructs a path $p$ and, at each encountered cycle $C$ in $\mathcal{G}$, may be interrupted due to cycle-acceptance, with probability $1 - \cos(\theta_C)$, in which case the ensuing estimation at node $i$ is $\widetilde{f}(i) = 0$.
    \\
    Similarly, consider the modified Feynman-Kac estimator, which performs successive Bernoulli trials with success probability $\cos(\theta_C)$, and notice that the result of the estimation at node $i$ is decided \emph{as soon as a trial fails}, in which case it is going to be $0$. Effectively, the resulting path $p$ upon which the parallel transport will be applied describes a part of a multi-type spanning forest: the branch on a unicycle, and the associated cycle.
    \\
    In both cases, if no cycle is accepted (resp. no Bernoulli trial fails), the walk ends upon reaching node $\Gamma$, and the estimation at node $i$ is described by the parallel transport from $\Gamma$ to $i$ along the path $p \subseteq \mathcal{E}$.}
    \end{paragraph}
    
    \vspace{0.2cm}

    \color{col}{The question is then the following. Can such a construction be extended to \emph{all} the nodes in $\mathcal{G}$, so as to recover the (unbiasedness of the) MTSF-based estimator $\widetilde{f}$? This question is important for a few reasons.
    \begin{itemize}
        \item First and foremost, this construction could both allow to extend MTSF-based estimators to non-weakly-inconsistent connections, and provide a much more intuitive construction of multi-type spanning forests.
        \item Second, it is likely to be more flexible: instead of interpreting the factors $\frac{\psi_C + \psi_{C^*}}{2}$ probabilistically, one could, \emph{e.g.}, systematically apply the map $\frac{\psi_C + \psi_{C^*}}{2}$, which would result in connection-aware MTSF-based estimators \emph{with no unicycles}.
    \end{itemize}
    
    \begin{rem}
    A similar cycle-equivalence argument has also been used to prove the determinantality of the spanning-forest-of-unicycles measure of~\cite{forman1993determinants,kenyon_spanning_2011} in~\cite{kassel2015learning} (in the weakly-inconsistent setting) and, more recently, to derive the exact expected number of steps for Algorithm~\ref{alg:sampling_mtsf_bis} algorithm in~\cite{FanBar2024}. 
    \end{rem}}

\end{document}